\documentclass[11pt]{article}
\usepackage[margin=1in]{geometry}
\usepackage{amsmath,amssymb,amsthm,mathtools,microtype,booktabs,tikz,hyperref,placeins}
\usepackage[font=small,labelfont=bf]{caption}
\usepackage{needspace}
\usetikzlibrary{arrows.meta,calc,decorations.pathreplacing,patterns,shapes.geometric}
\hypersetup{colorlinks=true,linkcolor=blue!45!black,citecolor=blue!45!black,urlcolor=blue!45!black,
 pdftitle={Cylindrical dimers and monomial immanants},pdfauthor={Pavlo Pylyavskyy}}
\newtheorem{theorem}{Theorem}[section]
\newtheorem{proposition}[theorem]{Proposition}
\newtheorem{lemma}[theorem]{Lemma}
\newtheorem{corollary}[theorem]{Corollary}
\newtheorem{conjecture}[theorem]{Conjecture}
\theoremstyle{definition}
\newtheorem{example}[theorem]{Example}
\theoremstyle{remark}

\newcommand{\NN}{\mathbb Z_{\ge0}}
\newcommand{\wt}{\operatorname{wt}}
\newcommand{\sgn}{\operatorname{sgn}}

\newcommand{\type}{\operatorname{type}}
\newcommand{\tr}{\operatorname{tr}}
\newcommand{\Imm}{\operatorname{Imm}}
\newcommand{\ch}{\operatorname{ch}}

\newcommand{\cyl}{\mathcal A}
\newcommand{\xx}{\mathbf x}
\newcommand{\yy}{\mathbf y}
\newcommand{\mm}{\mathbf m}
\definecolor{cA}{HTML}{0072B2}
\definecolor{cB}{HTML}{D55E00}
\definecolor{cC}{HTML}{009E73}
\tikzset{arr/.style={-{Stealth[length=2.0mm]},thick},base/.style={black!23,line width=.65pt},
  bv/.style={circle,fill=black,draw=black,inner sep=1.7pt},
  wv/.style={circle,fill=white,draw=black,inner sep=1.7pt},
  aux/.style={cC!65!black,line width=1.4pt}}
\newcommand{\MatchingFigure}{%
\begin{tikzpicture}[scale=1.0,baseline=(current bounding box.center)]
\fill[black!3] (0,0) circle (.25);
\draw[black!20] (0,0) circle (.25);
\draw[cB!75,dashed,line width=.65pt] (.27,0)--(1.78,0);
\coordinate (v00) at (1.4,1.4);
\coordinate (v01) at (-1.4,1.4);
\coordinate (v02) at (-1.4,-1.4);
\coordinate (v03) at (1.4,-1.4);
\coordinate (v10) at (0.6,0.6);
\coordinate (v11) at (-0.6,0.6);
\coordinate (v12) at (-0.6,-0.6);
\coordinate (v13) at (0.6,-0.6);
\draw[base] (v00)--(v03);
\draw[base] (v00)--(v01);
\draw[base] (v00)--(v10);
\draw[base] (v02)--(v01);
\draw[base] (v02)--(v03);
\draw[base] (v02)--(v12);
\draw[base] (v11)--(v10);
\draw[base] (v11)--(v12);
\draw[base] (v11)--(v01);
\draw[base] (v13)--(v12);
\draw[base] (v13)--(v10);
\draw[base] (v13)--(v03);
\draw[cA,line width=1.65pt] (v00)--(v03);
\draw[cA,line width=1.65pt] (v13)--(v12);
\draw[cA,line width=1.65pt] (v02)--(v01);
\draw[cA,line width=1.65pt] (v11)--(v10);
\node[bv] at (v00) {};
\node[wv] at (v01) {};
\node[bv] at (v02) {};
\node[wv] at (v03) {};
\node[wv] at (v10) {};
\node[bv] at (v11) {};
\node[wv] at (v12) {};
\node[bv] at (v13) {};
\node[font=\small,align=center] at (0,-2.05) {$h=0$};
\end{tikzpicture}
\hspace{.5cm}
\begin{tikzpicture}[scale=1.0,baseline=(current bounding box.center)]
\fill[black!3] (0,0) circle (.25);
\draw[black!20] (0,0) circle (.25);
\draw[cB!75,dashed,line width=.65pt] (.27,0)--(1.78,0);
\coordinate (v00) at (1.4,1.4);
\coordinate (v01) at (-1.4,1.4);
\coordinate (v02) at (-1.4,-1.4);
\coordinate (v03) at (1.4,-1.4);
\coordinate (v10) at (0.6,0.6);
\coordinate (v11) at (-0.6,0.6);
\coordinate (v12) at (-0.6,-0.6);
\coordinate (v13) at (0.6,-0.6);
\draw[base] (v00)--(v03);
\draw[base] (v00)--(v01);
\draw[base] (v00)--(v10);
\draw[base] (v02)--(v01);
\draw[base] (v02)--(v03);
\draw[base] (v02)--(v12);
\draw[base] (v11)--(v10);
\draw[base] (v11)--(v12);
\draw[base] (v11)--(v01);
\draw[base] (v13)--(v12);
\draw[base] (v13)--(v10);
\draw[base] (v13)--(v03);
\draw[cA,line width=1.65pt] (v11)--(v01);
\draw[cA,line width=1.65pt] (v00)--(v10);
\draw[cA,line width=1.65pt] (v13)--(v03);
\draw[cA,line width=1.65pt] (v02)--(v12);
\node[bv] at (v00) {};
\node[wv] at (v01) {};
\node[bv] at (v02) {};
\node[wv] at (v03) {};
\node[wv] at (v10) {};
\node[bv] at (v11) {};
\node[wv] at (v12) {};
\node[bv] at (v13) {};
\node[font=\small,align=center] at (0,-2.05) {$h=1$};
\end{tikzpicture}
\hspace{.5cm}
\begin{tikzpicture}[scale=1.0,baseline=(current bounding box.center)]
\fill[black!3] (0,0) circle (.25);
\draw[black!20] (0,0) circle (.25);
\draw[cB!75,dashed,line width=.65pt] (.27,0)--(1.78,0);
\coordinate (v00) at (1.4,1.4);
\coordinate (v01) at (-1.4,1.4);
\coordinate (v02) at (-1.4,-1.4);
\coordinate (v03) at (1.4,-1.4);
\coordinate (v10) at (0.6,0.6);
\coordinate (v11) at (-0.6,0.6);
\coordinate (v12) at (-0.6,-0.6);
\coordinate (v13) at (0.6,-0.6);
\draw[base] (v00)--(v03);
\draw[base] (v00)--(v01);
\draw[base] (v00)--(v10);
\draw[base] (v02)--(v01);
\draw[base] (v02)--(v03);
\draw[base] (v02)--(v12);
\draw[base] (v11)--(v10);
\draw[base] (v11)--(v12);
\draw[base] (v11)--(v01);
\draw[base] (v13)--(v12);
\draw[base] (v13)--(v10);
\draw[base] (v13)--(v03);
\draw[cA,line width=1.65pt] (v00)--(v01);
\draw[cA,line width=1.65pt] (v13)--(v10);
\draw[cA,line width=1.65pt] (v02)--(v03);
\draw[cA,line width=1.65pt] (v11)--(v12);
\node[bv] at (v00) {};
\node[wv] at (v01) {};
\node[bv] at (v02) {};
\node[wv] at (v03) {};
\node[wv] at (v10) {};
\node[bv] at (v11) {};
\node[wv] at (v12) {};
\node[bv] at (v13) {};
\node[font=\small,align=center] at (0,-2.05) {$h=2$};
\end{tikzpicture}%
}
\newcommand{\ColoringFigure}{%
\begin{tikzpicture}[scale=0.82,baseline=(current bounding box.center)]
\fill[black!3] (0,0) circle (.25);
\draw[black!20] (0,0) circle (.25);
\draw[cB!75,dashed,line width=.65pt] (.27,0)--(1.78,0);
\coordinate (v00) at (1.4,1.4);
\coordinate (v01) at (-1.4,1.4);
\coordinate (v02) at (-1.4,-1.4);
\coordinate (v03) at (1.4,-1.4);
\coordinate (v10) at (0.6,0.6);
\coordinate (v11) at (-0.6,0.6);
\coordinate (v12) at (-0.6,-0.6);
\coordinate (v13) at (0.6,-0.6);
\draw[base] (v00)--(v03);
\draw[base] (v00)--(v01);
\draw[base] (v00)--(v10);
\draw[base] (v02)--(v01);
\draw[base] (v02)--(v03);
\draw[base] (v02)--(v12);
\draw[base] (v11)--(v10);
\draw[base] (v11)--(v12);
\draw[base] (v11)--(v01);
\draw[base] (v13)--(v12);
\draw[base] (v13)--(v10);
\draw[base] (v13)--(v03);
\draw[cA,line width=1.65pt] (v00)--(v03);
\draw[cA,line width=1.65pt] (v13)--(v12);
\draw[cA,line width=1.65pt] (v02)--(v01);
\draw[cA,line width=1.65pt] (v11)--(v10);
\draw[cB,line width=1.65pt] (v11)--(v01);
\draw[cB,line width=1.65pt] (v00)--(v10);
\draw[cB,line width=1.65pt] (v13)--(v03);
\draw[cB,line width=1.65pt] (v02)--(v12);
\draw[cC,line width=1.65pt] (v00)--(v01);
\draw[cC,line width=1.65pt] (v13)--(v10);
\draw[cC,line width=1.65pt] (v02)--(v03);
\draw[cC,line width=1.65pt] (v11)--(v12);
\node[bv] at (v00) {};
\node[wv] at (v01) {};
\node[bv] at (v02) {};
\node[wv] at (v03) {};
\node[wv] at (v10) {};
\node[bv] at (v11) {};
\node[wv] at (v12) {};
\node[bv] at (v13) {};
\node[font=\small,align=center] at (0,-2.05) {$0,1,2$};
\end{tikzpicture}
\hspace{.19cm}
\begin{tikzpicture}[scale=0.82,baseline=(current bounding box.center)]
\fill[black!3] (0,0) circle (.25);
\draw[black!20] (0,0) circle (.25);
\draw[cB!75,dashed,line width=.65pt] (.27,0)--(1.78,0);
\coordinate (v00) at (1.4,1.4);
\coordinate (v01) at (-1.4,1.4);
\coordinate (v02) at (-1.4,-1.4);
\coordinate (v03) at (1.4,-1.4);
\coordinate (v10) at (0.6,0.6);
\coordinate (v11) at (-0.6,0.6);
\coordinate (v12) at (-0.6,-0.6);
\coordinate (v13) at (0.6,-0.6);
\draw[base] (v00)--(v03);
\draw[base] (v00)--(v01);
\draw[base] (v00)--(v10);
\draw[base] (v02)--(v01);
\draw[base] (v02)--(v03);
\draw[base] (v02)--(v12);
\draw[base] (v11)--(v10);
\draw[base] (v11)--(v12);
\draw[base] (v11)--(v01);
\draw[base] (v13)--(v12);
\draw[base] (v13)--(v10);
\draw[base] (v13)--(v03);
\draw[cA,line width=1.65pt] (v00)--(v03);
\draw[cA,line width=1.65pt] (v13)--(v10);
\draw[cA,line width=1.65pt] (v02)--(v01);
\draw[cA,line width=1.65pt] (v11)--(v12);
\draw[cB,line width=1.65pt] (v00)--(v01);
\draw[cB,line width=1.65pt] (v02)--(v03);
\draw[cB,line width=1.65pt] (v11)--(v10);
\draw[cB,line width=1.65pt] (v13)--(v12);
\draw[cC,line width=1.65pt] (v11)--(v01);
\draw[cC,line width=1.65pt] (v00)--(v10);
\draw[cC,line width=1.65pt] (v13)--(v03);
\draw[cC,line width=1.65pt] (v02)--(v12);
\node[bv] at (v00) {};
\node[wv] at (v01) {};
\node[bv] at (v02) {};
\node[wv] at (v03) {};
\node[wv] at (v10) {};
\node[bv] at (v11) {};
\node[wv] at (v12) {};
\node[bv] at (v13) {};
\node[font=\small,align=center] at (0,-2.05) {$1,1,1$};
\end{tikzpicture}
\hspace{.19cm}
\begin{tikzpicture}[scale=0.82,baseline=(current bounding box.center)]
\fill[black!3] (0,0) circle (.25);
\draw[black!20] (0,0) circle (.25);
\draw[cB!75,dashed,line width=.65pt] (.27,0)--(1.78,0);
\coordinate (v00) at (1.4,1.4);
\coordinate (v01) at (-1.4,1.4);
\coordinate (v02) at (-1.4,-1.4);
\coordinate (v03) at (1.4,-1.4);
\coordinate (v10) at (0.6,0.6);
\coordinate (v11) at (-0.6,0.6);
\coordinate (v12) at (-0.6,-0.6);
\coordinate (v13) at (0.6,-0.6);
\draw[base] (v00)--(v03);
\draw[base] (v00)--(v01);
\draw[base] (v00)--(v10);
\draw[base] (v02)--(v01);
\draw[base] (v02)--(v03);
\draw[base] (v02)--(v12);
\draw[base] (v11)--(v10);
\draw[base] (v11)--(v12);
\draw[base] (v11)--(v01);
\draw[base] (v13)--(v12);
\draw[base] (v13)--(v10);
\draw[base] (v13)--(v03);
\draw[cA,line width=1.65pt] (v00)--(v03);
\draw[cA,line width=1.65pt] (v13)--(v10);
\draw[cA,line width=1.65pt] (v02)--(v01);
\draw[cA,line width=1.65pt] (v11)--(v12);
\draw[cB,line width=1.65pt] (v00)--(v01);
\draw[cB,line width=1.65pt] (v13)--(v03);
\draw[cB,line width=1.65pt] (v02)--(v12);
\draw[cB,line width=1.65pt] (v11)--(v10);
\draw[cC,line width=1.65pt] (v11)--(v01);
\draw[cC,line width=1.65pt] (v13)--(v12);
\draw[cC,line width=1.65pt] (v00)--(v10);
\draw[cC,line width=1.65pt] (v02)--(v03);
\node[bv] at (v00) {};
\node[wv] at (v01) {};
\node[bv] at (v02) {};
\node[wv] at (v03) {};
\node[wv] at (v10) {};
\node[bv] at (v11) {};
\node[wv] at (v12) {};
\node[bv] at (v13) {};
\node[font=\small,align=center] at (0,-2.05) {$1,1,1$};
\end{tikzpicture}
\hspace{.19cm}
\begin{tikzpicture}[scale=0.82,baseline=(current bounding box.center)]
\fill[black!3] (0,0) circle (.25);
\draw[black!20] (0,0) circle (.25);
\draw[cB!75,dashed,line width=.65pt] (.27,0)--(1.78,0);
\coordinate (v00) at (1.4,1.4);
\coordinate (v01) at (-1.4,1.4);
\coordinate (v02) at (-1.4,-1.4);
\coordinate (v03) at (1.4,-1.4);
\coordinate (v10) at (0.6,0.6);
\coordinate (v11) at (-0.6,0.6);
\coordinate (v12) at (-0.6,-0.6);
\coordinate (v13) at (0.6,-0.6);
\draw[base] (v00)--(v03);
\draw[base] (v00)--(v01);
\draw[base] (v00)--(v10);
\draw[base] (v02)--(v01);
\draw[base] (v02)--(v03);
\draw[base] (v02)--(v12);
\draw[base] (v11)--(v10);
\draw[base] (v11)--(v12);
\draw[base] (v11)--(v01);
\draw[base] (v13)--(v12);
\draw[base] (v13)--(v10);
\draw[base] (v13)--(v03);
\draw[cA,line width=1.65pt] (v00)--(v03);
\draw[cA,line width=1.65pt] (v11)--(v01);
\draw[cA,line width=1.65pt] (v13)--(v10);
\draw[cA,line width=1.65pt] (v02)--(v12);
\draw[cB,line width=1.65pt] (v00)--(v01);
\draw[cB,line width=1.65pt] (v02)--(v03);
\draw[cB,line width=1.65pt] (v11)--(v10);
\draw[cB,line width=1.65pt] (v13)--(v12);
\draw[cC,line width=1.65pt] (v13)--(v03);
\draw[cC,line width=1.65pt] (v00)--(v10);
\draw[cC,line width=1.65pt] (v02)--(v01);
\draw[cC,line width=1.65pt] (v11)--(v12);
\node[bv] at (v00) {};
\node[wv] at (v01) {};
\node[bv] at (v02) {};
\node[wv] at (v03) {};
\node[wv] at (v10) {};
\node[bv] at (v11) {};
\node[wv] at (v12) {};
\node[bv] at (v13) {};
\node[font=\small,align=center] at (0,-2.05) {$1,1,1$};
\end{tikzpicture}%
}

\newcommand{\NetworkFigure}{%
\begin{tikzpicture}[scale=1.25]
\fill[black!3] (0,0) circle (.18);
\draw[black!20] (0,0) circle (.18);

\coordinate (a) at (1.5,0);
\coordinate (b) at (-1.5,0);
\coordinate (c) at (0,.65);
\coordinate (d) at (0,-.65);

\draw[arr,cA] (a) to[out=90,in=90,looseness=1.25] (b);
\draw[arr,cA] (b) to[out=-90,in=-90,looseness=1.25] (a);

\draw[arr,cB] (c)
  .. controls (-.47,.36) and (-.47,-.36) .. (d);
\draw[arr,cB] (d)
  .. controls (.47,-.36) and (.47,.36) .. (c);

\draw[arr,cC] (a)--(c);
\draw[arr,cC] (c)--(b);
\draw[arr,cC] (b)--(d);
\draw[arr,cC] (d)--(a);

\draw[black!70,dashed,line width=.75pt]
  (.20,.05)--(1.95,.4875);
\node[above,font=\small] at (1.92,.51) {$\gamma$};

\foreach \p/\lab/\pos in
  {a/1/right,b/2/left,c/3/above,d/4/below}{
  \node[circle,draw,fill=white,inner sep=1.7pt,
        label=\pos:{\small$\lab$}] at (\p) {};
}
\end{tikzpicture}%
}

\newcommand{\StarFigure}{%

\begin{tikzpicture}[x=1.25cm,y=.65cm]
\foreach \y/\lab in {3/1,2/2,1/3,0/4}{
 \node[font=\small,left] at (0,\y) {$\lab$};
 \node[font=\small,right] at (5,\y) {$\lab$};}
\foreach \y in {3,2}{\draw[arr,cA] (0,\y)--(1.2,2.5);\draw[arr,cA] (1.2,2.5)--(2.4,\y);}
\foreach \y in {1,0}{\draw[arr,black!50] (0,\y)--(2.4,\y);}
\draw[arr,black!50] (2.4,3)--(5,3);
\foreach \y in {2,1,0}{\draw[arr,cB] (2.4,\y)--(3.7,1);\draw[arr,cB] (3.7,1)--(5,\y);}
\fill[cA] (1.2,2.5) circle (2pt);\fill[cB] (3.7,1) circle (2pt);
\node[font=\small] at (1.2,-.65) {$J_1=[1,2]$};
\node[font=\small] at (3.7,-.65) {$J_2=[2,4]$};
\end{tikzpicture}
}

\newcommand{\SplitStripFigure}{%
\begin{tikzpicture}[x=.96cm,y=.84cm,
  baseline=(current bounding box.center)]

\fill[black!2] (-.05,.13) rectangle (5.95,2.87);
\draw[black!30,line width=.65pt]
  (-.05,.13) rectangle (5.95,2.87);
\node[font=\scriptsize] at (0,3.18) {$\gamma$};
\node[font=\scriptsize] at (5.9,3.18) {$\gamma$};

\coordinate (sone) at (0,.50);
\coordinate (stwo) at (0,1.50);
\coordinate (sthree) at (0,2.50);
\coordinate (tone) at (5.90,.50);
\coordinate (ttwo) at (5.90,1.50);
\coordinate (tthree) at (5.90,2.50);

\coordinate (cw) at (1.18,1.15);
\coordinate (cb) at (1.56,1.15);
\coordinate (bw) at (2.17,2.17);
\coordinate (bb) at (2.55,2.17);
\coordinate (dw) at (2.57,.80);
\coordinate (db) at (2.95,.80);
\coordinate (aw) at (3.93,1.72);
\coordinate (ab) at (4.31,1.72);

\foreach \w/\b in {cw/cb,bw/bb,dw/db,aw/ab}{
  \draw[black!75,line width=2.2pt] (\w)--(\b);
}

\draw[arr,cC!85!black,line width=1.05pt] (stwo)--(cw);
\draw[arr,cC!85!black,line width=1.05pt] (cb)--(bw);
\draw[arr,cC!85!black,line width=1.05pt] (bb)--(dw);
\draw[arr,cC!85!black,line width=1.05pt] (db)--(aw);
\draw[arr,cC!85!black,line width=1.05pt] (ab)--(ttwo);

\draw[arr,cB,line width=1.55pt] (sone)--(cw);
\draw[arr,cB,line width=1.55pt] (cb)--(dw);
\draw[arr,cB,line width=1.55pt] (db)--(tone);

\draw[arr,cA,line width=1.55pt] (sthree)--(bw);
\draw[arr,cA,line width=1.55pt] (bb)--(aw);
\draw[arr,cA,line width=1.55pt] (ab)--(tthree);

\foreach \v in {sone,stwo,sthree,tone,ttwo,tthree}{
  \node[circle,draw=black!75,fill=white,
        inner sep=1.4pt] at (\v) {};
}

\foreach \v in {cw,bw,dw,aw}{
  \node[circle,draw=black!85,fill=white,
        inner sep=1.9pt] at (\v) {};
}
\foreach \v in {cb,bb,db,ab}{
  \node[circle,draw=black!85,fill=black,
        inner sep=1.9pt] at (\v) {};
}

\node[font=\scriptsize,above] at (1.37,1.42) {$3$};
\node[font=\scriptsize,above] at (2.36,2.44) {$2$};
\node[font=\scriptsize,below] at (2.76,.54) {$4$};
\node[font=\scriptsize,above] at (4.12,2.00) {$1$};

\end{tikzpicture}%
}

\newcommand{\NeutralSwitchPanel}[2]{%
\begin{tikzpicture}[scale=1.05,baseline=(current bounding box.center)]
  \useasboundingbox (-1.86,-1.77) rectangle (1.86,1.78);
  \coordinate (nw) at (-.76,.76);
  \coordinate (ne) at (.76,.76);
  \coordinate (se) at (.76,-.76);
  \coordinate (sw) at (-.76,-.76);

  \ifnum#2=1
    \fill[black!7] (-.76,-.76) rectangle (.76,.76);
    \foreach \p/\sx/\sy in
      {nw/-1/1,ne/1/1,se/1/-1,sw/-1/-1}{
      \fill[black!7]
        (\p)--({1.65*\sx},{.88*\sy})
            --({.88*\sx},{1.65*\sy})--cycle;
    }
  \fi

  \draw[web,webblue,line width=1.1pt]
    (nw)--(ne)--(se)--(sw)--cycle;

  \foreach \p/\sx/\sy/\beforeB in
    {nw/-1/1/0,ne/1/1/1,se/1/-1/0,sw/-1/-1/1}{

    \pgfmathtruncatemacro{\isB}{mod(\beforeB+#1,2)}
    \coordinate (endone) at ({1.65*\sx},{.88*\sy});
    \coordinate (endtwo) at ({.88*\sx},{1.65*\sy});

    \ifnum#2=0
      \ifnum\isB=0
        \draw[web,webblue,line width=1.1pt]
          (\p)--(endone) (\p)--(endtwo);
      \else
        \coordinate (outer) at ({1.12*\sx},{1.12*\sy});
        \draw[web,webblue,line width=1.1pt]
          (outer)--(endone) (outer)--(endtwo);
        \draw[web,line width=2.45pt] (\p)--(outer);

        \pgfmathtruncatemacro{\outerBlack}{mod(\isB+1,2)}
        \ifnum\outerBlack=1
          \node[bv] at (outer) {};
        \else
          \node[wv] at (outer) {};
        \fi
      \fi

      \pgfmathtruncatemacro{\cornerBlack}{mod(\beforeB+#1,2)}
      \ifnum\cornerBlack=1
        \node[bv] at (\p) {};
      \else
        \node[wv] at (\p) {};
      \fi

      \ifnum\isB=0
        \node[font=\scriptsize,fill=white,inner sep=.5pt]
          at ({.43*\sx},{.43*\sy}) {$A$};
      \else
        \node[font=\scriptsize,fill=white,inner sep=.5pt]
          at ({.43*\sx},{.43*\sy}) {$B$};
      \fi

    \else
      \draw[web,webblue,line width=1.1pt]
        (\p)--(endone) (\p)--(endtwo);

      \ifnum\isB=1
        \draw[web,webpurple,densely dashed,line width=1.3pt]
          (0,{1.52*\sy}) -- ({1.52*\sx},0);
      \fi

      \fill[black] (\p) circle (1.8pt);
    \fi
  }

  \ifnum#2=1
    \foreach \qx/\qy in
      {0/1.52,1.52/0,0/-1.52,-1.52/0}{
      \filldraw[draw=webpurple,fill=white,line width=.8pt]
        (\qx,\qy) circle (2pt);
    }
  \fi
\end{tikzpicture}%
}

\newcommand{\GadgetFigure}{%
\begin{tikzpicture}[x=1cm,y=1cm,baseline=(current bounding box.center)]
\foreach \y in {-0.65,0,.65}{
 \draw[arr,cA] (0,\y)--(1.2,0);
 \draw[arr,cA] (2.5,0)--(3.7,\y);
}
\draw[cB,thick] (1.2,0) to[bend left=27] (2.5,0);
\draw[cB,thick] (1.2,0)--(2.5,0);
\draw[cB,thick] (1.2,0) to[bend right=27] (2.5,0);
\node[wv] at (1.2,0) {};\node[bv] at (2.5,0) {};
\node[font=\small] at (1.85,.68) {$r-k$};
\draw[decorate,decoration={brace,mirror,amplitude=4pt}] (-.16,-.65)--(-.16,.65);
\node[font=\small,left] at (-.34,0) {$k$};
\node[font=\small] at (1.85,-1.05) {complementary colors};
\begin{scope}[xshift=6cm]
\foreach \y in {-0.65,0,.65}{
 \draw[arr,cA] (0,\y)--(1,0);
 \draw[arr,cA] (2.1,0)--(3.1,\y);
}
\node[wv] at (1,0) {};\node[bv] at (2.1,0) {};
\node[font=\small] at (1.55,-1.05) {$k=r=3$};
\end{scope}
\end{tikzpicture}%
}
\newcommand{\IntroWordFigure}[1]{%
\begin{tikzpicture}[x=1.05cm,y=.55cm,baseline=(current bounding box.center)]
\draw[black!45,dashed] (-.1,-.25)--(-.1,3.25);
\draw[black!45,dashed] (5.3,-.25)--(5.3,3.25);
\foreach \a/\lo/\hi/\lab in {0/1/2/{23},1.3/0/1/{34},2.6/2/3/{12},3.9/1/2/{23}}{
 \pgfmathsetmacro{\mid}{(\lo+\hi)/2}
 \foreach \y in {0,1,2,3}{
  \ifdim\y pt=\lo pt
   \ifnum#1=0
    \draw[cA,thick] (\a,\y)--({\a+.625},\mid)--({\a+1.3},\y);
   \else
    \draw[cA,thick] (\a,\y)--({\a+.35},\mid);
    \draw[cA,thick] ({\a+.9},\mid)--({\a+1.3},\y);
   \fi
  \else\ifdim\y pt=\hi pt
   \ifnum#1=0
    \draw[cA,thick] (\a,\y)--({\a+.625},\mid)--({\a+1.3},\y);
   \else
    \draw[cA,thick] (\a,\y)--({\a+.35},\mid);
    \draw[cA,thick] ({\a+.9},\mid)--({\a+1.3},\y);
   \fi
  \else
   \draw[cA,thick] (\a,\y)--({\a+1.3},\y);
  \fi\fi
 }
 \ifnum#1=0
  \node[circle,fill=cB,inner sep=1.5pt] at ({\a+.625},\mid) {};
 \else
  \draw[cB,thick] ({\a+.35},\mid)--({\a+.9},\mid);
  \node[wv,inner sep=1.4pt] at ({\a+.35},\mid) {};
  \node[bv,inner sep=1.4pt] at ({\a+.9},\mid) {};
 \fi
 \node[font=\scriptsize] at ({\a+.65},-.65) {$\lab$};
}
\foreach \y in {0,1,2,3}{\draw[cA,thick] (-.1,\y)--(0,\y);\draw[cA,thick] (5.2,\y)--(5.3,\y);}
\end{tikzpicture}%
}
\newcommand{\ThreeRulesFigure}{%
\begin{tikzpicture}[x=.85cm,y=.85cm,baseline=(current bounding box.center)]
\draw[thick] (0,0) circle (.37);
\node at (1,0) {$=3$};
\node[font=\small,anchor=west] at (2,0) {contractible circle};
\begin{scope}[yshift=-1.4cm]
\draw[thick] (-.75,0)--(0,0);
\draw[thick] (1.3,0)--(2.05,0);
\draw[thick] (0,0) to[bend left=35] (1.3,0);
\draw[thick] (0,0) to[bend right=35] (1.3,0);
\node[wv] at (0,0) {};\node[bv] at (1.3,0) {};
\node at (2.7,0) {$=2$};
\draw[thick] (3.35,0)--(5.4,0);
\end{scope}
\end{tikzpicture}\\[3pt]
\SquareColorsFig
}

\newcommand{\StretchFigure}{%
\begin{tikzpicture}[x=1cm,y=.72cm,
  baseline=(current bounding box.center)]
\foreach \a/\lab in {0/1,3.4/2,8.2/n}{
  \draw[fill=cA!4,rounded corners=2pt]
    ({\a+.75},-.3) rectangle ({\a+2.2},2.3);
  \node at ({\a+1.475},1) {$R^{(\lab)}$};
  \foreach \y in {0,1,2}{
    \draw[arr] (\a,\y)--({\a+.75},\y);
    \draw[arr] ({\a+2.2},\y)--({\a+2.9},\y);
  }
  \node[font=\small] at ({\a+1.475},-.8)
    {copy $\lab$};
}
\foreach \y in {0,1,2}{
  \draw (2.9,\y)--(3.4,\y);
}
\node at (7.25,1) {$\cdots$};
\end{tikzpicture}%
}

\definecolor{cD}{HTML}{7448C8}
\tikzset{st/.style={cA,thick},jn/.style={cB,thick},
 newe/.style={cD,line width=1.4pt},circ/.style={black,line width=1.5pt}}
\newcommand{\ladderbase}{%
\fill[black!3] (0,0) circle (.25);
\draw[black!20] (0,0) circle (.25);
\draw[cB!75,dashed,line width=.65pt] (.27,0)--(1.78,0);
\coordinate (v00) at (1.4,1.4);\coordinate (v01) at (-1.4,1.4);
\coordinate (v02) at (-1.4,-1.4);\coordinate (v03) at (1.4,-1.4);
\coordinate (v10) at (0.6,0.6);\coordinate (v11) at (-0.6,0.6);
\coordinate (v12) at (-0.6,-0.6);\coordinate (v13) at (0.6,-0.6);
\foreach \a/\b in {v00/v03,v00/v01,v00/v10,v02/v01,v02/v03,v02/v12,v11/v10,v11/v12,v11/v01,v13/v12,v13/v10,v13/v03}
 \draw[base] (\a)--(\b);}
\newcommand{\laddervertices}{%
\node[bv] at (v00) {};\node[wv] at (v01) {};\node[bv] at (v02) {};\node[wv] at (v03) {};
\node[wv] at (v10) {};\node[bv] at (v11) {};\node[wv] at (v12) {};\node[bv] at (v13) {};}
\newcommand{\MinimaFig}{%
\begin{tikzpicture}[scale=1.0,baseline=(current bounding box.center)]
\fill[black!3] (0,0) circle (.22);\draw[black!20] (0,0) circle (.22);
\foreach \j/\rr in {0/0.45,1/0.85,2/1.25,3/1.65}{
 \coordinate (p0\j) at (\rr,\rr);\coordinate (p1\j) at ({-\rr},\rr);
 \coordinate (p2\j) at ({-\rr},{-\rr});\coordinate (p3\j) at (\rr,{-\rr});}
\fill[cB!15] (p12)--(p22)--(p23)--(p13)--cycle;
\foreach \j in {0,1,2}{\foreach \i/\k in {0/1,1/2,2/3,3/0}\draw[base] (p\i\j)--(p\k\j);}
\foreach \i/\k in {1/2,2/3}\draw[base] (p\i3)--(p\k3);
\draw[black!30,dotted,line width=.9pt] (p03)--(p13);\draw[black!30,dotted,line width=.9pt] (p33)--(p03);
\foreach \i in {0,1,2,3}{\foreach \j/\k in {0/1,1/2,2/3}\draw[base] (p\i\j)--(p\i\k);}
\foreach \a/\b in {p00/p30,p02/p03,p11/p01,p13/p23,p20/p10,p22/p12,p31/p21,p33/p32}\draw[cA,line width=1.7pt] (\a)--(\b);
\draw[cB,line width=1.5pt,dashed] (p13)--(p12);\draw[cB,line width=1.5pt,dashed] (p22)--(p23);
\foreach \i in {0,1,2,3}{\foreach \j in {0,1,2,3}{
 \pgfmathtruncatemacro{\vpar}{mod(\i+\j,2)}
 \ifnum\vpar=0 \node[bv] at (p\i\j) {};\else\node[wv] at (p\i\j) {};\fi}}
\end{tikzpicture}%
}

\newcommand{\gadgetcoords}{%
\coordinate (w1) at (.35,1.5);\coordinate (b1) at (.9,1.5);
\coordinate (w2) at (1.65,.5);\coordinate (b2) at (2.2,.5);
\coordinate (w3) at (2.95,2.5);\coordinate (b3) at (3.5,2.5);
\coordinate (w4) at (4.25,1.5);\coordinate (b4) at (4.8,1.5);}
\newcommand{\gadgetframe}{%
\draw[black!45,dashed] (-.1,-.35)--(-.1,3.35);\draw[black!45,dashed] (5.3,-.35)--(5.3,3.35);
\foreach \i/\y in {1/3,2/2,3/1,4/0}\node[font=\scriptsize,left] at (-.12,\y) {$\i$};}
\newcommand{\sone}{\draw[st] (-.1,3)--(2.6,3)--(w3);\draw[st] (b3)--(3.9,3)--(5.3,3);}
\newcommand{\stwowrap}{\draw[st] (b4)--(5.2,2)--(5.3,2);\draw[st] (-.1,2)--(0,2)--(w1);}
\newcommand{\sthreewrap}{\draw[st] (b4)--(5.2,1)--(5.3,1);\draw[st] (-.1,1)--(0,1)--(w1);}
\newcommand{\stwomidA}{\draw[st] (b1)--(1.3,2)--(2.6,2)--(w3);}
\newcommand{\stwomidB}{\draw[st] (b3)--(3.9,2)--(w4);}
\newcommand{\sthreemidA}{\draw[st] (b1)--(1.3,1)--(w2);}
\newcommand{\sthreemidB}{\draw[st] (b2)--(2.6,1)--(3.9,1)--(w4);}
\newcommand{\sfour}{\draw[st] (b2)--(2.6,0)--(5.3,0);\draw[st] (-.1,0)--(1.3,0)--(w2);}
\newcommand{\ReduceA}{%
\begin{tikzpicture}[x=1.3cm,y=.6cm,baseline=(current bounding box.center)]
\gadgetframe\gadgetcoords
\fill[cC!20] (b4)--(5.2,2)--(5.3,2)--(5.3,1)--(5.2,1)--cycle;
\fill[cC!20] (-.1,2)--(0,2)--(w1)--(0,1)--(-.1,1)--cycle;
\sone\stwowrap\sthreewrap\stwomidA\stwomidB\sthreemidA\sthreemidB\sfour
\foreach \t in {1,2,3,4}\draw[jn] (w\t)--(b\t);
\foreach \v in {w1,w2,w3,w4}\node[wv] at (\v) {};
\foreach \v in {b1,b2,b3,b4}\node[bv] at (\v) {};
\node[font=\small] at (2.6,-.95) {(a) $H$: a bigon across the cut};
\end{tikzpicture}%
}
\newcommand{\ReduceB}{%
\begin{tikzpicture}[x=1.3cm,y=.6cm,baseline=(current bounding box.center)]
\gadgetframe\gadgetcoords
\fill[cA!12] (-.1,3)--(2.6,3)--(w3)--(2.6,2)--(1.3,2)--(b1)--(-.1,1.5)--cycle;
\fill[cA!12] (b3)--(3.9,3)--(5.3,3)--(5.3,1.5)--(w4)--(3.9,2)--cycle;
\sone\stwomidA\stwomidB\sthreemidA\sthreemidB\sfour
\draw[newe] (w4)--(5.3,1.5);\draw[newe] (-.1,1.5)--(b1);
\foreach \t in {2,3}\draw[jn] (w\t)--(b\t);
\foreach \v in {w2,w3,w4}\node[wv] at (\v) {};
\foreach \v in {b1,b2,b3}\node[bv] at (\v) {};
\node[font=\small] at (2.6,-.95) {(b) $2\times$: a square across the cut};
\end{tikzpicture}%
}
\newcommand{\ReduceC}{%
\begin{tikzpicture}[x=1.3cm,y=.6cm,baseline=(current bounding box.center)]
\gadgetframe\gadgetcoords
\fill[cC!20] (b2)--(2.6,1)--(3.9,1)--(4.25,1.25)--(5.3,1.25)--(5.3,0)--(2.6,0)--cycle;
\fill[cC!20] (-.1,1.25)--(.9,1.25)--(1.3,1)--(w2)--(1.3,0)--(-.1,0)--cycle;
\draw[circ] (w3)--(b3) to[out=10,in=180] (3.95,2.85)--(5.3,2.85);
\draw[circ] (-.1,2.85)--(2.5,2.85) to[out=0,in=170] (w3);
\draw[st] (b2)--(2.6,1)--(3.9,1)--(4.25,1.25)--(5.3,1.25);
\draw[st] (-.1,1.25)--(.9,1.25)--(1.3,1)--(w2);
\sfour
\draw[jn] (w2)--(b2);
\node[wv] at (w2) {};\node[bv] at (b2) {};
\node[font=\small] at (2.6,-.95) {(c) first smoothing: a circle and a bigon};
\end{tikzpicture}%
}
\newcommand{\ReduceD}{%
\begin{tikzpicture}[x=1.3cm,y=.6cm,baseline=(current bounding box.center)]
\gadgetframe\gadgetcoords
\fill[cC!20] (w2)--(b2) to[out=70,in=0] (1.925,1.75) to[out=180,in=110] (w2);
\draw[st] (b2) to[out=70,in=0] (1.925,1.75) to[out=180,in=110] (w2);
\sfour
\draw[jn] (w2)--(b2);
\node[wv] at (w2) {};\node[bv] at (b2) {};
\node[font=\small] at (2.6,-.95) {(d) second smoothing: a bigon};
\end{tikzpicture}%
}
\newcommand{\ReduceE}{%
\begin{tikzpicture}[scale=.62,baseline=(current bounding box.center)]
\fill[black!3] (0,0) circle (.3);\draw[black!20] (0,0) circle (.3);\draw[black!30] (0,0) circle (1.75);
\draw[circ] (0,0) circle (.75);\draw[circ] (0,0) circle (1.25);
\node[font=\small] at (0,-2.3) {from (c): $2\cdot e_2\cdot e_2$};
\end{tikzpicture}\qquad
\begin{tikzpicture}[scale=.62,baseline=(current bounding box.center)]
\fill[black!3] (0,0) circle (.3);\draw[black!20] (0,0) circle (.3);\draw[black!30] (0,0) circle (1.75);
\draw[circ] (0,0) circle (1.0);
\node[font=\small] at (0,-2.3) {from (d): $2\cdot e_1$};
\end{tikzpicture}%
}

\newcommand{\SquareColorsFig}{%
\begin{tikzpicture}[x=.9cm,y=.9cm,baseline=(current bounding box.center)]
\coordinate (B1) at (0,1.3);\coordinate (W1) at (1.3,1.3);\coordinate (B2) at (1.3,0);\coordinate (W2) at (0,0);
\draw[cB,line width=1.6pt] (B1)--(W1);\draw[cC,line width=1.6pt] (W1)--(B2);
\draw[cA,line width=1.6pt] (B2)--(W2);\draw[cC,line width=1.6pt] (W2)--(B1);
\draw[cA,line width=1.6pt] (B1)--(-.55,1.85);\draw[cA,line width=1.6pt] (W1)--(1.85,1.85);
\draw[cB,line width=1.6pt] (B2)--(1.85,-.55);\draw[cB,line width=1.6pt] (W2)--(-.55,-.55);
\node[bv] at (B1) {};\node[wv] at (W1) {};\node[bv] at (B2) {};\node[wv] at (W2) {};
\node[font=\small] at (-.75,2.05) {$a$};\node[font=\small] at (2.05,2.05) {$a$};
\node[font=\small] at (2.05,-.75) {$b$};\node[font=\small] at (-.75,-.75) {$b$};
\node at (2.75,.65) {$=$};
\begin{scope}[xshift=3.8cm]
\draw[cA,line width=1.6pt] (-.55,1.85) to[out=-45,in=-135] (1.85,1.85);
\draw[cB,line width=1.6pt] (-.55,-.55) to[out=45,in=135] (1.85,-.55);
\node[font=\small] at (.65,-1.05) {$1$};
\end{scope}
\node at (6.75,.65) {$+$};
\begin{scope}[xshift=7.7cm]
\draw[black!35,line width=1.2pt] (-.55,1.85) to[out=-45,in=45] (-.55,-.55);
\draw[black!35,line width=1.2pt] (1.85,1.85) to[out=-135,in=135] (1.85,-.55);
\node[font=\small] at (.65,-1.05) {$0$};
\end{scope}
\end{tikzpicture}%
}

\newcommand{\TwoSuperpositionsFigure}{%
\begin{tikzpicture}[scale=.82,baseline=(current bounding box.center)]
\ladderbase
\foreach \a/\b in {v00/v03,v02/v01,v11/v10,v13/v12}\draw[cA,line width=1.7pt] (\a)--(\b);
\foreach \a/\b in {v00/v01,v02/v03,v11/v12,v13/v10}\draw[cB,line width=1.7pt] (\a)--(\b);
\laddervertices
\node[font=\small,align=center] at (0,-2.1) {(a) two essential cycles\\ $F=(x_1+x_2)^2=e_1^2$};
\end{tikzpicture}
\hspace{1.1cm}
\begin{tikzpicture}[scale=.82,baseline=(current bounding box.center)]
\ladderbase
\fill[cC!14] (v00)--(v01)--(v11)--(v10)--cycle;
\foreach \a/\b in {v00/v01,v11/v10}\draw[cA,line width=1.7pt] (\a)--(\b);
\foreach \a/\b in {v00/v10,v11/v01}\draw[cB,line width=1.7pt] (\a)--(\b);
\draw[cA,line width=1.3pt] ($(v02)+(0,.07)$)--($(v03)+(0,.07)$);
\draw[cB,line width=1.3pt] ($(v02)+(0,-.07)$)--($(v03)+(0,-.07)$);
\draw[cA,line width=1.3pt] ($(v12)+(0,.07)$)--($(v13)+(0,.07)$);
\draw[cB,line width=1.3pt] ($(v12)+(0,-.07)$)--($(v13)+(0,-.07)$);
\laddervertices
\node[font=\small,align=center] at (0,-2.1) {(b) a contractible cycle\\ $F=2x_1x_2=2e_2$};
\end{tikzpicture}%
}
\newcommand{\StripFigure}{%
\begin{tikzpicture}[x=.96cm,y=.84cm,
  baseline=(current bounding box.center)]

\fill[black!2] (-.05,.13) rectangle (5.95,2.87);
\draw[black!30,line width=.65pt]
  (-.05,.13) rectangle (5.95,2.87);
\node[font=\scriptsize] at (0,3.18) {$\gamma$};
\node[font=\scriptsize] at (5.9,3.18) {$\gamma$};

\coordinate (sone) at (0,.50);
\coordinate (stwo) at (0,1.50);
\coordinate (sthree) at (0,2.50);
\coordinate (tone) at (5.90,.50);
\coordinate (ttwo) at (5.90,1.50);
\coordinate (tthree) at (5.90,2.50);

\coordinate (vc) at (1.37,1.15); % vertex 3
\coordinate (vb) at (2.36,2.17); % vertex 2
\coordinate (vd) at (2.76,.80);  % vertex 4
\coordinate (va) at (4.12,1.72); % vertex 1

\draw[arr,cC!85!black,line width=1.05pt]
  (stwo)--(vc);
\draw[arr,cC!85!black,line width=1.05pt]
  (vc)--(vb);
\draw[arr,cC!85!black,line width=1.05pt]
  (vb)--(vd);
\draw[arr,cC!85!black,line width=1.05pt]
  (vd)--(va);
\draw[arr,cC!85!black,line width=1.05pt]
  (va)--(ttwo);

\draw[arr,cB,line width=1.55pt] (sone)--(vc);
\draw[arr,cB,line width=1.55pt] (vc)--(vd);
\draw[arr,cB,line width=1.55pt] (vd)--(tone);

\draw[arr,cA,line width=1.55pt] (sthree)--(vb);
\draw[arr,cA,line width=1.55pt] (vb)--(va);
\draw[arr,cA,line width=1.55pt] (va)--(tthree);

\foreach \v in {sone,stwo,sthree,tone,ttwo,tthree}{
  \node[circle,draw=black!75,fill=white,
        inner sep=1.4pt] at (\v) {};
}

\foreach \v/\num/\where in
  {vc/3/above left,vb/2/above,
   vd/4/below,va/1/below}{
  \node[circle,draw=black!75,fill=white,
        inner sep=1.6pt] at (\v) {};
  \node[font=\scriptsize,\where] at (\v) {$\num$};
}

\end{tikzpicture}%
}
\newcommand{\RoadmapFigure}{%
\begin{tikzpicture}[x=1cm,y=1cm,
 box/.style={draw,rounded corners=2pt,align=center,font=\small,inner sep=5pt,text width=2.9cm},
 lab/.style={font=\scriptsize,align=center,fill=white,inner sep=1pt}]
\node[box,fill=cA!6] (D) at (0,0) {cylindrical dimers\\$a_0,a_1,\ldots$};
\node[box,fill=cA!6] (N) at (0,2.7) {cylindrical networks\\directed cycles};
\node[box,fill=cC!8] (Z) at (5.2,0) {annular colorings\\$Z_H$};
\node[box,fill=cB!6] (R) at (5.2,2.7) {planar path matrices\\immanants};
\node[box,fill=cB!6] (U) at (10.4,2.7) {interval products\\$U=X_{J_b}\cdots X_{J_1}$};
\node[box,fill=cC!8] (Q) at (10.4,0) {antisymmetrizer closures\\$\Phi(Q)$};
\draw[arr] (D)--node[lab,right]{complete to a grid\\and contract\\Secs.~\ref{sec:networks}--\ref{sec:annular}}(N);
\draw[arr] (N)--node[lab,above]{cut open, LGV\\Sec.~\ref{sec:networks}}(R);
\draw[arr] (U.west) to[bend right=15] node[lab,above]{stars\\Sec.~\ref{sec:schur}} (R.east);
\draw[arr] (R.east) to[bend right=15] node[lab,below]{Wolfgang\\Sec.~\ref{sec:schur}} (U.west);
\draw[arr] (U)--node[lab,left]{sign twist\\Sec.~\ref{sec:schur}}(Q);
\draw[arr] (D)--node[lab,above]{superpose $r$\\Sec.~\ref{sec:annular}}(Z);
\draw[arr] (Q.west) to[bend right=15] node[lab,above]{gadget\\Sec.~\ref{sec:colorings}} (Z.east);
\draw[arr] (Z.east) to[bend right=15] node[lab,below]{back to closures\\Sec.~\ref{sec:annular}} (Q.west);
\end{tikzpicture}%
}

\colorlet{webblue}{cA}\colorlet{weborange}{cB}
\colorlet{webgreen}{cC}\colorlet{webpurple}{cD}
\tikzset{web/.style={line width=.8pt,line cap=round,line join=round}}
\newcommand{\FourSmallFill}[4]{%
\begin{tikzpicture}[scale=.85,baseline=(current bounding box.center)]
 \coordinate (nw) at (-.6,.6);\coordinate (ne) at (.6,.6);
 \coordinate (se) at (.6,-.6);\coordinate (sw) at (-.6,-.6);
 \ifnum#4=0
  \pgfmathtruncatemacro{\ll}{#1-#3}\pgfmathtruncatemacro{\tt}{#3}
  \pgfmathtruncatemacro{\rr}{4-#2-#3}\pgfmathtruncatemacro{\bb}{#3}
  \foreach \v in {nw,se}{\node[wv] at (\v) {};}
  \foreach \v in {ne,sw}{\node[bv] at (\v) {};}
  \draw[web] (ne)--(1.0,1.0) (sw)--(-1.0,-1.0);
  \node[wv] at (1.0,1.0) {};\node[wv] at (-1.0,-1.0) {};
  \draw[web,weborange] (nw)--(-1.25,1.25) (1,1)--(1.25,1.25)
    (se)--(1.25,-1.25) (-1,-1)--(-1.25,-1.25);
 \else
  \pgfmathtruncatemacro{\ll}{4-#1-#3}\pgfmathtruncatemacro{\tt}{#3}
  \pgfmathtruncatemacro{\rr}{#2-#3}\pgfmathtruncatemacro{\bb}{#3}
  \foreach \v in {ne,sw}{\node[wv] at (\v) {};}
  \foreach \v in {nw,se}{\node[bv] at (\v) {};}
  \draw[web] (nw)--(-1.0,1.0) (se)--(1.0,-1.0);
  \node[wv] at (-1.0,1.0) {};\node[wv] at (1.0,-1.0) {};
  \draw[web,weborange] (-1,1)--(-1.25,1.25) (ne)--(1.25,1.25)
    (1,-1)--(1.25,-1.25) (sw)--(-1.25,-1.25);
 \fi
 \foreach \a/\b/\n/\pos in {nw/sw/\ll/left,nw/ne/\tt/above,ne/se/\rr/right,sw/se/\bb/below}{
  \ifnum\n>0
   \pgfmathsetlengthmacro{\lw}{(.4+.45*\n)*1pt}
   \draw[web,webblue,line width=\lw] (\a)--(\b)
     node[midway,\pos,font=\scriptsize,fill=white,inner sep=.4pt] {$\n$};
  \else
   \draw[densely dotted,black!25] (\a)--(\b);
  \fi}
 \pgfmathtruncatemacro{\pa}{4-#1}\pgfmathtruncatemacro{\pb}{4-#2}
 \node[font=\scriptsize,weborange,above left] at (-1.25,1.25) {$\pa$};
 \node[font=\scriptsize,weborange,above right] at (1.25,1.25) {$\pb$};
 \node[font=\scriptsize,weborange,below right] at (1.25,-1.25) {$#2$};
 \node[font=\scriptsize,weborange,below left] at (-1.25,-1.25) {$#1$};
\end{tikzpicture}}

\newcommand{\FourContractionFigure}{%
\begin{tikzpicture}[scale=.9]
 \draw[web,webblue] (-1.1,.7)--(-.45,0)--(-1.1,-.7);
 \draw[web,webblue] (1.1,.7)--(.45,0)--(1.1,-.7);
 \draw[web,line width=2.4pt] (-.45,0)--(.45,0)
   node[midway,below,font=\small] {$2$};
 \node[bv] at (-.45,0) {};\node[wv] at (.45,0) {};
 \draw[web,webpurple,dashed] (0,-.95)--(0,.95);
 \foreach \y in {-.95,.95}{\node[draw=webpurple,circle,fill=white,inner sep=2pt] at (0,\y) {};}
 \node[font=\small] at (0,-1.4) {two cubic vertices};
 \draw[arr] (1.7,0)--(2.6,0);
 \begin{scope}[xshift=4.1cm]
  \draw[web,webblue] (-1.1,.7)--(0,0)--(1.1,-.7);
  \draw[web,webblue] (-1.1,-.7)--(0,0)--(1.1,.7);
  \draw[web,webpurple,dashed] (0,-.95)--(0,.95);
  \foreach \y in {-.95,.95}{\node[draw=webpurple,circle,fill=white,inner sep=2pt] at (0,\y) {};}
  \fill (0,0) circle(2pt);
  \node[font=\small] at (0,-1.4) {one type-$B$ vertex};
 \end{scope}
\end{tikzpicture}%
}

\newcommand{\FourEmptyFigure}{%
\begin{tikzpicture}[scale=.85]
 \fill[webpurple!6] (0,0) ellipse (1.65 and 1.1);
 \draw[web,webpurple,line width=1.2pt] (0,0) ellipse (1.65 and 1.1);
 \foreach \a in {20,85,155,220,285,335}
  \fill[webpurple] ({1.65*cos(\a)},{1.1*sin(\a)}) circle(2pt);
 \foreach \x/\y in {-.8/.35,.1/.55,.85/.05,-.35/-.5,.55/-.55}
  {\node[draw=webgreen,circle,fill=white,inner sep=2pt] at (\x,\y) {};}
 \node[webpurple,font=\small] at (0,1.45) {a component $K$};
 \node[align=center,font=\small,anchor=north,text width=3.9cm] at (0,-1.4)
  {flip the shaded\\squares inside};
 \draw[arr] (2.15,0)--(3.15,0);
 \begin{scope}[xshift=6cm]
  \draw[web,black!30,dashed] (0,0) ellipse(1.65 and 1.1);
  \foreach \a in {20,85,155,220,285,335}
   {\node[draw=webpurple,circle,fill=white,inner sep=1.5pt]
      at ({1.65*cos(\a)},{1.1*sin(\a)}) {};}
  \foreach \a in {52,120,187,252,310,357}
   {\fill[weborange] ({1.65*cos(\a)-.06},{1.1*sin(\a)-.06}) rectangle ++(.12,.12);}
  \foreach \x/\y in {-.8/.35,.1/.55,.85/.05,-.35/-.5,.55/-.55}
   {\node[draw=webgreen,circle,fill=webgreen,inner sep=2pt] at (\x,\y) {};}
  \node[align=center,font=\small,anchor=north,text width=3.9cm] at (0,-1.4)
   {boundary heavy\\edges disappear};
 \end{scope}
\end{tikzpicture}%
}

\newcommand{\FourPermutohedronFigure}{%
\begin{tikzpicture}[scale=.88]
\coordinate (v2) at (0.000000,3.000000);
\coordinate (v3) at (-2.598076,1.500000);
\coordinate (v5) at (-2.598076,-1.500000);
\coordinate (v4) at (-0.000000,-3.000000);
\coordinate (v1) at (2.598076,-1.500000);
\coordinate (v0) at (2.598076,1.500000);
\coordinate (v6) at (1.657916,0.513619);
\coordinate (v7) at (1.657916,-0.513619);
\coordinate (v8) at (0.717756,0.554475);
\coordinate (v9) at (0.374042,0.250973);
\coordinate (v10) at (0.717756,-0.554475);
\coordinate (v11) at (0.374042,-0.250973);
\coordinate (v12) at (-0.384151,1.692607);
\coordinate (v13) at (-1.273765,1.178988);
\coordinate (v14) at (0.121311,0.898833);
\coordinate (v15) at (0.030328,0.449416);
\coordinate (v16) at (-0.839067,0.344358);
\coordinate (v17) at (-0.404370,0.198444);
\coordinate (v18) at (-0.384151,-1.692607);
\coordinate (v19) at (-1.273765,-1.178988);
\coordinate (v20) at (0.121311,-0.898833);
\coordinate (v21) at (0.030328,-0.449416);
\coordinate (v22) at (-0.839067,-0.344358);
\coordinate (v23) at (-0.404370,-0.198444);
\fill[weborange!10] (v9)--(v11)--(v21)--(v23)--(v17)--(v15)--cycle;
\draw[webblue,line width=.8pt] (v0)--(v6);
\draw[weborange,line width=1.8pt] (v0)--(v2);
\draw[webblue,line width=.8pt] (v0)--(v1);
\draw[webblue,line width=.8pt] (v1)--(v7);
\draw[weborange,line width=1.8pt] (v1)--(v4);
\draw[webblue,line width=.8pt] (v2)--(v12);
\draw[webblue,line width=.8pt] (v2)--(v3);
\draw[webblue,line width=.8pt] (v3)--(v13);
\draw[weborange,line width=1.8pt] (v3)--(v5);
\draw[webblue,line width=.8pt] (v4)--(v18);
\draw[webblue,line width=.8pt] (v4)--(v5);
\draw[webblue,line width=.8pt] (v5)--(v19);
\draw[weborange,line width=1.8pt] (v6)--(v8);
\draw[webblue,line width=.8pt] (v6)--(v7);
\draw[weborange,line width=1.8pt] (v7)--(v10);
\draw[webblue,line width=.8pt] (v8)--(v14);
\draw[webblue,line width=.8pt] (v8)--(v9);
\draw[webblue,line width=.8pt] (v9)--(v15);
\draw[weborange,line width=1.8pt] (v9)--(v11);
\draw[webblue,line width=.8pt] (v10)--(v20);
\draw[webblue,line width=.8pt] (v10)--(v11);
\draw[webblue,line width=.8pt] (v11)--(v21);
\draw[weborange,line width=1.8pt] (v12)--(v14);
\draw[webblue,line width=.8pt] (v12)--(v13);
\draw[weborange,line width=1.8pt] (v13)--(v16);
\draw[webblue,line width=.8pt] (v14)--(v15);
\draw[weborange,line width=1.8pt] (v15)--(v17);
\draw[webblue,line width=.8pt] (v16)--(v22);
\draw[webblue,line width=.8pt] (v16)--(v17);
\draw[webblue,line width=.8pt] (v17)--(v23);
\draw[weborange,line width=1.8pt] (v18)--(v20);
\draw[webblue,line width=.8pt] (v18)--(v19);
\draw[weborange,line width=1.8pt] (v19)--(v22);
\draw[webblue,line width=.8pt] (v20)--(v21);
\draw[weborange,line width=1.8pt] (v21)--(v23);
\draw[webblue,line width=.8pt] (v22)--(v23);
\node[bv] at (v0) {};
\node[wv] at (v1) {};
\node[wv] at (v2) {};
\node[bv] at (v3) {};
\node[bv] at (v4) {};
\node[wv] at (v5) {};
\node[wv] at (v6) {};
\node[bv] at (v7) {};
\node[bv] at (v8) {};
\node[wv] at (v9) {};
\node[wv] at (v10) {};
\node[bv] at (v11) {};
\node[bv] at (v12) {};
\node[wv] at (v13) {};
\node[wv] at (v14) {};
\node[bv] at (v15) {};
\node[bv] at (v16) {};
\node[wv] at (v17) {};
\node[wv] at (v18) {};
\node[bv] at (v19) {};
\node[bv] at (v20) {};
\node[wv] at (v21) {};
\node[wv] at (v22) {};
\node[bv] at (v23) {};
\node[font=\scriptsize,fill=weborange!10,inner sep=.4pt] at (0,0) {$e_1$};
\node[font=\small] at (0,3.5) {outer end $e_2$};
\end{tikzpicture}%
}

\newcommand{\FourCubeFigure}{%
\begin{tikzpicture}[scale=.95]
\begin{scope}[xshift=0.0cm]
\coordinate (n6) at (-0.5,0.5);
\coordinate (n9) at (-1.1,1.1);
\coordinate (n7) at (0.5,0.5);
\coordinate (n12) at (1.1,1.1);
\coordinate (n10) at (0.5,-0.5);
\coordinate (n11) at (1.1,-1.1);
\coordinate (n8) at (-0.5,-0.5);
\coordinate (n13) at (-1.1,-1.1);
\draw[web,cB,line width=1.4pt] (n6)--(n7);
\draw[web,cB,line width=1.4pt] (n6)--(n8);
\draw[web,cD] (n6)--(n9);
\draw[web,cB,line width=1.4pt] (n7)--(n10);
\draw[web,cB,line width=1.4pt] (n8)--(n10);
\draw[web,cD] (n10)--(n11);
\draw[web,cD] (n7)--(n12);
\draw[web,cD] (n9)--(n12);
\draw[web,cD] (n11)--(n12);
\draw[web,cD] (n8)--(n13);
\draw[web,cD] (n9)--(n13);
\draw[web,cD] (n11)--(n13);
\node[circle,draw=cD,fill=cD,inner sep=1.5pt] at (n6) {};
\node[circle,draw=cD,fill=cD,inner sep=1.5pt] at (n9) {};
\node[circle,draw=cD,fill=cD,inner sep=1.5pt] at (n7) {};
\node[circle,draw=cD,fill=cD,inner sep=1.5pt] at (n12) {};
\node[circle,draw=cD,fill=cD,inner sep=1.5pt] at (n10) {};
\node[circle,draw=cD,fill=cD,inner sep=1.5pt] at (n11) {};
\node[circle,draw=cD,fill=cD,inner sep=1.5pt] at (n8) {};
\node[circle,draw=cD,fill=cD,inner sep=1.5pt] at (n13) {};
\draw[black] (n6) circle(.12);
\node[font=\scriptsize,above right] at (-0.44,0.56) {$ e_1 $};
\draw[black] (n11) circle(.12);
\node[font=\scriptsize,above right] at (1.1600000000000001,-1.04) {$ e_2 $};
\node[draw=cC,fill=white,rectangle,inner sep=1.7pt] at (0,0) {};
\node[font=\small] at (0,1.5) {before};
\end{scope}
\begin{scope}[xshift=3.9cm]
\coordinate (n6) at (-0.5,0.5);
\coordinate (n9) at (-1.1,1.1);
\coordinate (n7) at (0.5,0.5);
\coordinate (n12) at (1.1,1.1);
\coordinate (n10) at (0.5,-0.5);
\coordinate (n11) at (1.1,-1.1);
\coordinate (n8) at (-0.5,-0.5);
\coordinate (n13) at (-1.1,-1.1);
\draw[web,cD] (n6)--(n9);
\draw[web,cD] (n10)--(n11);
\draw[web,cD] (n7)--(n12);
\draw[web,cD] (n9)--(n12);
\draw[web,cD] (n11)--(n12);
\draw[web,cD] (n8)--(n13);
\draw[web,cD] (n9)--(n13);
\draw[web,cD] (n11)--(n13);
\node[circle,draw=cD,fill=cD,inner sep=1.5pt] at (n6) {};
\node[circle,draw=cD,fill=cD,inner sep=1.5pt] at (n9) {};
\node[circle,draw=cB,fill=cB,inner sep=1.5pt] at (n7) {};
\node[circle,draw=cD,fill=cD,inner sep=1.5pt] at (n12) {};
\node[circle,draw=cB,fill=cB,inner sep=1.5pt] at (n10) {};
\node[circle,draw=cD,fill=cD,inner sep=1.5pt] at (n11) {};
\node[circle,draw=cB,fill=cB,inner sep=1.5pt] at (n8) {};
\node[circle,draw=cD,fill=cD,inner sep=1.5pt] at (n13) {};
\draw[black] (n6) circle(.12);
\node[font=\scriptsize,above right] at (-0.44,0.56) {$ e_1 $};
\draw[black] (n11) circle(.12);
\node[font=\scriptsize,above right] at (1.1600000000000001,-1.04) {$ e_2 $};
\node[font=\small] at (0,1.5) {after};
\end{scope}
\draw[arr] (1.45,0)--(2.45,0);
\node[font=\small,align=center] at (1.95,-1.65) {switch the square inside the inner circuit};
\end{tikzpicture}%
}

\title{Cylindrical dimers and monomial immanants}
\author{Pavlo Pylyavskyy}
\date{}
\begin{document}
\maketitle

\begin{abstract}
For a bipartite graph on a cylinder, group its perfect matchings by homological height and sum their weights in each group. We prove that every minor of the resulting Toeplitz matrix has nonnegative coefficients in the edge weights, settling a question recorded by Lai. We first prove conjectures of Galashin and the author concerning directed cycles in cylindrical networks, and then deduce the dimer theorem from them. A related, stronger positivity property for superpositions of cylindrical matchings turns out to be equivalent to Stembridge's interval-product conjecture on monomial immanants. Using annular webs, we prove this conjecture for partitions with at most four columns. A stretching argument extends this to all partitions with parts in $\{n,2n,3n,4n\}$.
\end{abstract}

\section{Introduction}

Perfect matchings of graphs on surfaces carry
topological information: the difference of two
matchings consists of alternating cycles, which
may wind around the holes of the surface.
We can therefore group matchings by homology
and study how the weights of these groups are
related. Kasteleyn's Pfaffian method made exact
enumeration of planar dimers possible, and its
extension to surfaces brought topology into
the picture \cite{Kasteleyn63,CimRes}.
For bipartite dimers on a torus, the two winding
directions give rise to a spectral curve.
Kenyon, Okounkov, and Sheffield proved that
this curve is Harnack and used its geometry
to describe the phase diagram \cite{KOS}.
Kenyon and Okounkov also showed that every
Harnack curve can arise this way \cite{KO}.
This connection has since led to cluster
integrable systems and explicit dimer models
for higher-genus Harnack curves
\cite{GonKen,BCdT}.

On a cylinder there is only one winding direction,
so the weights of the homology classes form a
polynomial in one variable. Saldanha and Tomei
proved that this polynomial has only real roots
for domino tilings of quadriculated annuli
\cite{ST}. By the theorem of Aissen, Schoenberg,
and Whitney \cite{ASW}, a polynomial with
nonnegative coefficients is real-rooted exactly
when its coefficient sequence is P\'olya frequency,
meaning that every minor of the associated
Toeplitz matrix is nonnegative.

We ask for a stronger property: if the edge weights
are independent variables, do these minors have
nonnegative coefficients? This kind of
coefficientwise total positivity has been studied
in \cite{Brenti95,CDDGS}. The author posed the question for cylindrical dimers at
the 2020 Twin Cities REU. Cheng, Peng, and Srinivasan
recorded it in their report \cite[Question~1.9]{CPS}.
Working with Chris Fraser and Eric Stucky, they proved
several small-minor cases and real-rootedness for
unweighted cylindrical grids. Lai later recorded
the general question in his survey on tilings
\cite[Problem~32]{Lai}. In this paper we prove coefficientwise
P\'olya frequency property for every bipartite dimer model
on a cylinder.

Our proof starts with a conjecture that Pavel Galashin
and the author made while studying linear recurrences
for cylindrical networks \cite{GP}. In that setting,
the coefficients count families of vertex-disjoint
directed cycles. Cutting the cylinder open turns
these cycles into paths in a planar network, where
we can use the classical theory of nonintersecting
paths \cite{Lindstrom,GV}. We prove the conjecture, and then deduce the dimer theorem by
embedding the original graph into a cylindrical
grid.

Planar path matrices also lead to the second
positivity question of this paper. Recall that
an immanant is defined like a determinant, but
with the sign of each permutation replaced by
a character value. The sign and trivial characters
give the determinant and permanent. Goulden--Jackson,
Greene, and Wolfgang established coefficientwise
positivity for irreducible immanants of compatible
planar path matrices \cite{GJ,Greene,Wolfgang}.
Stembridge conjectured the same for monomial
immanants \cite{Stembridge92}. Equivalently, he
asked whether monomial characters are nonnegative
on products of interval sums in the group algebra.
A stronger conjecture of Haiman concerns monomial
characters on Kazhdan--Lusztig basis elements
$C'_w$ \cite{Haiman}. At $q=1$, it would imply
Stembridge's conjecture, since interval sums are
specializations of parabolic Kazhdan--Lusztig
elements, and their products expand positively
in this basis \cite{EW}.

An important special case connects this story to the
Stanley--Stembridge conjecture \cite{SS93}.
For a smooth permutation $w$, the Kazhdan--Lusztig
element $C'_w(1)$ is, up to a positive scalar,
a product of interval sums \cite{Skandera08}.
For codominant permutations, the corresponding
monomial-character evaluations are the elementary
coefficients of Stanley's chromatic symmetric
functions of unit interval orders
\cite{StanleyChrom,CHSS}. Guay-Paquet reduced
the Stanley--Stembridge conjecture to this setting,
and Hikita recently proved it
\cite{GuayPaquet,Hikita}. Moreover, every smooth
permutation has the same character evaluations
as some codominant permutation \cite{CHSS,AN}.
Thus Hikita's result proves Haiman's conjecture
at $q=1$ for all smooth permutations and all
partitions. In another direction, McNulty recently
established the $q=1$ conjecture for all hook
shapes \cite{McNulty}.

The cylinder enters the story again through a
connection with edge colorings. Closing the
permutation diagrams for interval products gives
graphs on an annulus, whose proper edge colorings
encode the corresponding character evaluations.
Such colorings also arise naturally by superimposing
perfect matchings, with one color for each matching.
Recording the winding of each color gives a
symmetric polynomial. Its Schur positivity is
closely related to the dimer theorem, while the
stronger elementary positivity for $r$ colors
is equivalent to Stembridge's interval-product
conjecture for partitions with at most $r$ columns.

These colorings are closely related to the theory
of webs. Kuperberg introduced local diagrammatic
rules for tensor invariants, later developed in
type $A$ by Cautis, Kamnitzer, and Morrison
\cite{Kuperberg,CKM}. Webs can also be drawn on
surfaces, where they appear in skein algebras
and arise naturally from superpositions of dimers
\cite{FrohmanSikora,FLL,DKS}. For annular webs,
Queffelec and Rose conjectured a decomposition
that would imply elementary positivity.
Gorsky and Wedrich proved the weaker Schur
positivity, while Oh and Yoon gave a direct
description of the associated symmetric functions
in terms of edge colorings \cite{QR,GW,OY}.

Our local reductions build on the positive version
of Kuperberg's rules used by the author in
\emph{$A_2$-web immanants} \cite{PWeb}.
At $q=1$, these rules become identities between
coloring counts. Douglas, Kenyon, and Shi used
them to reduce three-color annular webs
\cite[Proposition~6.1]{DKS}.
For four colors, we use the relations of
Fraser--Lam--Le \cite{FLL}. The main difficulty
is that some square moves only rearrange the
web without simplifying it. We show that such
moves always lead to a genuine reduction.
This proves elementary positivity through four
colors, and hence Stembridge's conjecture for
partitions with at most four columns.
A separate path-matrix argument shows that
positivity for a partition $\mu$ implies
positivity for every stretch $n\mu$.

\subsection{The cylinder theorems}

Let $\cyl$ be an annulus, or equivalently a cylinder, and
fix a generator of $H_1(\cyl;\mathbb Z)\cong\mathbb Z$.
Throughout, graphs are finite and embedded without
crossings; parallel edges are allowed.

For two perfect matchings $M,M'$ of a bipartite graph
in $\cyl$, orient the edges of $M$ from black to white
and those of $M'$ from white to black. The resulting
integral $1$-chain is a cycle, since at every vertex
one edge enters and one leaves. Denote its homology
class by $h(M,M')\in\mathbb Z$.

We have
\[
h(M,M'')=h(M,M')+h(M',M''),
\]
simply because the corresponding $1$-chains add:
the contributions of $M'$ cancel. Thus changing the
reference matching shifts all heights by the same
integer. In particular, if we set
\[
h(M)=h(M,M_0)-\min_N h(N,M_0),
\]
then $h(M)$ is independent of the choice of $M_0$.

Give edge $e$ a variable $y_e$, and put
\begin{equation}\label{eq:measurements}
 a_k=\sum_{h(M)=k}\wt(M),\qquad
 \wt(M)=\prod_{e\in M}y_e,\qquad A(z)=\sum_{k\ge0}a_kz^k.
\end{equation}
Set $a_k=0$ outside the finite support. For a sequence $b$, set $b_j=0$ for $j<0$. We call it \emph{coefficientwise P\'olya frequency} if every finite minor of
\[
 T(b)=(b_{j-i})_{i,j\ge0}
\]
has nonnegative coefficients. Here and below total positivity statements allow minors to vanish.

Our first main result settles the cylindrical dimer question
in its stronger, coefficientwise form.

\begin{theorem}\label{thm:dimer}
For every bipartite graph embedded in a cylinder and
admitting a perfect matching, the sequence $(a_k)$
in \eqref{eq:measurements} is coefficientwise
P\'olya frequency.
\end{theorem}

In particular, $A(z)$ has only negative real roots
for positive numerical edge weights.

We also prove an analogous result for directed networks.
A \emph{planar cylindrical network} is a directed graph
in $\cyl$ for which every directed cycle has positive
winding number. Every simple directed cycle then winds
once. Define
\begin{equation}\label{eq:H}
 H_k=\sum_{\mathcal C}\prod_{e\in\mathcal C}y_e,
 \qquad H_0=1,
\end{equation}
where $\mathcal C$ runs over unordered families of
$k$ vertex-disjoint simple directed cycles.

\begin{theorem}\label{thm:network}
For a planar cylindrical network, $(H_k)$ is
coefficientwise P\'olya frequency.
\end{theorem}

Theorems~\ref{thm:dimer} and~\ref{thm:network} are
equivalent by vertex expansion in Section~\ref{sec:networks}
and grid completion in Section~\ref{sec:annular}.
Thus the coefficientwise dimer question recorded by Lai
is equivalent to Conjectures~6.1 and~6.2 of
Galashin--Pylyavskyy \cite{GP}. Their real-rootedness
Conjecture~6.3 follows as well.

\begin{example} \label{ex:ladder}
The graph $C_4\times P_2$, drawn as two nested squares joined at corresponding corners, has nine perfect matchings: one of height $0$, seven of height $1$, and one of height $2$. Thus, at unit weights,
\[
 A(z)=1+7z+z^2.
\]
\end{example}

A convenient way to compute heights is to choose a
\emph{cut}: a simple curve connecting the two boundary
components of the cylinder, crossing edges transversely
and avoiding vertices. The height difference between
two matchings can then be read off from their signed
edge crossings with the cut.

\begin{figure}[htbp]
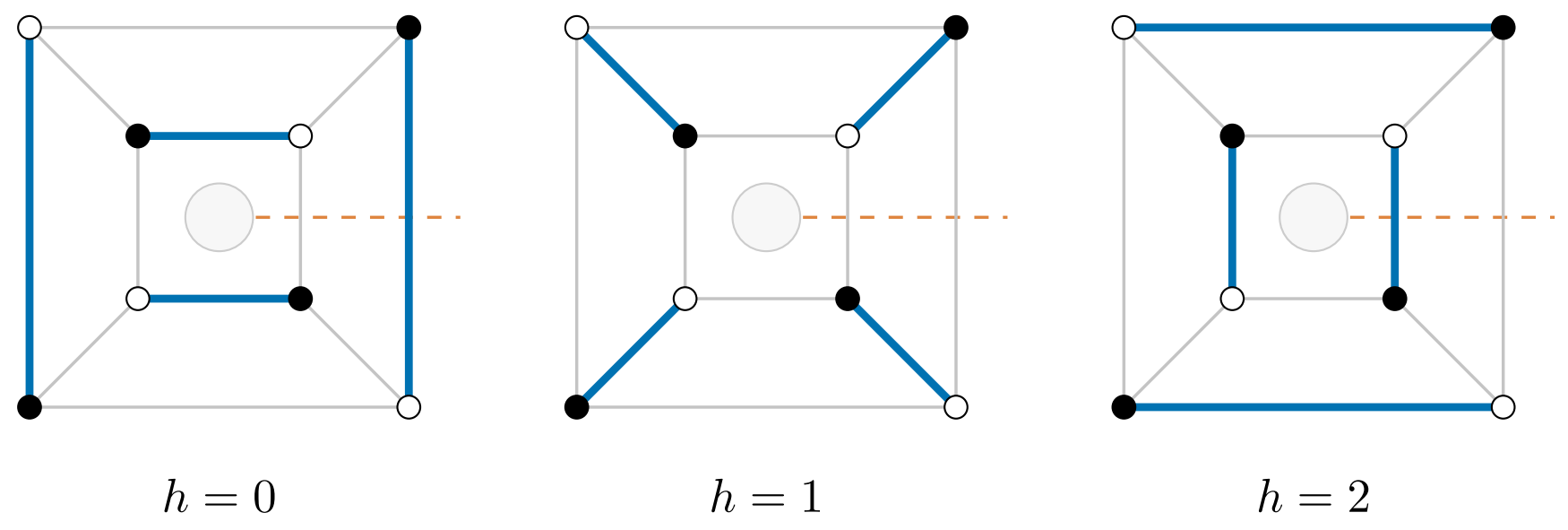

\centering\MatchingFigure
\caption{Three matchings of the ladder graph. The dashed radial arc marks the homology cut. The first and last are the unique height extrema; the middle is one of seven matchings of height $1$.}\label{fig:matchings}
\end{figure}

For the cut shown in Figure~\ref{fig:matchings}, this gives
\[
h(M)=1+[\text{inner cut edge}\in M]
       -[\text{outer cut edge}\in M],
\]
where brackets denote the indicator function.

We will return to this example as an edge-coloring
problem and as a directed network.

\subsection{Monomial characters and stretching}

We now recall a positivity conjecture of Stembridge.
For an interval $J\subseteq[N]$, let $S_J$ be the
subgroup of $S_N$ that permutes the elements of $J$
and fixes everything else. Set
\[
X_J=\sum_{\sigma\in S_J}\sigma.
\]
For a partition $\mu\vdash N$, let $\phi^\mu$ be
the virtual character whose Frobenius characteristic
is the monomial symmetric function $m_\mu$.
Stembridge conjectured that
\begin{equation}\label{eq:Stembridge}
\phi^\mu(X_{J_b}\cdots X_{J_1})\ge0
\end{equation}
for any choice of intervals $J_1,\ldots,J_b$
\cite[Conjecture~5.2]{Stembridge92}.
We write $\mathcal S(\mu)$ for this conjecture
with $\mu$ fixed.

It will be useful to package these character
evaluations into a single symmetric function:
\[
\Phi(U)=\sum_{\mu\vdash N}\phi^\mu(U)h_\mu.
\]
In Section~\ref{sec:schur}, we will see how to
construct $\Phi(U)$ by closing permutation diagrams.

There are three related versions of Stembridge's
conjecture that are worth keeping separate.
The interval-product statement $\mathcal S(\mu)$ is
Conjecture~5.2 of \cite{Stembridge92}.
Wolfgang's reduction, together with the star test
(Theorem~\ref{thm:reduction} and Lemma~\ref{lem:stars}),
shows that it is equivalent to coefficientwise
nonnegativity on all compatible acyclic path matrices
\cite[Conjecture~6.3]{Stembridge92}. Either statement implies
nonnegativity on all real totally nonnegative matrices
\cite[Conjecture~2.1]{Stembridge92}.

Tang's recent preprint \cite{Tang} claims to prove
this last statement for every shape $\mu$.
The result is about numerical nonnegativity;
the stronger interval-product and coefficientwise
network statements are not addressed there.
Our immanant theorems concern precisely these
stronger statements.

\begin{figure}[htbp]
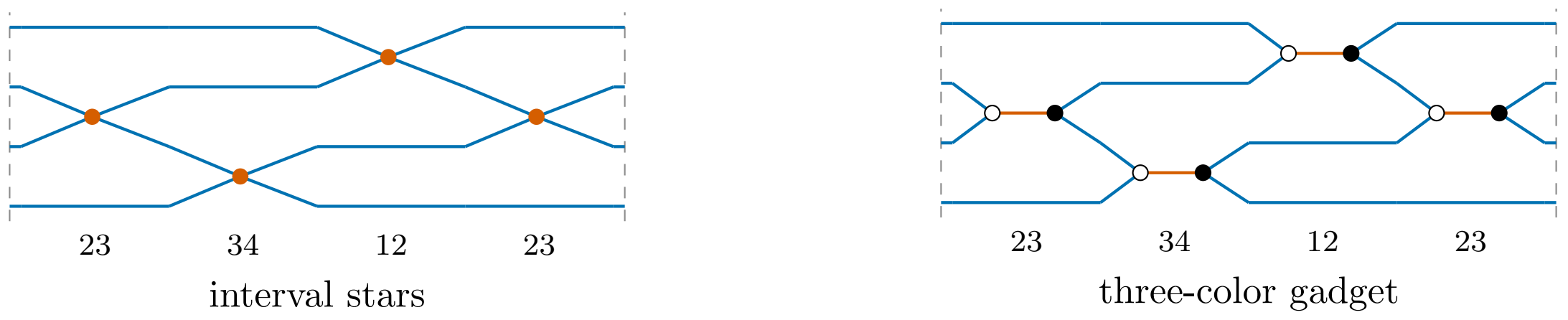

\centering
\begin{minipage}{.48\linewidth}\centering
\IntroWordFigure{0}\\[2pt]
\small interval stars
\end{minipage}\hfill
\begin{minipage}{.48\linewidth}\centering
\IntroWordFigure{1}\\[2pt]
\small three-color gadget
\end{minipage}
\caption{Two ways to draw the interval product from
Example~\ref{ex:word-intro}. The dashed sides are
identified, giving a cylinder. On the left, each
crossing allows two neighboring strands to be permuted.
On the right, each crossing is replaced by a white
and a black vertex joined by an edge. We will see
how this turns questions about permutation characters
into questions about edge colorings.}\label{fig:intro-bridge}
\end{figure}

\begin{example}\label{ex:word-intro}
For four positions, consider
\[
U=X_{[2,3]}X_{[1,2]}X_{[3,4]}X_{[2,3]}.
\]
A direct calculation gives
\begin{equation}\label{eq:word-intro}
\Phi(U)=8h_4+4h_3h_1+4h_2^2,
\end{equation}
so all monomial-character evaluations are nonnegative,
as Stembridge conjectured. Coincidentally,
$U=C'_{3412}(1)$, the Kazhdan--Lusztig element
of the smallest singular permutation.
Thus this example also illustrates Haiman's conjecture
\cite{Haiman}, which predicts
$\phi^\mu(C'_w(1))\ge0$ for all $w$ and $\mu$.

Figure~\ref{fig:intro-bridge} shows the interval product
and its three-color gadget; we return to the coloring
calculation in Example~\ref{ex:3412}.
\end{example}

There is a useful connection between Stembridge's conjecture
and colorings of bipartite graphs on the cylinder.
Superimposing $r$ perfect matchings gives an $r$-regular
bipartite multigraph $H$. Conversely, a proper
$r$-edge-coloring of $H$ decomposes it into $r$ ordered
perfect matchings $M_1,\ldots,M_r$.
We record their heights in the polynomial
\[
Z_H(x_1,\ldots,x_r)
=\sum_{(M_1,\ldots,M_r)}
x_1^{h(M_1)}\cdots x_r^{h(M_r)},
\]
where the sum is over all such decompositions; see
Section~\ref{sec:colorings} for details.
We say that $Z_H$ is \emph{elementary positive} if it
belongs to $\NN[e_1,\ldots,e_r]$, where $e_i$ are
the elementary symmetric functions.

Our first result shows that this is another way
to formulate Stembridge's conjecture.

\begin{theorem}\label{thm:dictionary-intro}
For each $r$, elementary positivity of $Z_H$ for every
$r$-regular annular bipartite multigraph is equivalent
to $\mathcal S(\mu)$ for every $\mu$ with $\mu_1\le r$.
It is also equivalent to elementary positivity of all
$r$-matching superposition polynomials
(Section~\ref{sec:annular}).
\end{theorem}

We establish positivity for up to four colors, extending
the two-column case of Stembridge's conjecture proved
by Clearman, Shelton, and Skandera \cite{CSS}.

\begin{theorem}\label{thm:fourcolumns-intro}
For $r\le4$, every annular $r$-regular bipartite
multigraph has $Z_H\in\NN[e_1,\ldots,e_r]$.
Consequently, Stembridge's conjecture holds
whenever $\mu_1\le4$.
\end{theorem}

We also give a way to obtain new cases of the
conjecture from known ones, by stretching partitions.

\begin{theorem}\label{thm:stretching-intro}
If $\mathcal S(\mu)$ holds, then $\mathcal S(n\mu)$
holds for every positive integer $n$, where
$n\mu=(n\mu_1,n\mu_2,\ldots)$.
\end{theorem}

Combining stretching with our four-color theorem and
McNulty's result for hooks \cite{McNulty}, we obtain
several infinite families.

\begin{corollary}\label{cor:families-intro}
Stembridge's conjecture holds for all partitions whose
parts lie in $\{n,2n,3n,4n\}$, and for all stretched
hooks $(an,n^b)$, where $n,a\ge1$ and $b\ge0$.
\end{corollary}

The corresponding results for path matrices, totally
nonnegative matrices, and Jacobi--Trudi matrices follow
as well; see Corollary~\ref{cor:matrix-consequences}.

Stembridge's conjecture was previously known for
rectangles \cite{Stembridge92}, shapes with at most
two columns \cite{CSS}, and hooks \cite{McNulty}.
Our results extend these to new infinite families.
Some individual small shapes outside these families
had already been verified computationally
\cite{Haiman,McNultyCheck}, but our results are
uniform and give coefficientwise positivity.

\subsection{Guide to the proofs}

Section~\ref{sec:schur} collects the path-matrix and
character identities used throughout the paper.
Sections~\ref{sec:colorings} and~\ref{sec:fourcolors}
use annular colorings to prove Stembridge's conjecture
for partitions with at most four columns.
Section~\ref{sec:stretching} extends these results by
stretching partitions. Sections~\ref{sec:networks}
and~\ref{sec:annular} prove the cylinder theorems,
and Section~\ref{sec:annular} completes the equivalence
with Stembridge's conjecture.

\begin{figure}[htbp]
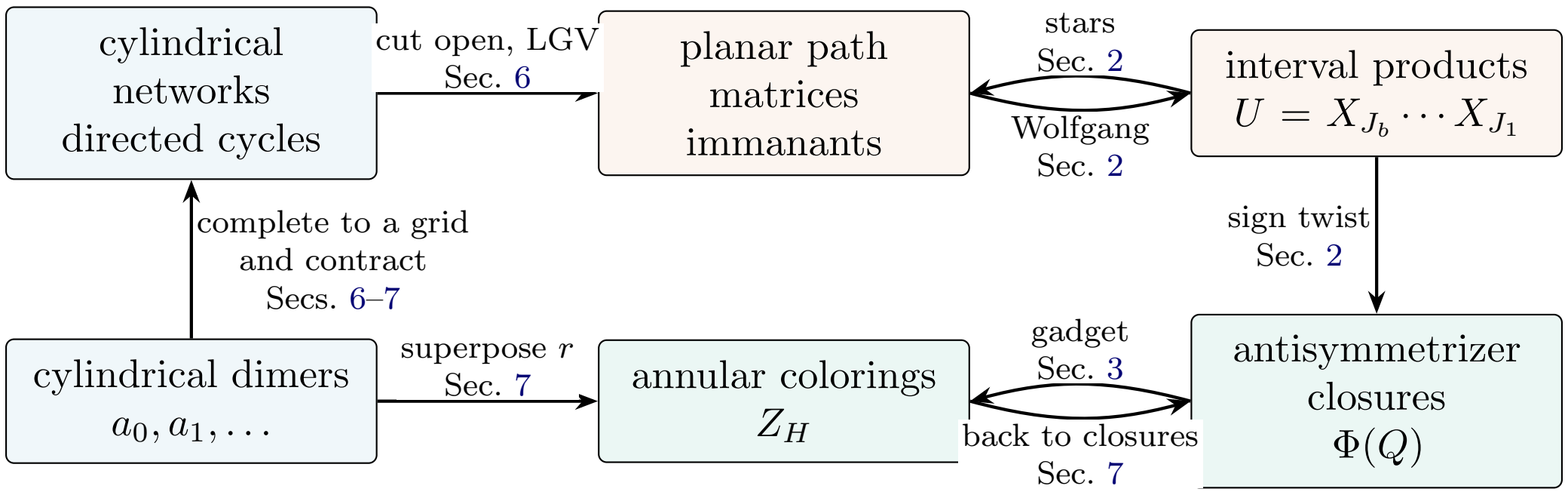

\centering\RoadmapFigure
\caption{The constructions linking cylindrical dimers,
planar path matrices, and annular colorings.}
\label{fig:roadmap}
\end{figure}

\subsection*{Acknowledgements}
The author is grateful to Chris Fraser for feedback on the draft. 

\subsection*{AI statement}
The author received substantial assistance from ChatGPT
(OpenAI) and Claude (Anthropic) in developing the results
and preparing the manuscript. The author takes full responsibility for the content.

%\Needspace{15\baselineskip}
\section{Closures and planar path matrices}\label{sec:schur}

We use the symmetric-function conventions of
\cite[Chapter~7]{Stanley}. 

\subsection{Toeplitz minors as Schur coefficients}

The first step is to express Toeplitz minors as
Schur coefficients.

\begin{lemma}\label{lem:jt}
For every partition $\lambda$ with at most $r$ parts,
\begin{equation}\label{eq:alternant}
[s_\lambda(\xx)]\prod_{i=1}^r A(x_i)
=\det(a_{\lambda_i-i+j})_{i,j=1}^r.
\end{equation}
In particular, if these Schur coefficients are
coefficientwise nonnegative for every $r$, then
every minor of $T(a)$ is coefficientwise nonnegative.
\end{lemma}

\begin{proof}
We use the dual Cauchy identity. First assume $a_0$ is
invertible, and introduce a formal alphabet $\mathbf y$
by setting $e_k(\mathbf y)=a_k/a_0$. Then
\[
\prod_{i=1}^r A(x_i)
=a_0^r\sum_\nu s_\nu(\xx)s_{\nu'}(\mathbf y),
\]
where the sum is over partitions with at most $r$ parts.
By the dual Jacobi--Trudi identity,
\[
a_0^r s_{\lambda'}(\mathbf y)
=\det(a_{\lambda_i-i+j})_{i,j=1}^r,
\]
which proves \eqref{eq:alternant}.

The same argument handles arbitrary minors. Choose rows
$i_1<\cdots<i_r$ and columns $j_1<\cdots<j_r$ of
$T(a)$. Shifting all indices by the same integer if
necessary, define partitions
\[
\mu_k=i_{r+1-k}-(r-k),\qquad
\lambda_k=j_{r+1-k}-(r-k).
\]
The skew dual Jacobi--Trudi identity gives
\[
\det(a_{j_\ell-i_k})_{k,\ell=1}^r
=a_0^r s_{\lambda'/\mu'}(\mathbf y).
\]
By Littlewood--Richardson, the skew Schur function
$s_{\lambda'/\mu'}$ is a nonnegative integer combination
of ordinary Schur functions. Thus every Toeplitz minor
is a nonnegative integer combination of the Schur
coefficients in \eqref{eq:alternant}.

Finally, these are polynomial identities in the $a_k$,
so they remain valid without assuming $a_0$ invertible.
\end{proof}

\subsection{Closures and the sign twist}

Let $S_N$ act on words of length $N$ by permuting
their positions. For an interval $J\subseteq[N]$, set
\[
X_J=\sum_{\sigma\in S_J}\sigma,\qquad
P_J=\frac{1}{|J|!}\sum_{\sigma\in S_J}
\sgn(\sigma)\sigma.
\]
The first operator symmetrizes the positions in $J$,
while the second antisymmetrizes them. For example,
\[
P_{[1,2]}(abc)=\frac12(abc-bac),\qquad
P_{[1,2]}(aac)=0.
\]
More generally, $P_J$ is a projection onto the
subspace of linear combinations of words that change
sign when two positions in $J$ are exchanged.
We will use antisymmetrizers because they kill
repeated colors, naturally enforcing the distinct-color
condition in our graphs. A sign twist will relate
them to the symmetrizers appearing in Stembridge's
conjecture.

For $g=\sum_{\sigma\in S_N}g_\sigma\sigma
\in\mathbb Q[S_N]$, define its \emph{color closure} by
\begin{equation}\label{eq:Phi}
\Phi(g)=\sum_{\sigma\in S_N}
g_\sigma p_{\type(\sigma)}.
\end{equation}
For example,
\[
\Phi(P_{[1,2]})
=\frac12(p_1^2-p_2)=e_2=s_{11}.
\]
There is a simple picture behind this definition.
Draw a permutation diagram and identify its top
and bottom, closing the strands into cycles.
A coloring that is constant along each strand
must assign one color to each cycle. A cycle of
length $j$ contributes $p_j=\sum_i x_i^j$,
giving exactly \eqref{eq:Phi}.

\begin{figure}[htbp]
\centering
\begingroup
\newcommand{\ClosureRouting}[2]{%
  \begingroup
  \def\cOne{cA}\def\cTwo{cB}\def\cThree{cC}%
  \def\opname{straight/straight}
  \def\outterm{p_1^3}\def\qsign{+}%
  \ifnum#1=0
    \ifnum#2=1
      \def\cOne{cA}\def\cTwo{cA}\def\cThree{cB}%
      \def\opname{straight/cross}
      \def\outterm{p_2p_1}\def\qsign{-}%
    \fi
  \else
    \ifnum#2=0
      \def\cOne{cB}\def\cTwo{cA}\def\cThree{cA}%
      \def\opname{cross/straight}
      \def\outterm{p_2p_1}\def\qsign{-}%
    \else
      \def\cOne{cA}\def\cTwo{cA}\def\cThree{cA}%
      \def\opname{cross/cross}
      \def\outterm{p_3}\def\qsign{+}%
    \fi
  \fi
  \begin{tikzpicture}[
    x=.77cm,y=.77cm,
    baseline=(current bounding box.center)]
    \fill[black!3,rounded corners=2pt]
      (1.73,1.92) rectangle (3.26,2.51);
    \draw[black!21,densely dashed,rounded corners=2pt]
      (1.73,1.92) rectangle (3.26,2.51);
    \fill[black!3,rounded corners=2pt]
      (.73,.86) rectangle (2.26,1.45);
    \draw[black!21,densely dashed,rounded corners=2pt]
      (.73,.86) rectangle (2.26,1.45);
    \draw[black!35,densely dashed,line width=.55pt]
      (.55,3.04)--(3.45,3.04);
    \draw[black!35,densely dashed,line width=.55pt]
      (.55,.30)--(3.45,.30);
    \draw[cB!75,thin,-{Stealth[length=1.5mm]}]
      (.30,.39)--(.30,2.96);
    \node[font=\scriptsize,anchor=west,color=black!65]
      at (3.30,2.20) {$[2,3]$};
    \node[font=\scriptsize,anchor=east,color=black!65]
      at (.66,1.16) {$[1,2]$};
    \foreach \i/\c in
      {1/\cOne,2/\cTwo,3/\cThree}{
       \pgfmathtruncatemacro{\midp}{
         \i + (#1)*((\i==2)-(\i==3))}
       \pgfmathtruncatemacro{\endp}{
         \midp + (#2)*((\midp==1)-(\midp==2))}
       \draw[white,line width=3.6pt,line cap=round]
         (\i,3.04)--(\i,2.49)
         --(\midp,1.94)--(\midp,1.43)
         --(\endp,.88)--(\endp,.30);
       \draw[draw=\c,line width=1.65pt,line cap=round]
         (\i,3.04)--(\i,2.49)
         --(\midp,1.94)--(\midp,1.43)
         --(\endp,.88)--(\endp,.30);
       \fill[\c] (\i,3.04) circle (1.5pt);
       \fill[\c] (\endp,.30) circle (1.5pt);
    }
    \foreach \i in {1,2,3}{
      \node[font=\scriptsize,color=black!70]
        at (\i,3.33) {\i};
    }
    \node[font=\scriptsize] at (2,3.79) {\opname};
    \node[font=\small] at (2,-.13) {$\outterm$};
    \node[font=\scriptsize,color=black!70]
      at (2,-.53) {$Q:\;\qsign$};
  \end{tikzpicture}%
  \endgroup
}
\begin{minipage}{.24\linewidth}
\centering\ClosureRouting{0}{0}
\end{minipage}\hfill%
\begin{minipage}{.24\linewidth}
\centering\ClosureRouting{0}{1}
\end{minipage}\hfill%
\begin{minipage}{.24\linewidth}
\centering\ClosureRouting{1}{0}
\end{minipage}\hfill%
\begin{minipage}{.24\linewidth}
\centering\ClosureRouting{1}{1}
\end{minipage}
\endgroup
\caption{Closing a product of two interval operators.
At the upper interval $[2,3]$ and the lower interval
$[1,2]$, strands can go straight or cross.
The dashed top and bottom boundaries are identified;
the colors distinguish the resulting cycles.
Each routing contributes the indicated power sum
to $\Phi(U)$. For $\Phi(Q)$, use the displayed signs
and divide by $4$.}
\label{fig:closure-example}
\end{figure}
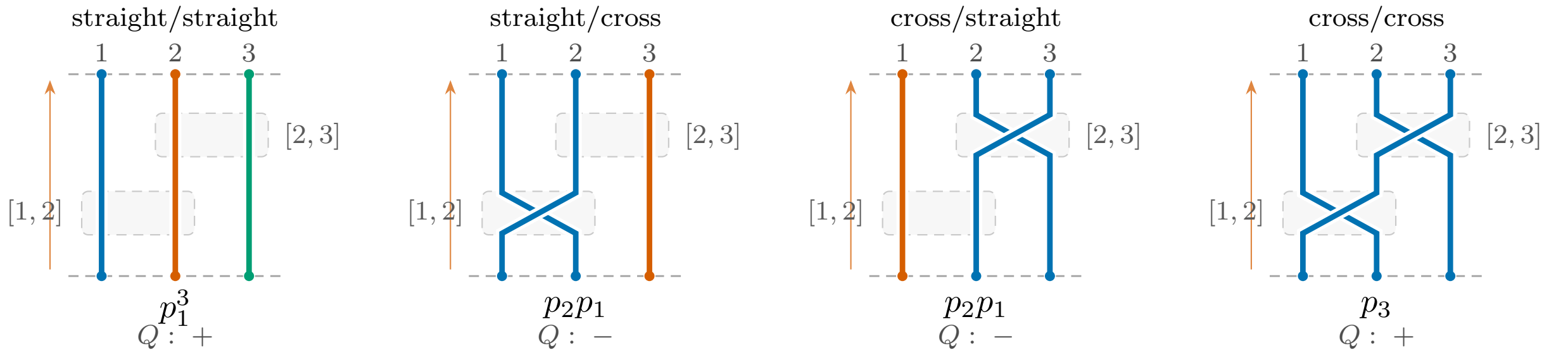

For a single permutation $\sigma$, closing its diagram
joins the strands into cycles, so
$\Phi(\sigma)=p_{\type(\sigma)}$.
For a product of interval operators, each interval
allows several ways to route the strands. Consider
\[
U=X_{[1,2]}X_{[2,3]},\qquad
Q=P_{[1,2]}P_{[2,3]}.
\]
The four routings in Figure~\ref{fig:closure-example}
give cycle types $(1,1,1)$, $(2,1)$, $(2,1)$,
and $(3)$, respectively. Thus
\[
\Phi(U)=p_1^3+2p_2p_1+p_3,\qquad
\Phi(Q)=\frac14(p_1^3-2p_2p_1+p_3).
\]
For $Q$, each crossing contributes a minus sign,
and the two antisymmetrizer normalizations contribute
the factor $1/4$.

Equivalently, consider the vector space with basis
all words of length $N$ over $[r]$. Let $D_x$ be
the diagonal operator defined by
\[
D_x(w_1\cdots w_N)
=x_{w_1}\cdots x_{w_N}\,w_1\cdots w_N.
\]
Since $g$ also acts linearly on this space, we have
\[
\Phi(g)(x_1,\ldots,x_r)=\tr(D_xg),
\]
where $D_xg$ is the composition of the two operators.
Indeed, a word contributes to the trace precisely
when it is constant on each cycle of the permutation.
Notice that $D_x$ commutes with permutations of
positions, and hence with all $P_J$. 

Let us now connect these closures to characters.
Recall that the Frobenius characteristic is
\[
\ch(f)=\frac1{N!}\sum_{\sigma\in S_N}
f(\sigma)p_{\type(\sigma)}.
\]
We have $\ch(\chi^\lambda)=s_\lambda$ and
$\ch(\phi^\mu)=m_\mu$. Since the Schur functions
are orthonormal and the bases $h_\mu,m_\mu$
are dual under the Hall inner product,
\begin{equation}\label{eq:charclosure}
[s_\lambda]\Phi(g)=\chi^\lambda(g),\qquad
[h_\mu]\Phi(g)=\phi^\mu(g).
\end{equation}

The key observation is that multiplying each
permutation by its sign induces the involution
$\omega$ on symmetric functions, with
$\omega(h_j)=e_j$ and $\omega(s_\lambda)=s_{\lambda'}$.
Thus, if
\[
U=X_{J_b}\cdots X_{J_1},\qquad
Q=P_{J_b}\cdots P_{J_1},
\]
then, writing $\kappa=\prod_t|J_t|!$, we obtain
\begin{equation}\label{eq:sign-twist}
\Phi(Q)=\frac{\omega\Phi(U)}{\kappa},
\qquad
[e_\mu]\Phi(Q)=\frac{\phi^\mu(U)}{\kappa}.
\end{equation}
In $r$ variables, the products $e_\mu$ with
$\mu_1\le r$ form a basis in each degree.
We therefore get the following dictionary:
\begin{center}
\begin{tabular}{lll}
\toprule
Coefficient & Character evaluation & Positivity\\
\midrule
$[s_\lambda]\Phi(Q)$ &
$\chi^{\lambda'}(U)/\kappa$ &
Greene's theorem\\
$[e_\mu]\Phi(Q)$ &
$\phi^\mu(U)/\kappa$ &
Stembridge's conjecture\\
\bottomrule
\end{tabular}
\end{center}
In particular, Schur positivity of these closures
follows from Greene's theorem \cite{GJ,Greene},
while elementary positivity is precisely Stembridge's
interval-product conjecture
\cite[Conjecture~5.2]{Stembridge92}.
Tang's recent preprint \cite{Tang} claims the weaker
numerical nonnegativity of monomial immanants on
totally nonnegative matrices, but does not establish
the coefficientwise statement.

\subsection{The path-matrix reduction}

We now recall how products of interval sums arise
from planar networks. For a class function $f$ on
$S_N$ and an $N\times N$ matrix $R$, define the
\emph{immanant}
\[
\Imm_f(R)=\sum_{\sigma\in S_N}
f(\sigma)\prod_{i=1}^N R_{i,\sigma(i)}.
\]
Let $G$ be an acyclic directed graph with marked
vertices $u_1,\ldots,u_N$ and $v_1,\ldots,v_N$,
and assign an independent variable $y_e$ to each edge.
Its \emph{path matrix} is
\[
R_{ij}=\sum_{\gamma:u_i\to v_j}
\prod_{e\in\gamma}y_e.
\]
We call the endpoints \emph{compatible} if any two
paths $u_i\to v_j$ and $u_k\to v_\ell$ must meet
whenever $i<k$ and $j>\ell$. Shared endpoints
count as intersections. This holds, for example,
when the sources and sinks are ordered along
opposite sides of a disk. Throughout, a disk path matrix has its sources and sinks
on opposite boundary arcs, numbered in compatible order.

The key point of Wolfgang's reduction is that
positivity on products of interval sums already
implies coefficientwise positivity for these
path matrices.

\begin{theorem}[{Wolfgang's reduction
\cite[Section~2.6]{Wolfgang}}]\label{thm:reduction}
If a class function $f$ is nonnegative on every
product of interval sums in $S_N$, then
$\Imm_f(R)$ has nonnegative coefficients for
every compatible acyclic path matrix $R$.
\end{theorem}

Indeed, Wolfgang's reduction expresses each
edge-monomial coefficient of $\Imm_f(R)$ as a
nonnegative combination of evaluations of $f$
on products of interval sums. Conversely, these
evaluations can themselves be recovered as
coefficients of immanants of planar networks.

\begin{lemma}[The star test]\label{lem:stars}
Suppose that $\Imm_f(R)$ has nonnegative coefficients
for every acyclic planar disk network with compatible
endpoints. Then
\[
f(X_{J_b}\cdots X_{J_1})\ge0
\]
for every choice of intervals $J_1,\ldots,J_b
\subseteq[N]$.
\end{lemma}

\begin{figure}[htbp]
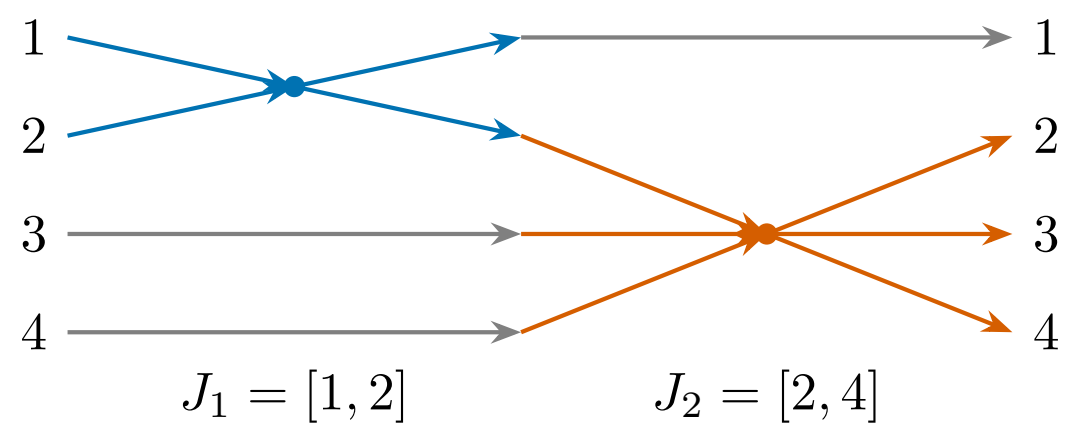

\centering\StarFigure
\caption{Realizing an interval product by a planar
network. At each common vertex, incoming and outgoing
edges can be paired arbitrarily, giving a permutation
of the corresponding interval. The two stages realize
$X_{[2,4]}X_{[1,2]}$.}
\label{fig:stars}
\end{figure}

\begin{proof}
Let $U=X_{J_b}\cdots X_{J_1}$. Start with $N$
parallel directed strands, and for each interval
$J_t$, bring the corresponding strands together
at a single vertex, as in Figure~\ref{fig:stars}.
All other strands pass straight through. Since
the intervals are consecutive, the resulting
network is planar. Give every edge an independent
variable, and let $R$ be its path matrix.

Consider the coefficient of $\prod_e y_e$ in
$\Imm_f(R)$. To contribute, a collection of paths
must use every edge exactly once. At each common
vertex, this amounts to choosing a bijection
between the incoming and outgoing edges, or
equivalently a permutation in $S_{J_t}$.
Thus these path collections correspond exactly
to the terms in $U$, and
\begin{equation}\label{eq:star-test}
[\textstyle\prod_e y_e]\Imm_f(R)=f(U).
\end{equation}
The assumed coefficientwise nonnegativity therefore
gives $f(U)\ge0$, as required.
\end{proof}

Greene proved that interval sums $X_J$ act by
nonnegative matrices in Young's seminormal basis
\cite{GJ,Greene}; see also
\cite[Theorem~5.1]{Stembridge92}.
As observed by Haiman, the same holds in the
Kazhdan--Lusztig basis
\cite[footnote to Theorem~5.1]{Stembridge92}.
Their products therefore have nonnegative traces,
so $\chi^\lambda(U)\ge0$ for every product $U$
of interval sums.
 Together with Wolfgang's reduction,
this gives the following classical result.

\begin{theorem}[{Wolfgang
\cite[Corollary~2.6.5]{Wolfgang}}]\label{thm:wolfgang}
Every irreducible immanant of a compatible acyclic
path matrix has nonnegative coefficients.
\end{theorem}

We also get a consequence directly for our closures.
If $U=X_{J_b}\cdots X_{J_1}$ and
$Q=P_{J_b}\cdots P_{J_1}$, then the
involution $\omega$ gives
\[
[s_\lambda]\Phi(Q)
=\frac{\chi^{\lambda'}(U)}{\kappa}\ge0,
\qquad \kappa=\prod_t |J_t|!.
\]
Thus we have:

\begin{corollary}\label{cor:closure}
The closure of any product of antisymmetrizers
on consecutive intervals is Schur positive
over $\mathbb Q$.
\end{corollary}

For example, the closure from
Figure~\ref{fig:closure-example} is
\[
\Phi(P_{[1,2]}P_{[2,3]})
=\frac14s_{21}+s_{111}.
\]

\subsection{From immanants to symmetric functions}

We can use immanant positivity to study symmetric
functions evaluated at matrices. The connection is
given by the following identity.

\begin{lemma}\label{lem:frobenius-immanant}
Let $f$ be a class function on $S_N$, and let $R$
be an $m\times m$ matrix. Then
\begin{equation}\label{eq:SchurImm}
(\ch f)(R)=\frac1{N!}
\sum_{i_1,\ldots,i_N=1}^m \Imm_f(R[i,i]),
\end{equation}
where $i=(i_1,\ldots,i_N)$, the matrix $R[i,i]$
allows repeated rows and columns, and the
symmetric function on the left is evaluated
by $p_j(R)=\tr(R^j)$.
\end{lemma}

\begin{proof}
Expanding the immanants and interchanging the sums gives
\[
\begin{aligned}
\sum_{i_1,\ldots,i_N=1}^m \Imm_f(R[i,i])
&=\sum_{i_1,\ldots,i_N=1}^m
  \sum_{\sigma\in S_N}f(\sigma)
  \prod_{a=1}^N R_{i_a,i_{\sigma(a)}}\\
&=\sum_{\sigma\in S_N}f(\sigma)
  \sum_{i_1,\ldots,i_N=1}^m
  \prod_{a=1}^N R_{i_a,i_{\sigma(a)}}.
\end{aligned}
\]
For each cycle $(a_1\,\cdots\,a_j)$ of $\sigma$,
the corresponding indices contribute
\[
\sum_{i_{a_1},\ldots,i_{a_j}=1}^m
R_{i_{a_1},i_{a_2}}\cdots R_{i_{a_j},i_{a_1}}
=\tr(R^j).
\]
Since different cycles involve disjoint indices,
their contributions multiply. Thus the inner sum
is $p_{\type(\sigma)}(R)$, and dividing by $N!$
gives 
\[
\frac1{N!}\sum_{\sigma\in S_N}
f(\sigma)p_{\type(\sigma)}(R)
=(\ch f)(R).
\]
\end{proof}

\begin{corollary}\label{cor:pathschur}
For every planar acyclic disk path matrix $R$,
the symmetric functions $s_\lambda(R)$ have
nonnegative coefficients in the edge variables.
Moreover, for any fixed partition $\mu$,
the assertion $\mathcal S(\mu)$ is equivalent
to coefficientwise nonnegativity of $m_\mu(R)$
for all such matrices $R$.
\end{corollary}

\begin{proof}
By Lemma~\ref{lem:frobenius-immanant},
\[
s_\lambda(R)=\frac1{N!}
\sum_{i_1,\ldots,i_N}
\Imm_{\chi^\lambda}(R[i,i]),
\qquad N=|\lambda|.
\]
Each $R[i,i]$ is again a compatible path matrix,
even if some indices repeat. Indeed, we may sort
the indices, since simultaneously permuting
rows and columns preserves an immanant associated
with a class function. Repeated sources and sinks
cause no problem: paths sharing an endpoint
automatically intersect.

By Wolfgang's theorem
(Theorem~\ref{thm:wolfgang}), every summand has
nonnegative coefficients, proving the first
assertion. If $\mathcal S(\mu)$ holds, the same
argument applies to $m_\mu(R)$, using Wolfgang's
reduction (Theorem~\ref{thm:reduction}) for the
monomial character $\phi^\mu$.

For the converse, suppose $m_\mu(R)$ has
nonnegative coefficients for every planar
disk path matrix $R$. Take such a matrix of
size $N=|\mu|$ and put
\[
S=\operatorname{diag}(t_1,\ldots,t_N)R.
\]
This is again a planar path matrix, obtained
by giving source $i$ the additional weight $t_i$.
Thus $m_\mu(S)$ has nonnegative coefficients.
In \eqref{eq:SchurImm} for $S$, the coefficient of
$t_1\cdots t_N$ comes from the $N!$ tuples using each
index once. Each gives the same immanant of $R$. Thus
\[
\Imm_{\phi^\mu}(R)
=[t_1\cdots t_N]m_\mu(S),
\]
so the immanant also has nonnegative coefficients.
The star test (Lemma~\ref{lem:stars}) now gives
$\mathcal S(\mu)$.
\end{proof}

\section{Annular colorings and monomial characters}\label{sec:colorings}

We now turn to edge colorings of bipartite graphs on
the annulus. We will relate their coloring polynomials
to the closures introduced above, and prove elementary
positivity for two and three colors. The four-color
case is treated in Section~\ref{sec:fourcolors}.

Let $H$ be an $r$-regular bipartite multigraph
embedded in an annulus. A proper $r$-edge-coloring
assigns the colors $1,\ldots,r$ to its edges,
with every color appearing exactly once at each
vertex. Thus each color class is a perfect matching,
and a coloring amounts to an ordered decomposition
of $H$ into $r$ matchings $M_1,\ldots,M_r$.
Define
\begin{equation}\label{eq:ZH}
Z_H(x_1,\ldots,x_r)
=\sum_c\prod_{i=1}^r x_i^{h(M_i)},
\end{equation}
where the sum runs over all proper colorings and
the heights are normalized to have minimum zero.

The polynomial $Z_H$ is symmetric, since permuting
the colors permutes the colorings. It is also
homogeneous: the sum of the heights of the color
classes is fixed, since together they use every
edge of $H$ exactly once. Finally, such colorings
always exist. Indeed, every regular bipartite
multigraph has a perfect matching, and removing
one leaves an $(r-1)$-regular bipartite multigraph.

\subsection{A positive local rule}

We want to interpret antisymmetrizer closures as
proper edge colorings of bipartite graphs. The
antisymmetrizers $P_J$ already have one useful
property: they kill words with repeated letters
in $J$, enforcing the distinct-color condition.
The problem is that their other terms come with
signs. We now show that a simple change of basis
removes these signs without changing the closure.

For a word $w=w_1\cdots w_N$ over $[r]$, let
$\operatorname{inv}(w)$ count the pairs $i<j$
with $w_i>w_j$, and define the diagonal operator
\[
D(w)=(-1)^{\operatorname{inv}(w)}w.
\]

\begin{lemma}\label{lem:word-gauge}
The operator $DP_JD^{-1}$ kills any word with
repeated letters in positions $J$. Otherwise,
it sends the word to the uniform positive average
of all reorderings of those letters.
\end{lemma}

\begin{proof}
If two letters in $J$ are equal, their terms
cancel in the antisymmetrizer. Otherwise,
reordering the letters in $J$ does not change
the inversions involving positions outside $J$,
since $J$ is consecutive. The change in inversion
parity is therefore the sign of the reordering,
which cancels the sign in $P_J$.
\end{proof}

For example,
\[
\begin{aligned}
P_{[1,2]}(21)&=\tfrac12(21-12),&
DP_{[1,2]}D^{-1}(21)&=\tfrac12(21+12),\\
P_{[1,2]}(11)&=0,&
DP_{[1,2]}D^{-1}(11)&=0.
\end{aligned}
\]

The important point is that this change of basis
does not affect the closure. We obtain the following
positive formula.

\begin{lemma}[Positive closure rule]
\label{lem:positive-closure}
Let $Q=P_{J_b}\cdots P_{J_1}$ be a product of
interval antisymmetrizers. Then
\[
\Phi(Q)(x_1,\ldots,x_r)
=\frac{1}{\prod_{t=1}^b|J_t|!}
\sum_{w^{(0)},\ldots,w^{(b)}}
x_{w^{(0)}_1}\cdots x_{w^{(0)}_N},
\]
where the sum is over sequences of words over $[r]$
satisfying $w^{(b)}=w^{(0)}$ and the following
rule at each step $t$: the letters in positions
$J_t$ of $w^{(t-1)}$ must be distinct, and
$w^{(t)}$ is obtained by reordering those letters,
leaving all other positions unchanged.
\end{lemma}

\begin{proof}
Recall that $\Phi(Q)=\tr(D_xQ)$. Since $D$ and
$D_x$ are both diagonal, they commute, and
conjugation preserves the trace. Thus
\[
\Phi(Q)=\tr\bigl(D_x(DP_{J_b}D^{-1})\cdots
(DP_{J_1}D^{-1})\bigr).
\]
By Lemma~\ref{lem:word-gauge}, each factor
allows exactly the stated reorderings, with
weight $1/|J_t|!$. The trace requires the
final word to equal the initial one, while
$D_x$ supplies its monomial weight.
\end{proof}

For example, take $Q=P_{[1,2]}P_{[2,3]}$
with two colors. The only closed sequences are
\[
121\to121\to121,\qquad
212\to212\to212,
\]
each with weight $1/4$. Thus
\[
\Phi(Q)(x_1,x_2)
=\frac14(x_1^2x_2+x_1x_2^2),
\]
in agreement with the signed formula
$\Phi(Q)=\frac14(p_1^3-2p_1p_2+p_3)$
from Section~2.2. Notice that repeated colors
are allowed: they only need to be distinct
within each interval being antisymmetrized.

The next subsection realizes this positive local
rule as proper edge colorings of a bipartite graph.

\subsection{The bipartite gadget}

We now turn the positive closure rule into a graph
coloring problem. Recall that at an interval $J$
of size $k$, the incoming colors must be distinct,
and the same colors may leave in any order.
We can realize this rule using two vertices of
a bipartite graph.

\begin{figure}[htbp]
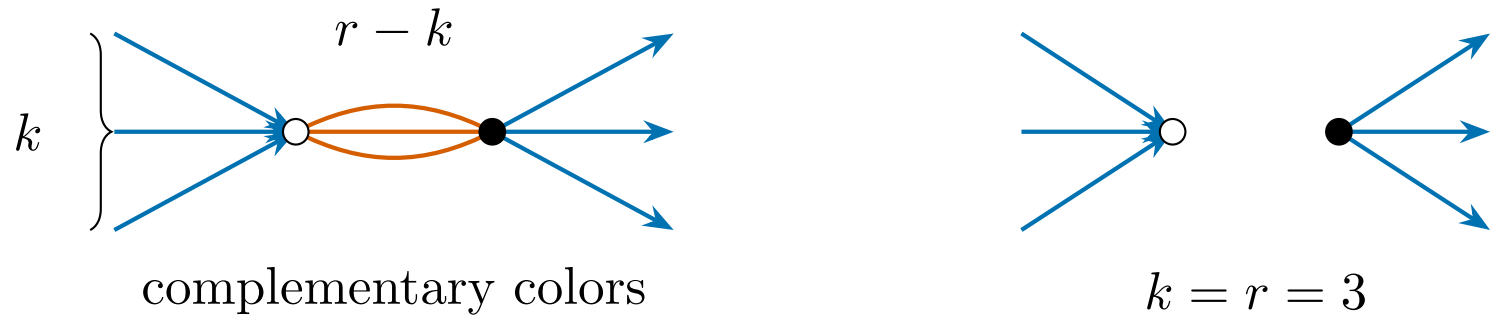

\centering\GadgetFigure
\caption{The bipartite gadget for an interval of
size $k$. The incoming strands meet a white vertex,
the outgoing strands start at a black vertex,
and the $r-k$ joining edges carry the remaining
colors. Right: the case $k=r=3$, with no joining
edges. The horizontal direction follows the
cylinder around.}
\label{fig:gadget}
\end{figure}

Draw $N$ parallel strands on the cylinder, with
the intervals $J_1,\ldots,J_b$ appearing in cyclic
order. At an interval $J$ of size $k\le r$,
bring the $k$ incoming strands together at a white
vertex, and let the $k$ outgoing strands start
at a black vertex. Join these two vertices by
$r-k$ parallel edges, as in Figure~\ref{fig:gadget}.
In a proper $r$-edge-coloring, the joining edges
must carry exactly the colors missing from the
incoming strands. The outgoing strands therefore
have the same colors as the incoming ones, in
any order. This is precisely the rule of
Lemma~\ref{lem:positive-closure}.

Since the intervals are consecutive, the
construction is planar, and every vertex has
degree $r$. A strand untouched by all intervals
closes into a circle and contributes $e_1$.
Remove these circles, say $u$ of them, and
denote the remaining bipartite graph by $H$.
Choose a cut meeting only the forward strands,
and let $f_0\ge0$ be the minimum number of
matching edges crossing it. Thus the number
of crossings of a matching is its normalized
height plus $f_0$.

\begin{lemma}[The gadget]\label{lem:gadget}
If $|J_t|\le r$ for every $t$, then
\begin{equation}\label{eq:gadget}
\Phi(P_{J_b}\cdots P_{J_1})(x_1,\ldots,x_r)
=\frac{e_1^u e_r^{f_0}}
{\prod_t |J_t|!(r-|J_t|)!}\,
Z_H(x_1,\ldots,x_r).
\end{equation}
If some $|J_t|>r$, the closure vanishes
in $r$ variables.
\end{lemma}

\begin{proof}
Consider an interval $J$ of size $k$.
By Lemma~\ref{lem:positive-closure}, its incoming
colors must be distinct, and any reordering is
allowed, with weight $1/k!$.

The gadget realizes exactly these transitions.
The $r-k$ joining edges must carry the colors
missing from the incoming strands, and hence
the outgoing strands have the same set of colors.
For each allowed transition, there are
$(r-k)!$ ways to color the joining edges.
Thus the gadget contributes a factor of
$k!(r-k)!$ relative to the positive closure rule.

Multiplying over all intervals gives the formula
with raw flux in place of height. Each untouched
strand contributes $e_1$. Finally, raw flux
equals normalized height plus $f_0$ for every
color class, giving the factor $e_r^{f_0}$.
If $k>r$, the distinct-color condition cannot
be satisfied, so the closure vanishes.
\end{proof}

\begin{lemma}[Restoring $e_r$]\label{lem:elementary-lift}
Let $F$ be a homogeneous symmetric polynomial
in $r$ variables. If
\[
F|_{e_r=1}\in\NN[e_1,\ldots,e_{r-1}],
\]
then $F\in\NN[e_1,\ldots,e_r]$.
\end{lemma}

\begin{proof}
Expand $F$ in the elementary symmetric functions.
Since $F$ is homogeneous and $\deg e_i=i$,
the exponents of $e_1,\ldots,e_{r-1}$ uniquely
determine the exponent of $e_r$.
Thus setting $e_r=1$ does not combine any terms,
and every coefficient remains unchanged.
\end{proof}

\subsection{Two and three colors}

Let us start with two colors. Since every vertex
of $H$ has degree two, $H$ is a disjoint union
of even cycles. Each cycle has exactly two proper
colorings, obtained by exchanging its colors.
Switching the coloring of a contractible cycle
does not change either matching's homological height.
Thus, after normalizing heights, each such cycle
contributes a factor of $2$, with no additional
powers of $x_1$ or $x_2$.
For an essential cycle, their fluxes differ by
one, giving a factor of $x_1+x_2=e_1$.

Since the cycles can be colored independently,
normalizing the minimum height to zero gives
\begin{equation}\label{eq:two-color}
Z_H=2^c e_1^\ell,
\end{equation}
where $c$ and $\ell$ are the numbers of
contractible and essential cycles, respectively.
If we use an inherited height normalization,
the expression may acquire an additional
nonnegative power of $e_2$.

Thus elementary positivity is immediate for
two colors. Our next result establishes it
for three colors as well.

\begin{proposition}\label{prop:threecolor}
For every $3$-regular annular bipartite
multigraph $H$, the polynomial $Z_H$
is elementary positive.
\end{proposition}

\begin{figure}[htbp]
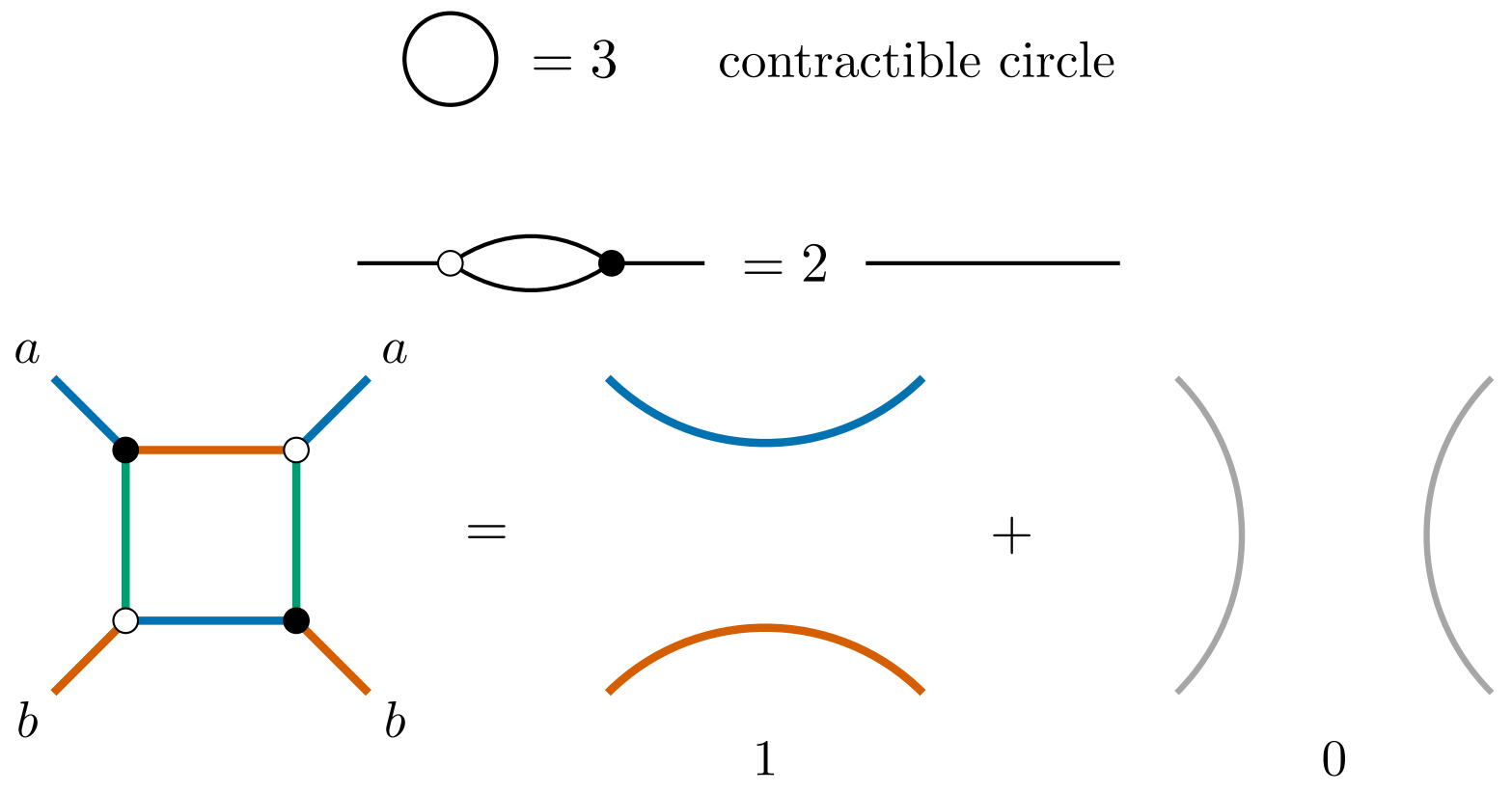

\centering\ThreeRulesFigure
\caption{The three-color rules: a contractible circle is $3$, a bigon is twice an edge, and a square is the sum of its two smoothings. For the illustrated boundary word $aabb$, the square has the unique displayed coloring; its smoothing counts are $1$ and $0$.}\label{fig:three-rules}
\end{figure}

\begin{proof}
We first work on $e_3=x_1x_2x_3=1$.
Orient every edge from black to white, and give
an edge of color $i$ weight $x_i^\delta$, where
$\delta$ is its signed crossing number with a
chosen cut. Moving the cut past a vertex changes
the weight of a coloring by
$(x_1x_2x_3)^{\pm1}=1$, since every vertex sees
all three colors. Thus we may move the cut away
from any contractible disk and apply local
relations there. 
Let $F(M)$ be the signed number of crossings of
the cut by a matching $M$, and let
$f_0=\min_M F(M)$. Then $h(M)=F(M)-f_0$.
Thus, for a coloring with matchings $M_1,M_2,M_3$,
\[
\prod_{i=1}^3 x_i^{F(M_i)}
=(x_1x_2x_3)^{f_0}
\prod_{i=1}^3 x_i^{h(M_i)}.
\]
Since $x_1x_2x_3=1$, the crossing weights
agree with those in $Z_H$.

We use three positive local rules, shown in
Figure~\ref{fig:three-rules}. These are the
rules of \cite[Figure~1 and Remark~2.1]{PWeb},
obtained from Kuperberg's spider \cite{Kuperberg}
at $q=1$. A contractible circle contributes $3$,
since it can have any of the three colors.
A bigon contributes twice a single edge:
the two external colors must agree, and the
remaining two colors can be assigned to the
bigon in either order. Finally, a square is
the sum of its two smoothings.

Let us check the square rule. Fix the colors
on its four external edges. Each color must
occur an even number of times: it appears once
at each of the four vertices, while every
internal edge accounts for two occurrences.
Thus, up to rotation and renaming colors,
the only possible boundary patterns are
$aaaa$, $aabb$, and $abab$.
The square has respectively $2$, $1$, and $0$
colorings. The two smoothings have respectively
$(1,1)$, $(1,0)$, and $(0,0)$ colorings,
proving the rule. In particular, all three
local identities hold for each fixed boundary
coloring.

We now show that the local rules reduce every closed
annular web to a nonnegative sum of collections of
circles. Set aside any isolated circles and consider one of the connected components that still has vertices.
It is trivalent and bipartite, so Euler's formula
on the sphere gives
\[
6=3F-E
=\sum_f\left(3-\frac{\deg f}{2}\right),
\]
where $\deg f$ counts the edges along the boundary
of $f$, with multiplicity.

Since the graph is bipartite, every face has even
degree, at least $2$. At most two faces contain
the ends of the annulus, and together they contribute
at most $4$ to the sum. If all other faces had
degree at least $6$, the sum would be at most $4$,
a contradiction. Thus some contractible face is
a bigon or a square. The component is bridgeless, since every edge lies in
a two-color cycle, and a facial walk in a bridgeless
cubic graph cannot repeat a vertex. Applying the corresponding
local rule decreases the number of vertices.
Repeating this process leaves only circles.
This is the annular $SL_3$ reduction of
\cite[Proposition~6.1]{DKS}.

It remains to evaluate the circles. A contractible
circle contributes $3$, while an essential circle
contributes either
\[
x_1+x_2+x_3=e_1
\quad\text{or}\quad
x_1^{-1}+x_2^{-1}+x_3^{-1}=e_2,
\]
depending on its orientation, since $e_3=1$.
All reductions have nonnegative coefficients,
so we obtain
\[
Z_H\big|_{e_3=1}
=\sum_{a,b}c_{ab}e_1^ae_2^b,
\qquad c_{ab}\in\NN.
\]
Lemma~\ref{lem:elementary-lift} now restores
$e_3$, proving elementary positivity.
\end{proof}

We can now prove Stembridge's conjecture for
partitions with at most three columns.

\begin{corollary}\label{cor:threecolumns}
Stembridge's interval-product conjecture
$\mathcal S(\mu)$ holds whenever $\mu_1\le3$.
\end{corollary}

\begin{proof}
Let $U=X_{J_b}\cdots X_{J_1}$ and
$Q=P_{J_b}\cdots P_{J_1}$.
If some $|J_t|>3$, then $\Phi(Q)$ vanishes
in three variables. Otherwise, the gadget
lemma expresses $\Phi(Q)$ as a positive
multiple of $e_1^u e_3^{f_0}Z_H$, which is
elementary positive by
Proposition~\ref{prop:threecolor}.
Thus $\Phi(Q)$ is elementary positive in
three variables in either case.

By the sign-twist identity \eqref{eq:sign-twist},
\[
\phi^\mu(U)=\kappa[e_\mu]\Phi(Q)\ge0
\]
for every $\mu_1\le3$, as desired.
\end{proof}

\begin{figure}[htbp]
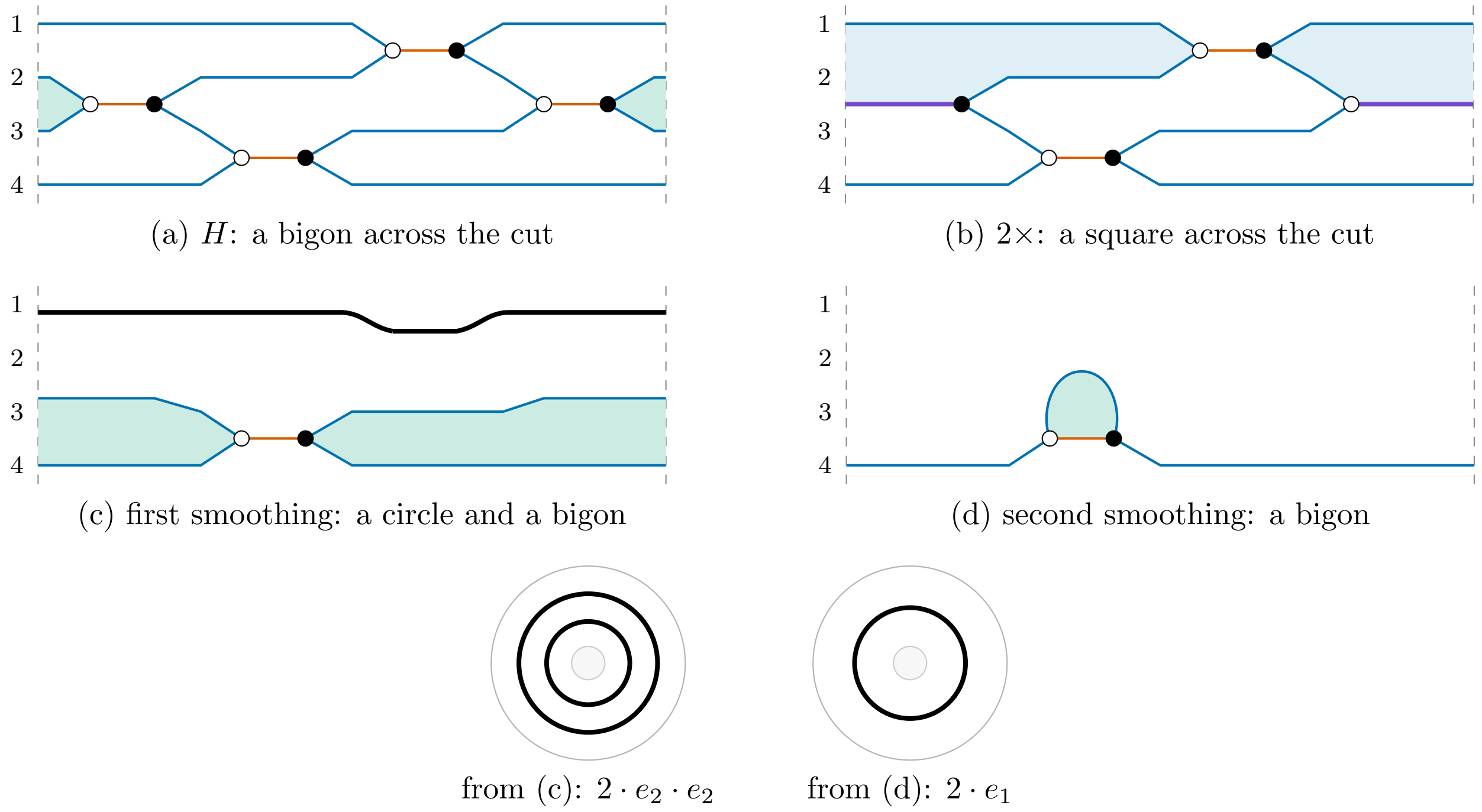

\centering
\ReduceA\hfill\ReduceB\\[6pt]
\ReduceC\hfill\ReduceD\\[6pt]
\ReduceE
\caption{Reducing the three-color gadget from
Example~\ref{ex:3412}. The dashed sides of
each strip are identified. Removing the bigon
in (a) gives a factor of $2$ and the graph
in (b). Its shaded square has two smoothings,
shown in (c) and (d). Each leaves a bigon,
whose reduction gives the essential circles
shown below.}
\label{fig:3412}\label{fig:word-reduction}
\end{figure}

\begin{example}[Four positions, three colors]\label{ex:3412}
Let us return to the interval product $U$ from
Example~\ref{ex:word-intro}. Its corresponding
antisymmetrizer product is
\[
Q=P_{[2,3]}P_{[1,2]}P_{[3,4]}P_{[2,3]}.
\]
Recall that $\Phi(U)=8h_4+4h_3h_1+4h_2^2$.
Applying the sign twist gives
\[
16\Phi(Q)=8e_4+4e_3e_1+4e_2^2.
\]
Now construct the three-color gadget $H$.
Each interval has size two, so each gadget has
one joining edge. The minimum flux is zero,
and the gadget lemma gives
\[
Z_H=16\Phi(Q)(x_1,x_2,x_3)
=4e_3e_1+4e_2^2,
\]
where $e_4=0$ in three variables.

Let us see how the local rules give the same
answer. Work with $e_3=1$ and follow
Figure~\ref{fig:word-reduction}.
The bigon across the cut in (a) contributes
a factor of $2$, leaving the graph in (b).
Applying the square rule gives the two
branches (c) and (d). In (c), one essential
circle separates off, and the remaining
bigon reduces to a second essential circle,
with a factor of $2$. In (d), the bigon
reduces to a single essential circle,
again with a factor of $2$.

An essential circle contributes either $e_1$
or $e_2$. Since $Z_H$ is homogeneous of
degree $4$, its terms on $e_3=1$ must have
degree congruent to $4$ modulo $3$.
Thus the two circles in (c) contribute
$e_2^2$, while the single circle in (d)
contributes $e_1$. Altogether,
\[
Z_H\big|_{e_3=1}
=2(2e_2^2+2e_1)
=4e_2^2+4e_1.
\]
Restoring $e_3$ by homogeneity gives
$Z_H=4e_2^2+4e_3e_1$, as expected.
\end{example}

\section{Four colors}\label{sec:fourcolors}

\subsection{The local rules and the two measures of progress}

The four-color proof follows the same general idea
as the three-color proof: use positive local relations
to reduce a web to circles. The difficulty is that
some of the four-color relations only rearrange the
web, without simplifying it. We will show that such
neutral moves always lead to a genuine reduction.
For this, we need a slightly larger class of webs
and two measures of progress.

We allow an edge labeled $a$ to carry an unordered
$a$-element subset of $[4]$, with the sets on edges
incident to each vertex partitioning $[4]$.
We call such a graph a \emph{weighted web} $W$,
and define its coloring polynomial $Z_W$ using
matching heights as before. An ordinary
$4$-regular bipartite multigraph has all edge
labels equal to $1$. Expanding an edge labeled
$a$ into $a$ distinguished parallel edges
multiplies the coloring polynomial by $a!$.
We use the subset-coloring normalization for
the intermediate webs.

We will prove the following stronger statement,
which allows us to apply local reductions without
leaving the class of objects under consideration.

\begin{proposition}\label{prop:fourcolor}
Every four-color annular weighted web has
$Z_W\in\NN[e_1,e_2,e_3,e_4]$.
\end{proposition}

To prove the proposition, we will simplify webs
using local identities that preserve their coloring
polynomials. We begin with a few elementary
identities, all valid for each fixed coloring
of the exterior edges.

A bivalent vertex with labels $a$ and $4-a$
can be removed by identifying its neighbors,
since the two color sets are complementary.
This preserves planarity and bipartiteness;
if the neighbors coincide, an isolated circle
results. A multiplicity-$4$ dipole contributes
$1$. Parallel edges of labels $a,b$ bounding
a disk merge with factor $\binom{a+b}{a}$,
counting the ways to split their union.
In particular, $(1,1)$ and $(1,2)$ bigons
contribute factors $2$ and $3$.
A contractible label-$a$ circle contributes
$\binom4a$.

Using these simplifications, we may assume that
every remaining vertex has type $112$ or $1111$,
where the numbers are the incident edge labels.
We call edges labeled $2$ \emph{heavy}, and edges
labeled $1$ \emph{unit edges}. Every trivalent
vertex is incident to exactly one heavy edge.

The three square relations we need are shown
in Figure~\ref{fig:four-rules}. They are special
cases of the relations of Fraser--Lam--Le
\cite[Section~6]{FLL}. All coefficients are
nonnegative, and the identities hold for each
fixed coloring of the exterior edges.

\par\medskip\noindent\begin{minipage}{\linewidth}
\captionsetup{type=figure}
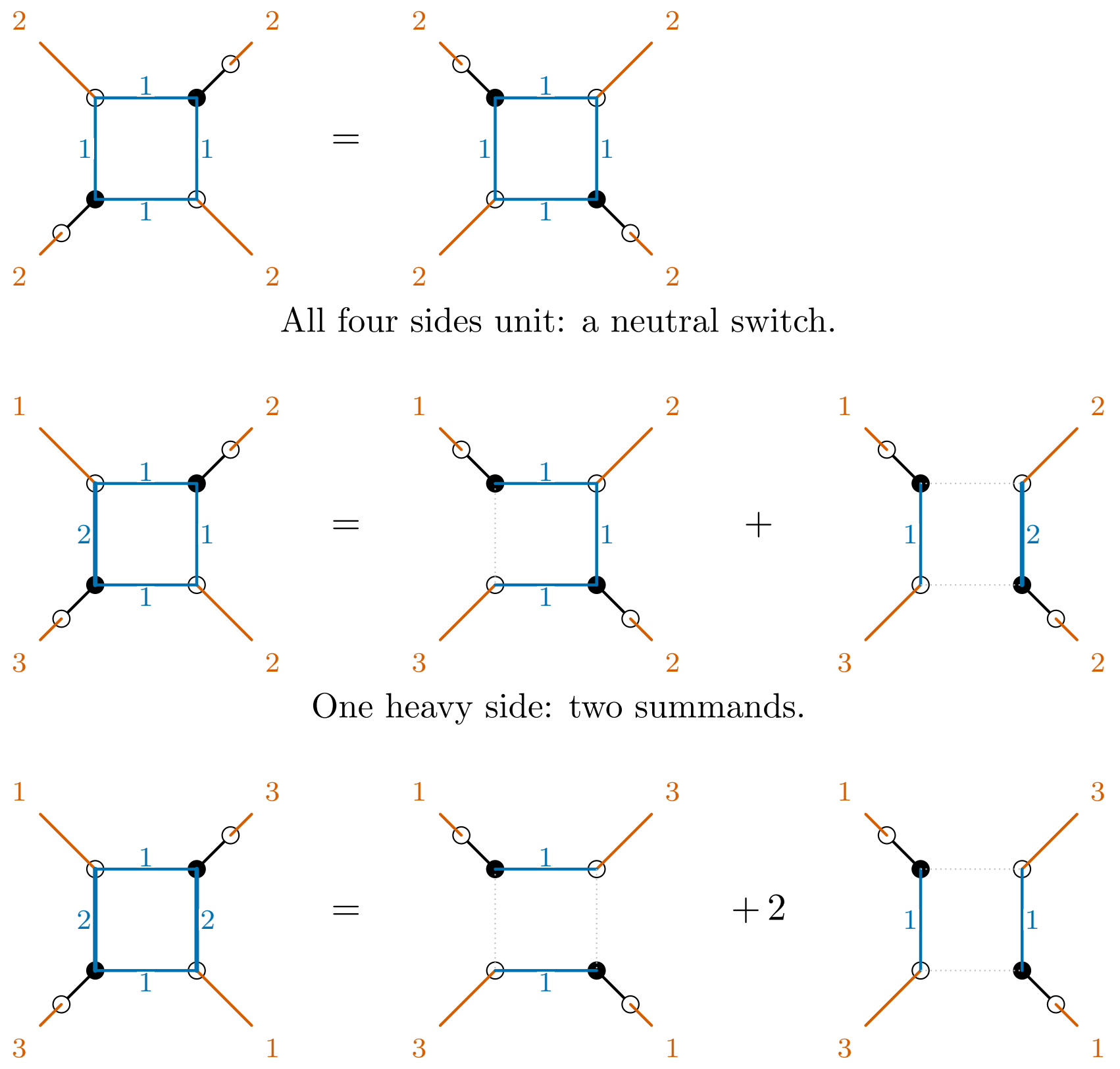
\centering
\begin{tabular}{c@{\quad}c@{\quad}c@{\quad}c@{\quad}c}
\FourSmallFill{2}{2}{1}{0}&$=$&
\FourSmallFill{2}{2}{1}{1}&&\\[.6em]
\multicolumn{5}{c}{\small All four sides unit: a neutral switch.}\\[1em]
\FourSmallFill{3}{2}{1}{0}&$=$&
\FourSmallFill{3}{2}{1}{1}&$+$&
\FourSmallFill{3}{2}{0}{1}\\[.6em]
\multicolumn{5}{c}{\small One heavy side: two summands.}\\[1em]
\FourSmallFill{3}{1}{1}{0}&$=$&
\FourSmallFill{3}{1}{1}{1}&$+\,2$&
\FourSmallFill{3}{1}{0}{1}\\[.6em]
\multicolumn{5}{c}{\small Two opposite heavy sides: coefficients $1$ and $2$.}
\end{tabular}
\caption{The four-color square rules. The four
exterior color sets are fixed. Orange numbers
are port sizes; blue numbers label square edges.
Dotted segments are absent. Each unlabeled short
segment has the unique multiplicity making the
total at a vertex $4$. Bivalent vertices are
removed afterward.}
\label{fig:four-rules}
\end{minipage}\par\medskip

\begin{lemma}[The four-color square rules]
\label{lem:four-square}
The three relations in Figure~\ref{fig:four-rules}
hold for every fixed coloring of the exterior edges.
\end{lemma}

\begin{proof}
Fix the colors on the exterior edges and consider
one color at a time. On the square, each color
must occupy a matching: it is either absent,
forced onto one edge, or free to choose between
the two pairs of opposite edges. These are the
only possibilities; all other boundary patterns
give zero on both sides.

The switch exchanges absent and free colors,
while forced colors move to opposite edges.
If $p$ and $q$ are the numbers of free colors
before and after switching, the numbers of
fillings are related by
\[
\binom pk=\sum_{t=0}^m\binom mt\binom q{k-t},
\qquad m=p-q.
\]
Here $k$ is the bottom label minus the number of colors
forced onto that edge. Since the top and bottom labels
agree, their forced counts agree. Switching and lowering
the bottom label by $t$ therefore replaces $k$ by $k-t$.
Counting two sides per free color and one per forced color,
the four side labels sum to $2p+(4-p-q)=4+m$. Thus
$m=0,1,2$ in the three rows of Figure~\ref{fig:four-rules},
respectively. Thus
Vandermonde's identity gives exactly the displayed
coefficients. In the last row, the term $t=2$
vanishes because the bottom label is $1$.
\end{proof}

We now explain how to use these rules for induction.
Let $C(W)=Z_W(1,1,1,1)$ be the number of
proper colorings of $W$. It is positive,
since every regular bipartite multigraph
decomposes into perfect matchings.

Call a positive local identity
$W=\sum_i c_iW_i$ \emph{branching} if
$\sum_i c_i\ge2$. Since
\[
C(W)=\sum_i c_iC(W_i),
\qquad C(W_i)>0,
\]
every resulting web satisfies $C(W_i)<C(W)$.
Thus branching gives us an induction on the
number of colorings, even if the number of
vertices does not decrease.

The first square rule is different: it replaces
one web by another with the same number of
colorings. We call such a move \emph{neutral}.
Our task is therefore to show that neutral
switches can always lead to a branching move.
For this we use a second measure, the number
of heavy edges. It need not decrease at every
switch, but will let us prove that branching
cannot be avoided.

\subsection{The contracted graph and its heavy-edge dual}

To understand how neutral switches can change the
heavy edges, we introduce two auxiliary graphs.
First, contract every heavy edge of $W$. The
resulting embedded graph $U$ is four-valent.
Its vertices are of two types: type $A$ comes
from an original four-valent vertex of $W$,
while type $B$ comes from contracting a heavy
edge between two trivalent vertices.
Contraction does not change the faces.
Viewing the annulus as a sphere with two
points removed, we can color the faces of
$U$ in checkerboard fashion, since every
vertex has even degree. We call the two
face colors \emph{shaded} and \emph{unshaded},
to distinguish them from the bipartition of $W$.

\begin{figure}[htbp]
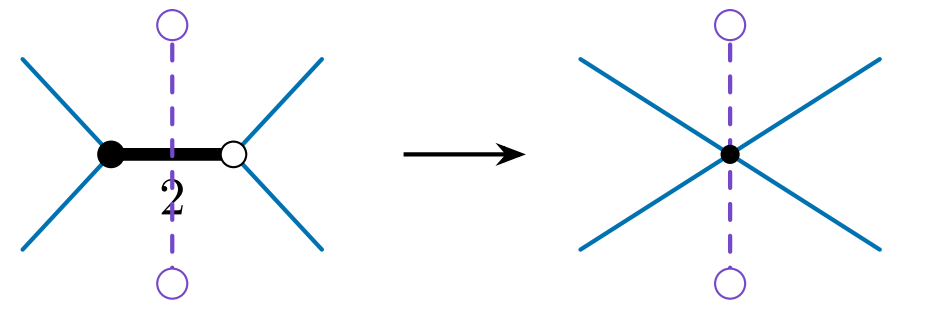

\centering\FourContractionFigure
\caption{Contracting a heavy edge gives a
type-$B$ vertex. The dashed dual edge joins
opposite faces of the same checkerboard color.
An original four-valent vertex (type $A$)
has no heavy dual edge.}
\label{fig:four-contraction}
\end{figure}

Next, define the \emph{heavy-edge dual} $R$.
Place a vertex in each face and draw an edge
across every heavy edge of $W$. After contraction,
this dual edge joins two opposite sectors at
a type-$B$ vertex of $U$, as in
Figure~\ref{fig:four-contraction}.
Opposite sectors have the same checkerboard
color. Thus every connected component of $R$
is entirely shaded or entirely unshaded,
and components of opposite colors do not cross.

For a face $f$, let $u(f)$ and $t(f)$ count
the incidences of unit and heavy edges along
its boundary. Thus $u(f)$ is the degree of
$f$ in $U$, while $t(f)$ is its degree in
the dual graph $R$. The original boundary
has length $u(f)+t(f)$.
Since no vertex is incident to two heavy edges,
heavy edges cannot occur consecutively along
a face boundary. Also, $W$ is bipartite.
Consequently,
\[
t(f)\le u(f),\qquad u(f)+t(f)\text{ is even}.
\]

We will use the fact that facial boundaries
are cycles, allowing two-vertex digons.
Indeed, a connected weighted web has no cut
vertices. To see this, delete a black vertex
$v$ and consider a component $K$ of the
remaining graph. If $b_K,w_K$ are its numbers
of black and white vertices, the total
multiplicity of edges from $v$ to $K$ is
\[
4(w_K-b_K)>0,
\]
by counting degrees in $K$. These are
positive multiples of $4$, but their sum
over all components is $4$. Thus there
is only one component. The argument for
a white vertex is identical.

To use Euler's formula, cap the annulus to obtain a
sphere. Mark the two faces containing its ends; if
both ends lie in the same face, mark one other face
as well. Every unmarked face is then a disk where
we can apply the local rules.

\begin{figure}[htbp]
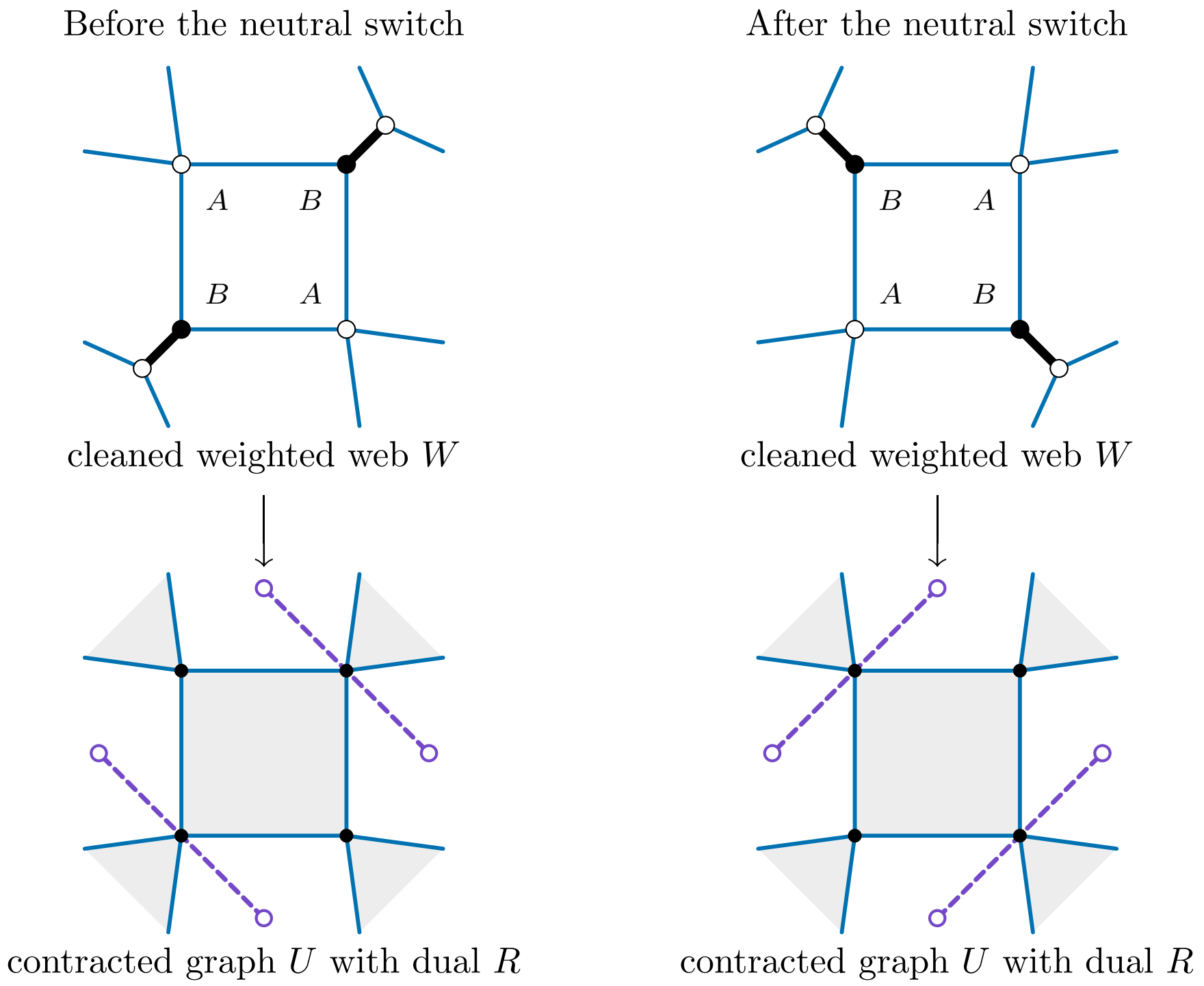

\centering
\begin{tabular}{@{}c@{\qquad\qquad}c@{}}
  \small Before the neutral switch &
  \small After the neutral switch \\[1pt]

  \NeutralSwitchPanel{0}{0} &
  \NeutralSwitchPanel{1}{0} \\[-2pt]

  \small cleaned weighted web $W$ &
  \small cleaned weighted web $W$ \\[4pt]

  $\Big\downarrow$ & $\Big\downarrow$ \\[-2pt]

  \NeutralSwitchPanel{0}{1} &
  \NeutralSwitchPanel{1}{1} \\[-2pt]

  \small contracted graph $U$ with dual $R$ &
  \small contracted graph $U$ with dual $R$
\end{tabular}

\caption{A neutral switch after bivalent cleanup.
The square corners have types $A$ (quartic) or
$B$ (two cubic vertices joined by a heavy edge);
the alternating arrangement shown is only an example.
Below, the contracted graph $U$ is unchanged,
while the purple dashed edges of its heavy-edge
dual $R$ are toggled. Gray indicates one
checkerboard color.}
\label{fig:neutral-explained}
\end{figure}

\begin{lemma}\label{lem:four-local}
Suppose no branching move is available.
Every unmarked face meeting a heavy edge has
boundary length at least $6$. An unmarked face
with no heavy edges is either an all-unit square
or has even length at least $6$.

Moreover, a neutral switch at an all-unit square
leaves the contracted graph $U$ unchanged and
only changes heavy dual connections of the
opposite checkerboard color. In particular,
other all-unit squares of the same color
remain available for switching.
\end{lemma}

\begin{proof}
Recall that a face has $u$ unit-edge incidences
and $t$ heavy-edge incidences, with $t\le u$
and $u+t$ even. The possibilities with
$u+t<6$ are
\[
\begin{array}{c|ccccc}
(u,t)&(1,1)&(2,0)&(3,1)&(2,2)&(4,0)\\ \hline
\text{move}&3\times&2\times&
\text{branch}&\text{branch}&\text{neutral}.
\end{array}
\]
The first two cases are digons, handled by our elementary
reductions. To apply the square rules at a quartic corner,
separate the two square edges from the two exterior edges
by a path of two weight-two edges. This reverses bivalent
cleanup. If an exterior heavy edge joins two square corners,
subdivide it twice so the two ports can be treated separately.

The squares with one or two heavy sides admit branching.
This leaves only the all-unit square and proves the first
assertion.

Now consider a neutral switch at an all-unit
square. As illustrated in
Figure~\ref{fig:neutral-explained}, the switch
may exchange types $A$ and $B$ at the corners,
but after contracting heavy edges the graph
$U$ is unchanged. Each such exchange toggles
a heavy dual connection between opposite
sectors of the checkerboard color opposite
to the square. 

Thus connections of the square's own color remain unchanged,
and other all-unit squares of that color remain available.
If an exterior heavy edge joins two switched corners,
its dual connection is toggled once at each end, so
the two toggles cancel.
\end{proof}

\subsection{Switching inside an empty hole}

We now use the heavy-edge dual $R$ to find neutral
switches that decrease the number of heavy edges.
Recall that isolated vertices of $R$ correspond
to faces with no heavy edges. Discard those
corresponding to unmarked all-unit squares, and
call each remaining connected component of $R$
an \emph{object}. Thus isolated vertices
corresponding to marked faces or larger faces
are also objects.

Draw a connected component $K$ of $R$ in the plane,
ignoring the other components. Its bounded faces
are called \emph{holes}. Their boundaries need
not be simple cycles: entire trees may hang into
a hole, with each of their edges appearing twice
in its boundary walk.

A hole is \emph{empty} if it contains no other
objects, although it may contain discarded
vertices corresponding to all-unit squares.
We first show how an empty hole can be used to
make progress. In the next subsection, we will
prove that such a hole must exist.

\begin{lemma}\label{lem:four-empty}
Suppose no branching move is available, and let
$D$ be an empty hole of an object of $R$.
Then neutral switches inside $D$ either expose
branching or strictly decrease the number of
heavy edges, without changing $U$.
\end{lemma}

\begin{proof}
Let $K$ be the object surrounding the hole,
and suppose it is unshaded. Every shaded
vertex of $R$ inside the hole must correspond
to an unmarked all-unit square; otherwise
it would belong to another object.
Moreover, at least one such square exists.
Indeed, a nonbridge edge on the boundary of
the hole has a shaded sector on its inner side,
giving a shaded vertex of $R$ inside the hole.

We switch all these squares, stopping if
branching becomes available. By
Lemma~\ref{lem:four-local}, switching one
shaded square does not prevent us from
switching the others.

Let us see what these switches do to the
unshaded heavy dual connections. Each such
connection can be toggled from either of
its two opposite shaded sectors. If both
sectors lie inside the hole, it is toggled
twice and remains unchanged. If neither
lies inside, it is not affected.

It remains to consider a connection with exactly
one shaded sector inside the hole. We claim that
this connection is already present in $K$.
Otherwise, a short arc through that vertex
of $U$ would connect the two shaded sectors
without crossing $K$. They would then lie
in the same complementary region of $K$,
a contradiction. Thus every connection
toggled exactly once is removed, and no new
connection is created.

In fact, every nonbridge edge on the boundary
of the hole is removed: it has the hole on
one side and another face on the other, so
it is toggled exactly once. Bridges, including
edges of trees hanging into the hole, have
the hole on both sides and are toggled twice.
Since the boundary contains a cycle, at least
one heavy edge is removed. No heavy edges
are added, and $U$ remains unchanged.
\end{proof}

\begin{figure}[htbp]
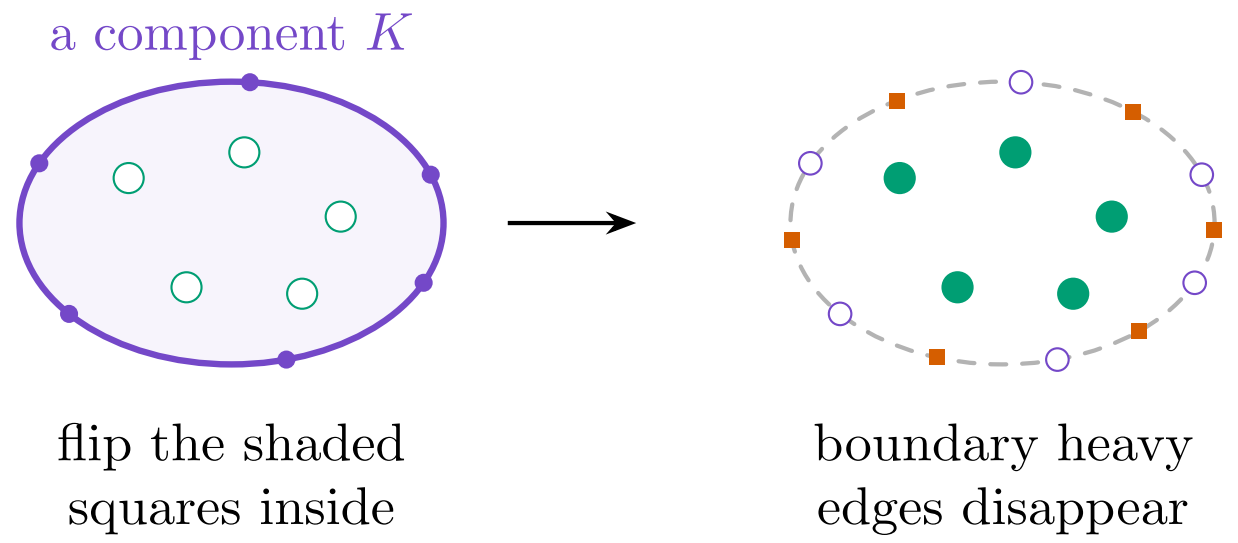

\centering
\FourEmptyFigure
\caption{Switching inside an empty hole.
The green circles represent available
opposite-color squares. Switching them
removes the heavy edges along the boundary
of the hole, shown here for a simple cycle.
Orange squares mark the resulting quartic
vertices.}
\label{fig:four-empty}
\end{figure}

\begin{example}
\label{ex:empty-hole-preview}
The mechanism of Lemma~\ref{lem:four-empty} is
illustrated by the permutohedron web shown 
in Figure~\ref{fig:four-permutohedron}. In the
heavy-edge dual, the inner circuit bounds an
empty hole containing one available all-unit
square. Switching that square removes all four
heavy edges on the boundary of the hole.

After the switch, the outer circuit bounds a
larger hole, and the four edges joining the
inner and outer circuits become dangling bridges.
Thus the new hole has a nonsimple boundary:
its boundary walk traverses each of those bridges
twice. This is exactly the situation discussed
in the proof above.
\end{example}

\subsection{Finding an empty hole}

We have seen how to make progress once an empty hole
is available. We now show that such a hole must exist
whenever branching is unavailable. The proof uses
Euler's formula to show that the heavy-edge dual
has too many holes for all of them to contain
other objects.

\begin{lemma}\label{lem:four-hole}
If branching is unavailable in a connected
nonterminal web, some object of $R$ has an empty hole.
\end{lemma}

\begin{proof}
Cap the annulus to a sphere, and recall that exactly
two faces are marked (choosing an arbitrary second
face if both ends lie in the same one).

Let $F$ be the number of faces of $U$, let $s$ be
the number of discarded all-unit squares, and let
$h$ be the number of heavy edges. Since $U$ is
connected and four-valent, $E(U)=2V(U)$, so Euler's
formula gives $V(U)=F-2$.

Recall that the original boundary length of a face
is $\ell(f)=u(f)+t(f)$. Every unit edge contributes
two boundary incidences, as does every heavy edge.
Thus
\[
\sum_f\ell(f)=2E(U)+2h=4(F-2)+2h.
\]
On the other hand, the $s$ discarded squares have
length $4$, every other unmarked face has length
at least $6$ by Lemma~\ref{lem:four-local}, and
the two marked faces have length at least $2$.
Hence
\[
4(F-2)+2h\ge 4s+6(F-s-2)+4,
\]
which simplifies to
\[
h\ge F-s.
\]

Now discard the $s$ isolated square vertices
from $R$. The remaining graph is planar, with
$F-s$ vertices, $h$ edges, and say $N$ connected
components, which are precisely our objects.
Let $B$ be the total number of their holes.
Euler's formula gives
\[
B=h-(F-s)+N\ge N.
\]
Thus there are at least as many holes as objects.

We claim that some hole must be empty. Otherwise,
choose infinity in a marked face and consider the
nesting of the objects in the plane. Every hole
contains another object, and hence an immediate
child in the nesting forest. Different holes give
different parent--child edges, since each object
has at most one parent. The forest therefore has
at least $B$ edges. But a forest on $N$ vertices
has at most $N-1$ edges, contradicting $B\ge N$.
\end{proof}

We can now finish the argument that neutral switches
always lead to branching.

\begin{proposition}\label{prop:four-availability}
Every connected nonterminal four-color web can reach
branching by neutral square switches and bivalent
cleanup.
\end{proposition}

\begin{proof}
Suppose no reachable state admits branching,
and choose a reachable cleaned state with the
fewest heavy edges. By Lemma~\ref{lem:four-hole},
it has an empty hole. Lemma~\ref{lem:four-empty}
then gives a sequence of neutral switches that
decreases the number of heavy edges, contradicting
our choice.
\end{proof}

The following example shows why neutral switches
are needed: branching may not be available
in the original web.

\begin{figure}[tbp]
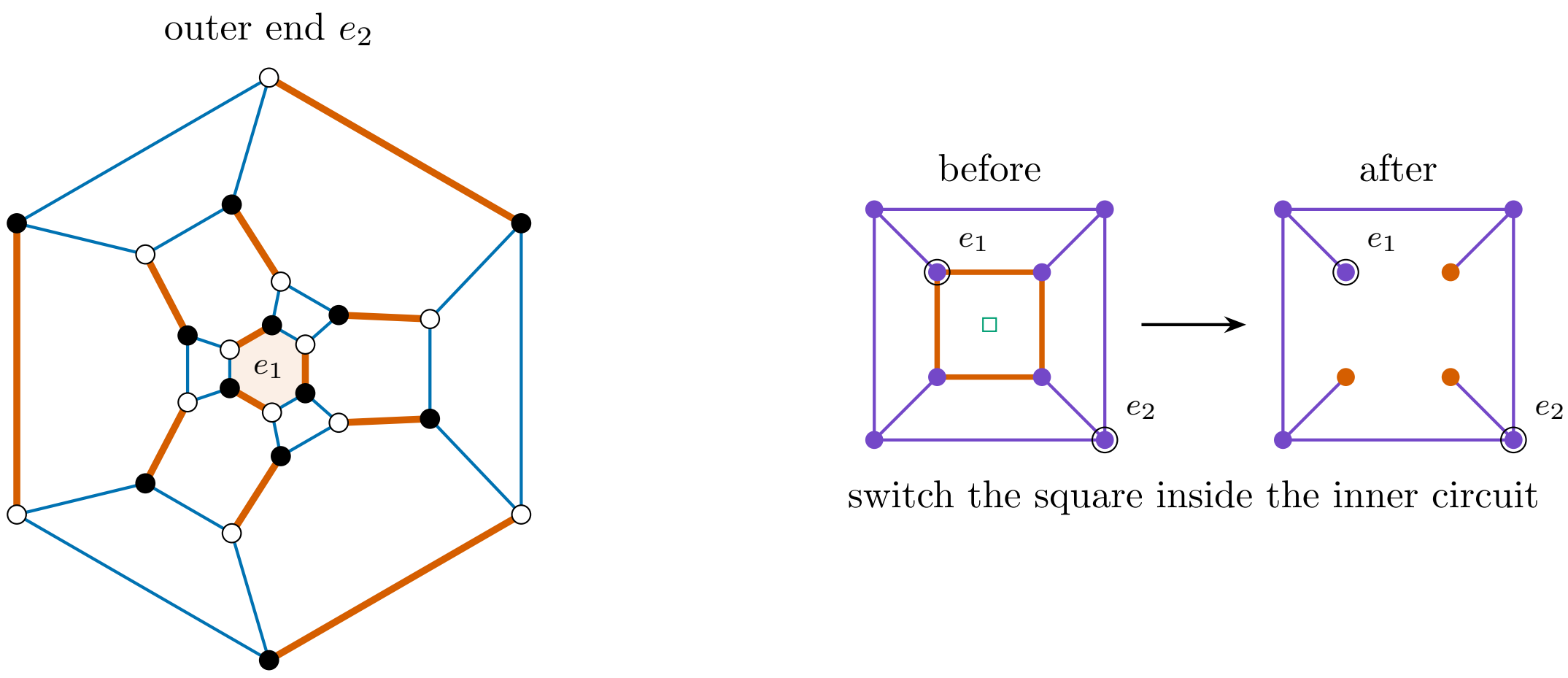

\centering
\begin{minipage}{.39\linewidth}
\centering\FourPermutohedronFigure
\end{minipage}\hfill
\begin{minipage}{.58\linewidth}
\centering\FourCubeFigure
\end{minipage}
\caption{Left: an annular web with no initial
branching move. Orange edges have label $2$,
and blue edges label $1$. The central hexagon
and the exterior are the two ends of the annulus.
Right: its cube-shaped heavy-edge dual before
and after a neutral switch. The inner circuit
disappears, leaving dangling edges in the
larger hole and exposing branching squares.
Ringed vertices mark the two ends.}
\label{fig:four-permutohedron}
\end{figure}

\begin{example}[A neutral move is necessary]
\label{ex:four-permutohedron}
Consider the permutohedron graph, with one vertex
for each permutation in $S_4$ and an edge of
type $i$ whenever two permutations differ by
swapping positions $i,i+1$. Give edges of
types $1,3$ label $1$ and edges of type $2$
label $2$. Every vertex then has type $112$.
The six square faces are all-unit, while the
eight hexagons have alternating labels $1,2$.
Puncture two opposite hexagons to obtain an
annular web. Initially, no branching move
is available.

Its heavy-edge dual is a cube. Switching the
square inside one face of the cube removes
all four edges around that face. Three other
unmarked dual vertices now have $u=3,t=1$:
their original hexagons have become branching
$2111$ squares. Figure~\ref{fig:four-permutohedron}
shows this first step.
\end{example}

\subsection{Positive reduction and the four-column theorem}

We can now finish the four-color argument.

\begin{proof}[Proof of Proposition~\ref{prop:fourcolor}]
 Work with $e_4=1$,
as in the three-color proof. The local relations
hold for each fixed boundary coloring, and moving
the cut across a vertex changes the weight by
a power of $e_4=1$. Thus we can apply these
relations in any contractible disk, using signed
crossings instead of normalized heights.

We induct on $C(W)$, allowing disconnected webs.
After cleanup, choose a component that still has vertices.
By Proposition~\ref{prop:four-availability}, neutral switches
expose a branching move in that component. Each summand has
a smaller coloring count on the part that was changed.
Since coloring counts multiply over components, each
summand also has a smaller total $C(W)$. The induction
therefore reduces the web to a positive sum of collections
of circles.

A contractible circle labeled $a$ contributes
$\binom4a$. An essential circle contributes
$e_a$ in one orientation and $e_{4-a}$ in the
other, since $e_4=1$. Thus every resulting
collection has an elementary-positive evaluation.

Finally, $Z_W$ is homogeneous, so
Lemma~\ref{lem:elementary-lift} restores $e_4$
and proves the proposition.
\end{proof}

We now obtain our main immanant result.

\begin{proof}[Proof of Theorem~\ref{thm:fourcolumns-intro}]
Elementary positivity for one and two colors is
immediate, and the three-color case was proved
in Section~\ref{sec:colorings}. The four-color
case is Proposition~\ref{prop:fourcolor}.

It remains to translate this into Stembridge's
conjecture. Let
\[
U=X_{J_b}\cdots X_{J_1},\qquad
Q=P_{J_b}\cdots P_{J_1}.
\]
If some interval has size greater than $4$,
then $\Phi(Q)$ vanishes in four variables.
Otherwise, the gadget lemma expresses $\Phi(Q)$
as a positive multiple of a four-color coloring
polynomial, together with powers of $e_1$ and
$e_4$. Hence $\Phi(Q)$ is elementary positive.

By the sign-twist identity \eqref{eq:sign-twist},
\[
[e_\mu]\Phi(Q)
=\frac{\phi^\mu(U)}{\prod_t|J_t|!}.
\]
In four variables, the $e_\mu$ with $\mu_1\le4$
form a basis. Their coefficients are therefore
nonnegative, proving $\mathcal S(\mu)$ for
every partition with at most four columns.
\end{proof}

Let us return once more to the interval product
from Example~\ref{ex:word-intro}. Its compressed
four-color gadget has one heavy edge at each
of the four events, and
\[
Z_W=16\Phi(Q)=8e_4+4e_3e_1+4e_2^2.
\]
Expanding the four heavy edges into distinguished
pairs multiplies the coloring polynomial by $2^4$.
Unlike in the three-color calculation of
Example~\ref{ex:3412}, the term $8e_4$ now survives.

The two-column case was previously proved using
Temperley--Lieb immanants \cite{CSS}. Our
three-column proof uses the classical annular
$SL_3$ reduction, while the fourth column
requires the neutral-switch argument developed
above.

\Needspace{11\baselineskip}
%\begin{samepage}
\section{Stretching}\label{sec:stretching}

The stretching theorem has a simple interpretation
in terms of planar networks. Recall that
\[
m_{n\mu}(\rho_1,\rho_2,\ldots)
=m_\mu(\rho_1^n,\rho_2^n,\ldots).
\]
On the other hand, concatenating $n$ copies of
a planar network replaces its path matrix $R$
by $R^n$, whose eigenvalues are the $n$th powers
of those of $R$. Combining these observations
with Corollary~\ref{cor:pathschur} gives the proof.

\begin{proof}[Proof of Theorem~\ref{thm:stretching-intro}]
Put $M=n|\mu|$, and let $R$ be the path matrix
of any acyclic planar disk network with $M$
sources and sinks. Concatenate $n$ copies of
this network, initially giving each copy
independent edge variables. The resulting
network is again planar and acyclic, so
$\mathcal S(\mu)$ and
Corollary~\ref{cor:pathschur} imply that
$m_\mu$ of its path matrix has nonnegative
coefficients.

Identifying corresponding edge variables in
all copies preserves coefficientwise positivity
and turns the path matrix into $R^n$. Thus
\[
m_\mu(R^n)\in\mathbb Q_{\ge0}[\yy].
\]
But the eigenvalues of $R^n$ are the $n$th
powers of those of $R$, so
\[
m_\mu(R^n)=m_{n\mu}(R).
\]
We have therefore proved coefficientwise
nonnegativity of $m_{n\mu}(R)$ for every planar
disk path matrix $R$. By the equivalence in
Corollary~\ref{cor:pathschur}, this is exactly
$\mathcal S(n\mu)$.
\end{proof}

\begin{figure}[htbp]
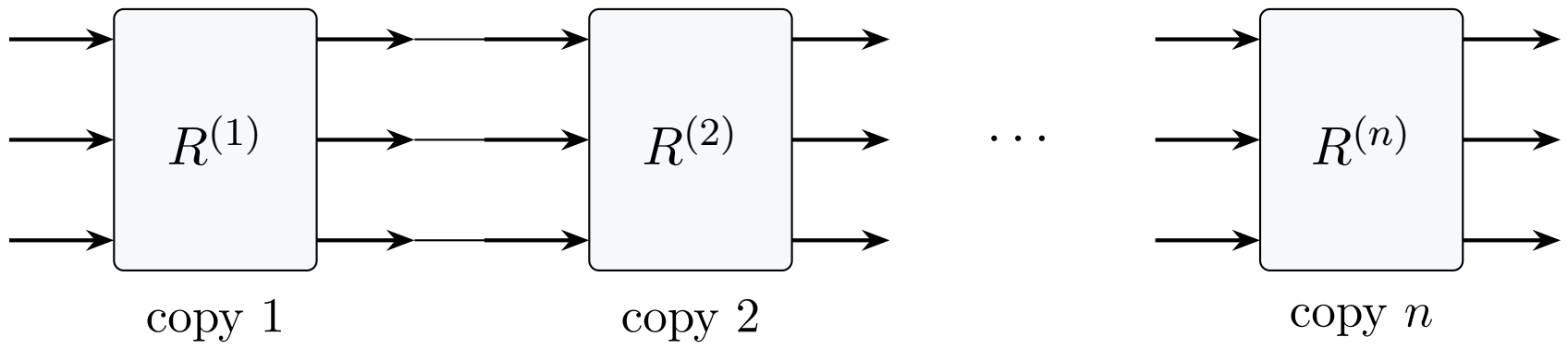

\centering\StretchFigure
\caption{Concatenating $n$ copies of a planar
network multiplies their path matrices.}
\label{fig:stretch}
\end{figure}

\begin{corollary}\label{cor:four-stretches}
Stembridge's assertion holds for every partition whose parts lie in $\{n,2n,3n,4n\}$ for some $n\ge1$.
\end{corollary}
\begin{proof}
Use Theorem~\ref{thm:fourcolumns-intro}, then stretch.
\end{proof}

\begin{corollary}\label{cor:hook-stretches}
Stembridge's assertion holds for all stretched hooks $(an,n^b)$.
\end{corollary}
\begin{proof}
McNulty's hook theorem \cite[Theorem~1.3]{McNulty} gives nonnegative evaluations of $\phi^{(a,1^b)}$ on Kazhdan--Lusztig basis elements at $q=1$. Each $X_J$ is $C'_{w_0(J)}(1)$, where $w_0(J)$ is the longest permutation in $S_J$. Nonnegativity of the Kazhdan--Lusztig structure constants \cite{EW} therefore makes every product of interval sums a nonnegative combination of these basis elements. Thus $\mathcal S((a,1^b))$ holds. Apply stretching.
\end{proof}
The interval-product hypothesis here is supplied by the stronger hook theorem \cite{McNulty}; the numerical hook theorem of Dong--Wu--Gao \cite{DWG} alone is not enough.

The interval-product form of Stembridge's conjecture
also implies several other positivity statements.

\begin{corollary}[Matrix consequences]
\label{cor:matrix-consequences}
If $\mathcal S(\mu)$ holds for $\mu\vdash N$,
then:
\begin{enumerate}
\item The monomial immanant $\Imm_{\phi^\mu}(R)$
has nonnegative coefficients for every compatible
acyclic path matrix $R$
\cite[Conjecture~6.3]{Stembridge92}.
\item $\Imm_{\phi^\mu}(R)\ge0$ for every real
totally nonnegative $N\times N$ matrix $R$
\cite[Conjecture~2.1]{Stembridge92}.
\item Monomial immanants of Jacobi--Trudi matrices
have nonnegative coefficients in the symmetric-function
variables \cite[Conjecture~4.2(b)]{Stembridge92}.
\end{enumerate}
In particular, all three statements hold for
$\mu_1\le4$, for partitions with parts in
$\{n,2n,3n,4n\}$, and for stretched hooks.
\end{corollary}

\begin{proof}
The first assertion is Wolfgang's reduction,
Theorem~\ref{thm:reduction}. For the second,
recall that every totally nonnegative matrix
can be obtained as a limit of path matrices
of planar networks with nonnegative edge weights
\cite{Brenti95}. The first assertion gives
nonnegativity for these matrices, and the
inequality is preserved under limits.

Finally, Jacobi--Trudi matrices have the usual
lattice-path realization, with edge weights
given by monomials in the symmetric-function
variables. Applying the first assertion gives
coefficientwise nonnegativity. We may work
with a finite alphabet, since each coefficient
stabilizes.
\end{proof}

Part~(2) follows from the stronger coefficientwise
statement. Tang's recent preprint
\cite[Theorem~1.1]{Tang} claims part~(2) for
all shapes, but our coefficientwise results
are independent of his work.

\section{Directed cycles}\label{sec:networks}

We begin with the proof of Theorem~\ref{thm:network}.
Recall that $H_k$ is the sum of weights of all
families of $k$ vertex-disjoint directed cycles
in a cylindrical network. Orient the cylinder
so that a positively winding simple loop has
the inner boundary on its left. We may discard
edges that belong to no directed cycle, since
they do not contribute to any $H_k$.

\begin{lemma}[A directed cut]\label{lem:cut}
There is a simple arc joining the two boundary
components of the cylinder that crosses every
directed edge it meets from left to right.
Moreover, every simple directed cycle crosses
this arc exactly once.
\end{lemma}

\begin{proof}
Consider the planar dual graph, with each dual
edge directed from the face on the left of the
original edge to the face on its right.
We claim that there is a directed dual path
from the inner face to the outer face.

Otherwise, let $S$ be the union of faces
reachable from the inner face. No dual edge
leaves $S$, so every edge along its boundary
is directed with $S$ on its right. In particular,
the boundary separating the two ends of the
cylinder contains a directed closed walk of
negative winding. Decomposing this walk into
directed simple cycles contradicts our assumption
that every directed cycle winds positively.

Thus a directed dual path exists. Taking it
simple gives the desired arc. Every crossing
with a directed cycle has the same sign.
Since a simple positively winding cycle on
a cylinder has winding number $1$, its
intersection number with the arc is $1$.
Hence it crosses the arc exactly once.
\end{proof}

Cut the cylinder along the arc $\gamma$ from
Lemma~\ref{lem:cut}. The resulting network lies
in a disk and is acyclic: a directed cycle
avoiding $\gamma$ would have winding zero.
Suppose $\gamma$ crosses $N$ edges, labeled
$1,\ldots,N$ from inner to outer boundary.
Cutting each edge creates a source and a sink
with the same label, in the same order along
the two copies of $\gamma$. Keep its original
weight on one half and give the other half
weight $1$. Let $R$ be the path matrix of
this disk network.

Every directed cycle crosses $\gamma$ exactly
once. Thus a family of $k$ vertex-disjoint
cycles becomes a family of $k$ vertex-disjoint
paths, with the same subset $S$ of source and
sink labels. Conversely, such paths glue back
to disjoint cycles. Since the network is planar,
disjoint paths must preserve the order of their
endpoints. The Lindstr\"om--Gessel--Viennot
lemma \cite{Lindstrom,GV} therefore gives
\begin{equation}\label{eq:detIP}
H_k=\sum_{|S|=k}\det R[S,S],
\qquad
\sum_k H_kz^k=\det(\mathrm{Id}+zR).
\end{equation}

Let $\rho_1,\ldots,\rho_N$ be the eigenvalues
of $R$, counted with multiplicity. Then
\[
\det(\mathrm{Id}+x_iR)
=\prod_{j=1}^N(1+x_i\rho_j).
\]
Thus the dual Cauchy identity gives
\[
\begin{aligned}
\prod_{i=1}^r\left(\sum_k H_kx_i^k\right)
&=\prod_{i=1}^r\det(\mathrm{Id}+x_iR)\\
&=\prod_{i=1}^r\prod_{j=1}^N(1+x_i\rho_j)\\
&=\sum_\lambda s_\lambda(\xx)s_{\lambda'}(R),
\end{aligned}
\]
where $s_{\lambda'}(R)$ means evaluation on
the eigenvalues of $R$, as in Section~\ref{sec:schur}.

By Corollary~\ref{cor:pathschur}, every
$s_{\lambda'}(R)$ has nonnegative coefficients
in the edge weights. Thus the product on the
left is coefficientwise Schur positive, and
Lemma~\ref{lem:jt} gives coefficientwise
nonnegativity of all Toeplitz minors.
This proves Theorem~\ref{thm:network}, which is
\cite[Conjecture~6.1]{GP}.

The same argument proves Conjecture~6.2 of
Galashin--Pylyavskyy \cite{GP}. Indeed, their
specialization of symmetric functions is defined
by $e_k\mapsto H_k$. Since $H_k=e_k(R)$ by
\eqref{eq:detIP}, this specialization sends
$s_\lambda$ to $s_\lambda(R)$, which has
nonnegative coefficients by
Corollary~\ref{cor:pathschur}.

We also obtain their real-rootedness Conjecture~6.3.
For positive numerical edge weights, the P\'olya frequency
criterion implies that $H(z)=\sum_k H_kz^k$ has only
negative real roots. Let $d=\deg H(z)$. The recurrence
polynomial
\[
Q_N(t)=\sum_{k=0}^d(-t)^{d-k}H_k
\]
therefore has only positive real roots. As shown in
\cite{GP}, the recurrence polynomials for
tuples of paths have roots that are products
of these roots, with or without repetition.
These are therefore positive real as well.

\begin{figure}[ht]
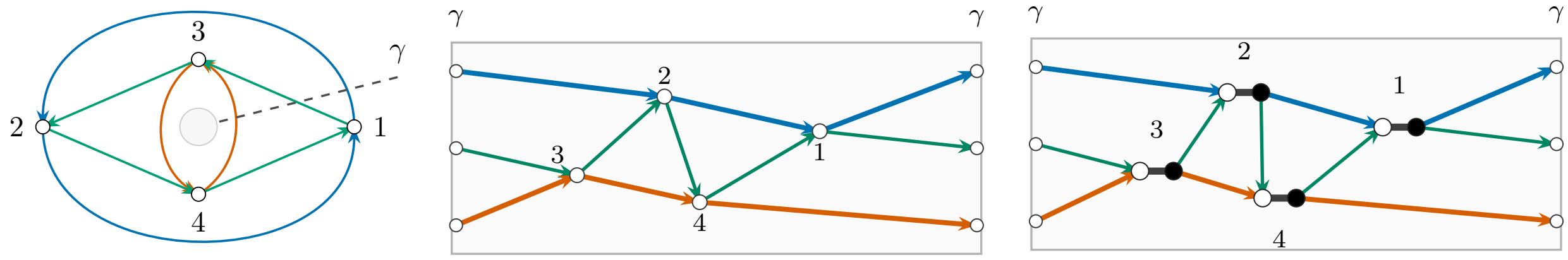

\centering
\begin{minipage}[c]{.265\linewidth}
  \centering\resizebox{\linewidth}{!}{\NetworkFigure}
\end{minipage}\hfill
\begin{minipage}[c]{.355\linewidth}
  \centering\resizebox{\linewidth}{!}{\StripFigure}
\end{minipage}\hfill
\begin{minipage}[c]{.355\linewidth}
  \centering\resizebox{\linewidth}{!}{\SplitStripFigure}
\end{minipage}
\caption{Three views of the same example.
Left: the cylindrical network obtained by
contracting the minimum matching of the ladder,
with cut $\gamma$.
Middle: after cutting along $\gamma$, the orange
and blue paths form the disjoint family with
$S=\{1,3\}$.
Right: each internal vertex is split into a white
incoming vertex and a black outgoing vertex,
joined by a thick matching edge.
Gluing the two copies of the cut recovers the original bipartite graph.}
\label{fig:network}
\end{figure}

For this network, setting all edge weights to $1$
gives the path matrix
\begin{equation}\label{eq:cutexample}
R=\begin{pmatrix}
2&3&3\\
2&3&3\\
1&2&2
\end{pmatrix},
\qquad
\det(\mathrm{Id}+zR)=1+7z+z^2.
\end{equation}
Thus the original cylindrical network has seven
simple directed cycles and exactly one pair of
vertex-disjoint cycles.

\subsection{Minimum-height matchings}

We now explain why directed networks and dimers
with a unique minimum-height matching are essentially
the same objects.

\begin{proposition}[Dimers and directed networks]
\label{prop:unique-minimum}
Let $G$ be a bipartite graph on a cylinder with
a unique minimum-height matching $M_0$.
Then $G$ gives rise to a planar cylindrical
directed network whose cycle measurements satisfy
\[
a_k=\wt(M_0)H_k
\]
after a Laurent monomial change of edge weights.
Conversely, after we delete the edges belonging to no
directed cycle, any planar cylindrical network arises
this way from a bipartite dimer model with a unique
minimum-height matching.
\end{proposition}

\begin{proof}
Orient the edges of $M_0$ from white to black,
and all other edges from black to white.
Contract the edges of $M_0$. Directed cycles
in the resulting network correspond to
alternating cycles along which we can switch
$M_0$. A negatively winding cycle would
contradict minimality, while a zero-winding
cycle would give another minimum-height matching.
Thus every directed cycle winds positively.

For a nonmatching edge $e$ leaving a black
vertex $b$, give the corresponding directed
edge weight $y_e/y_{m(b)}$, where $m(b)$ is
the matching edge at $b$. Switching $M_0$
along a family of vertex-disjoint alternating
cycles multiplies its weight by the product
of the corresponding directed edge weights.
Since every such cycle winds once, we obtain
$a_k=\wt(M_0)H_k$.

For the converse, start with a planar cylindrical
network and discard edges that belong to no
directed cycle. Cut it open as in
Lemma~\ref{lem:cut}. 
We claim that at every vertex, the incoming
edges form one consecutive group in the cyclic
order around the vertex, and the outgoing edges
form another. In other words, the two kinds of
edges do not alternate. This allows us to
separate the incoming and outgoing edges into
two vertices without introducing crossings.

Indeed, join the cut sources to a common source
along the boundary of the disk. If incoming
and outgoing edges alternated at a vertex $v$,
two paths from the common source to $v$ would
enclose a region, using their portions after
their last common vertex. One of the outgoing
edges would enter this region. Its continuation
to a sink would have to meet one of the
bounding paths, creating a directed cycle
in the disk network, a contradiction.

We can therefore split each vertex into a
white vertex receiving its incoming edges
and a black vertex carrying its outgoing
edges, joined by a new edge of weight $1$.
The consecutive-block property makes this
construction planar. The new edges form a
perfect matching $M_0$, and switching it
corresponds exactly to choosing disjoint
directed cycles. Since all these cycles
wind positively, $M_0$ is the unique
minimum-height matching, and $a_k=H_k$.
\end{proof}

\begin{corollary}\label{cor:unique-minimum}
If a bipartite graph on a cylinder has a unique
minimum-height matching, then its matching
measurements $(a_k)$ are coefficientwise
P\'olya frequency.
\end{corollary}

\begin{proof}
By Proposition~\ref{prop:unique-minimum},
the measurements satisfy
$a_k=\wt(M_0)H_k$ after a Laurent monomial
substitution of edge weights. Thus every
Toeplitz minor of order $r$ is $\wt(M_0)^r$
times the corresponding minor for $(H_k)$,
which is coefficientwise nonnegative by
Theorem~\ref{thm:network}. After substitution,
we obtain a Laurent polynomial with nonnegative
coefficients. Since the original dimer minor
is a polynomial in the edge weights, its
coefficients are nonnegative as well.
\end{proof}

\begin{example}[The ladder in both directions]
\label{ex:ladder-correspondence}
Consider the ladder $C_4\times P_2$, with its unique
minimum-height matching $M_0$ shown in
Figure~\ref{fig:matchings}. Contracting the four
edges of $M_0$ gives the network on the left of
Figure~\ref{fig:network}. It has seven simple
directed cycles and one disjoint pair, so at
unit weights
\[
A(z)=1+7z+z^2.
\]

The weights work just as simply. Let $p_i$ be
the weight of the matching edge contracted to
vertex $i$, and $x_{ij}$ the weight of the edge
giving the arrow $i\to j$. We give this arrow
weight $x_{ij}/p_i$. For example, switching
$M_0$ along the cycle $1\to3\to2\to1$ gives
a matching of weight
\[
p_4x_{13}x_{32}x_{21}
=\wt(M_0)\frac{x_{13}x_{32}x_{21}}{p_1p_3p_2}.
\]

Conversely, the right panel of
Figure~\ref{fig:network} shows how to split each
vertex into a white incoming vertex and a black
outgoing vertex, joined by a matching edge.
Gluing the two sides of the cut recovers the
original ladder and its matching $M_0$.
\end{example}

The preceding corollary assumes that the minimum-height
matching is unique. In the next section we will handle graphs where this is
not the case by starting with a larger graph that does
have a unique minimum, and then deleting edges.

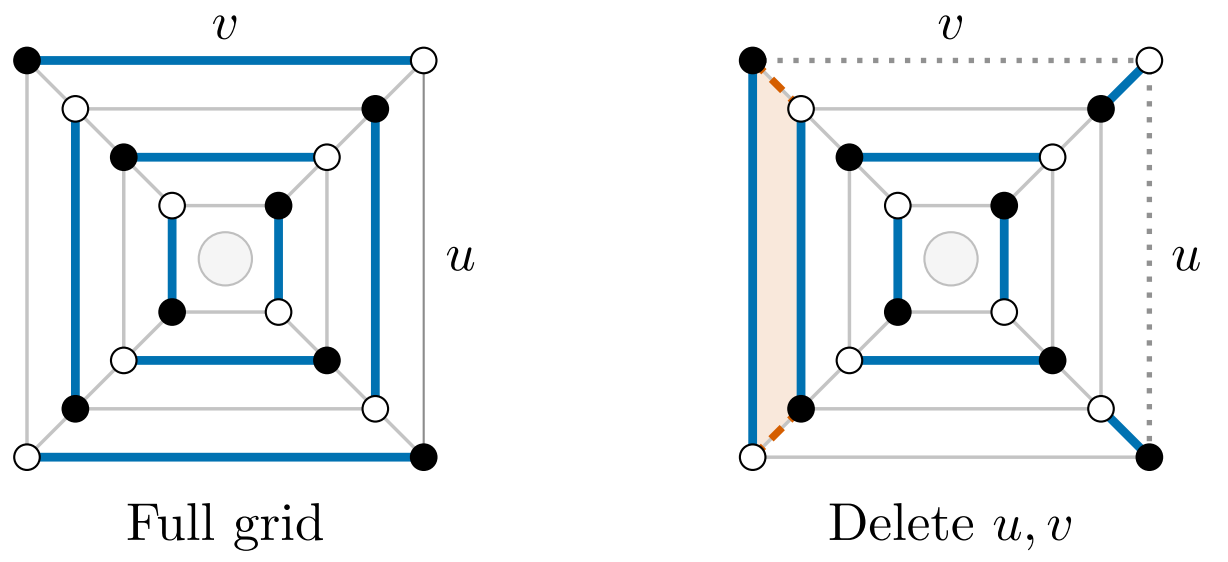
\begin{figure}[htbp]
\centering
\begin{tikzpicture}[x=.80cm,y=.80cm]
\begin{scope}
\fill[black!4] (0,0) circle (.22);\draw[black!25] (0,0) circle (.22);
\coordinate (p00) at (0.44,0.44);
\coordinate (p10) at (-0.44,0.44);
\coordinate (p20) at (-0.44,-0.44);
\coordinate (p30) at (0.44,-0.44);
\coordinate (p01) at (0.84,0.84);
\coordinate (p11) at (-0.84,0.84);
\coordinate (p21) at (-0.84,-0.84);
\coordinate (p31) at (0.84,-0.84);
\coordinate (p02) at (1.24,1.24);
\coordinate (p12) at (-1.24,1.24);
\coordinate (p22) at (-1.24,-1.24);
\coordinate (p32) at (1.24,-1.24);
\coordinate (p03) at (1.64,1.64);
\coordinate (p13) at (-1.64,1.64);
\coordinate (p23) at (-1.64,-1.64);
\coordinate (p33) at (1.64,-1.64);
\draw[base] (p00)--(p10);
\draw[base] (p01)--(p11);
\draw[base] (p02)--(p12);
\draw[black!45] (p03)--(p13);
\draw[base] (p10)--(p20);
\draw[base] (p11)--(p21);
\draw[base] (p12)--(p22);
\draw[base] (p13)--(p23);
\draw[base] (p20)--(p30);
\draw[base] (p21)--(p31);
\draw[base] (p22)--(p32);
\draw[base] (p23)--(p33);
\draw[base] (p30)--(p00);
\draw[base] (p31)--(p01);
\draw[base] (p32)--(p02);
\draw[black!45] (p33)--(p03);
\draw[base] (p00)--(p01);
\draw[base] (p01)--(p02);
\draw[base] (p02)--(p03);
\draw[base] (p10)--(p11);
\draw[base] (p11)--(p12);
\draw[base] (p12)--(p13);
\draw[base] (p20)--(p21);
\draw[base] (p21)--(p22);
\draw[base] (p22)--(p23);
\draw[base] (p30)--(p31);
\draw[base] (p31)--(p32);
\draw[base] (p32)--(p33);
\draw[cA,line width=1.65pt] (p00)--(p30);
\draw[cA,line width=1.65pt] (p02)--(p32);
\draw[cA,line width=1.65pt] (p11)--(p01);
\draw[cA,line width=1.65pt] (p13)--(p03);
\draw[cA,line width=1.65pt] (p20)--(p10);
\draw[cA,line width=1.65pt] (p22)--(p12);
\draw[cA,line width=1.65pt] (p31)--(p21);
\draw[cA,line width=1.65pt] (p33)--(p23);
\node[bv] at (p00) {};
\node[wv] at (p01) {};
\node[bv] at (p02) {};
\node[wv] at (p03) {};
\node[wv] at (p10) {};
\node[bv] at (p11) {};
\node[wv] at (p12) {};
\node[bv] at (p13) {};
\node[bv] at (p20) {};
\node[wv] at (p21) {};
\node[bv] at (p22) {};
\node[wv] at (p23) {};
\node[wv] at (p30) {};
\node[bv] at (p31) {};
\node[wv] at (p32) {};
\node[bv] at (p33) {};
\node[above,font=\small] at (0,1.65) {$v$};
\node[right,font=\small] at (1.65,0) {$u$};
\node[below,font=\small] at (0,-1.85) {Full grid};
\end{scope}
\begin{scope}[xshift=4.8cm]
\fill[black!4] (0,0) circle (.22);\draw[black!25] (0,0) circle (.22);
\coordinate (p00) at (0.44,0.44);
\coordinate (p10) at (-0.44,0.44);
\coordinate (p20) at (-0.44,-0.44);
\coordinate (p30) at (0.44,-0.44);
\coordinate (p01) at (0.84,0.84);
\coordinate (p11) at (-0.84,0.84);
\coordinate (p21) at (-0.84,-0.84);
\coordinate (p31) at (0.84,-0.84);
\coordinate (p02) at (1.24,1.24);
\coordinate (p12) at (-1.24,1.24);
\coordinate (p22) at (-1.24,-1.24);
\coordinate (p32) at (1.24,-1.24);
\coordinate (p03) at (1.64,1.64);
\coordinate (p13) at (-1.64,1.64);
\coordinate (p23) at (-1.64,-1.64);
\coordinate (p33) at (1.64,-1.64);
\fill[cB!14] (p12)--(p22)--(p23)--(p13)--cycle;
\draw[base] (p00)--(p10);
\draw[base] (p01)--(p11);
\draw[base] (p02)--(p12);
\draw[black!43,dotted,line width=1pt] (p03)--(p13);
\draw[base] (p10)--(p20);
\draw[base] (p11)--(p21);
\draw[base] (p12)--(p22);
\draw[base] (p13)--(p23);
\draw[base] (p20)--(p30);
\draw[base] (p21)--(p31);
\draw[base] (p22)--(p32);
\draw[base] (p23)--(p33);
\draw[base] (p30)--(p00);
\draw[base] (p31)--(p01);
\draw[base] (p32)--(p02);
\draw[black!43,dotted,line width=1pt] (p33)--(p03);
\draw[base] (p00)--(p01);
\draw[base] (p01)--(p02);
\draw[base] (p02)--(p03);
\draw[base] (p10)--(p11);
\draw[base] (p11)--(p12);
\draw[base] (p12)--(p13);
\draw[base] (p20)--(p21);
\draw[base] (p21)--(p22);
\draw[base] (p22)--(p23);
\draw[base] (p30)--(p31);
\draw[base] (p31)--(p32);
\draw[base] (p32)--(p33);
\draw[cA,line width=1.65pt] (p00)--(p30);
\draw[cA,line width=1.65pt] (p02)--(p03);
\draw[cA,line width=1.65pt] (p11)--(p01);
\draw[cA,line width=1.65pt] (p13)--(p23);
\draw[cA,line width=1.65pt] (p20)--(p10);
\draw[cA,line width=1.65pt] (p22)--(p12);
\draw[cA,line width=1.65pt] (p31)--(p21);
\draw[cA,line width=1.65pt] (p33)--(p32);
\draw[cB,densely dashed,line width=1.4pt] (p13)--(p12);
\draw[cB,densely dashed,line width=1.4pt] (p22)--(p23);
\node[bv] at (p00) {};
\node[wv] at (p01) {};
\node[bv] at (p02) {};
\node[wv] at (p03) {};
\node[wv] at (p10) {};
\node[bv] at (p11) {};
\node[wv] at (p12) {};
\node[bv] at (p13) {};
\node[bv] at (p20) {};
\node[wv] at (p21) {};
\node[bv] at (p22) {};
\node[wv] at (p23) {};
\node[wv] at (p30) {};
\node[bv] at (p31) {};
\node[wv] at (p32) {};
\node[bv] at (p33) {};
\node[above,font=\small] at (0,1.65) {$v$};
\node[right,font=\small] at (1.65,0) {$u$};
\node[below,font=\small] at (0,-1.85) {Delete $u,v$};
\end{scope}
\end{tikzpicture}
\par\medskip
\renewcommand{\arraystretch}{1.12}
\begin{tabular}{@{}c c c c c c@{}}
\toprule
height & $0$ & $1$ & $2$ & $3$ & $4$\\
\midrule
weight sum & $v$ & $4+2u+17v$ & $25+24u+24v$ & $4+17u+2v$ & $u$\\
after deletion & $0$ & $4$ & $25$ & $4$ & $0$\\
\bottomrule
\end{tabular}
\caption{
Left: the full grid with its unique minimum-height
matching shown in blue. Right: the graph after deleting
$u$ and $v$, with the original height labels retained.
The blue matching now has height $1$. Switching around
the shaded square gives another matching of the same height.
}
\label{fig:minima}
\end{figure}

\begin{example}

Consider the full grid $C_4\times P_4$. Give the two
circumferential edges $u,v$ in Figure~\ref{fig:minima}
variable weights, and give all other edges weight $1$.
The grid has a unique lowest matching, which uses $v$,
and a unique highest matching, which uses $u$.
Setting $u=v=1$, its matching polynomial is
\[
1+23z+73z^2+23z^3+z^4.
\]
Now delete the two edges by setting $u=v=0$.
The polynomial becomes
\[
4z+25z^2+4z^3.
\]
In particular, the remaining graph has four
minimum-height matchings, so the preceding corollary
no longer applies directly.

Nevertheless, positivity survives the deletion.
For example, the original grid has the minor
\[
a_2^2-a_1a_3
=609+1124(u+v)+542(u^2+v^2)+859uv.
\]
After deleting the two edges, this becomes
\[
609=25^2-4^2.
\]

This is the idea behind the general proof. We will
embed an arbitrary dimer graph into a grid with a
unique minimum-height matching, apply the corollary
there, and then delete the extra edges by setting
their weights to zero. The minimum matching need
not survive, but coefficientwise positivity does.
\end{example}

\FloatBarrier

\section{From networks to dimers}\label{sec:annular}\label{sec:superpositions}

\subsection{Grid completion}

We will reduce the general case to a full cylindrical
grid, where the minimum matching is unique.
Indeed, unroll a grid of even circumference so that
positive winding goes to the right. Each black vertex
is matched either west, east, or vertically, with
horizontal displacement $-1$, $1$, or $0$.
The sum of these displacements determines the height
up to a constant. It is uniquely minimized by
matching every black vertex westward.

The next lemma shows that every cylindrical dimer
model can be obtained from such a grid by deleting
edges.

\begin{lemma}[Grid completion]\label{lem:grid-completion}
Let $G$ be a bipartite graph on a cylinder with a
perfect matching. There is a full cylindrical grid
$K=C_{2p}\times P_{2q}$ and a specialization of its
edge weights to $0$, $1$, and the edge variables of
$G$ such that
\begin{equation}\label{eq:grid-completion}
A_K(z)\longmapsto z^q A_G(z).
\end{equation}
In fact, the surviving matchings of $K$ are in
weight-preserving bijection with those of $G$,
with every height increased by $q$.
\end{lemma}

\begin{figure}[htbp]
\centering

% --------------------------------------------------
% Vertex block
% --------------------------------------------------
\begin{minipage}{.31\linewidth}
\centering
\begin{tikzpicture}[x=1.05cm,y=1.0cm,
  path edge/.style={cB,line width=1.3pt},
  marked/.style={cA,line width=1.5pt},
  forced/.style={black!42,line width=1.4pt},
  meshpt/.style={circle,fill=black,inner sep=1.3pt},
  downarrow/.style={->,black!55,line width=.9pt}
]

% Before: mesh point
\begin{scope}[shift={(0,3.05)}]
  \node[meshpt,label={[font=\small]above left:$v$}]
    at (0,.5) {};
  \draw[path edge] (0,.5)--(0,1.15);
  \draw[path edge] (0,.5)--(-.55,.5);
  \draw[marked] (0,.5)--(.55,.5)
    node[right,font=\small] {$m\in M$};
\end{scope}

% Expansion arrow
\draw[downarrow] (0,2.40)--(0,1.95);

% After: 2x2 block
\begin{scope}[shift={(-.5,-.1)}]
  \draw[base,dashed,rounded corners=2pt]
    (-.3,-.3) rectangle (1.3,1.3);
  \draw[base,black!12] (0,0) rectangle (1,1);
  \draw[path edge] (-.66,1)--(0,1)
                    (0,1)--(0,1.66);
  \draw[marked] (0,1)--(1.66,1);
  \draw[forced] (0,0)--(1,0);
  \node[wv] at (0,1) {};
  \node[bv] at (1,1) {};
  \node[bv] at (0,0) {};
  \node[wv] at (1,0) {};
\end{scope}

\node[font=\small] at (0,-.75) {Vertex};

\end{tikzpicture}
\end{minipage}\hfill
%
% --------------------------------------------------
% Path block
% --------------------------------------------------
\begin{minipage}{.31\linewidth}
\centering
\begin{tikzpicture}[x=1.05cm,y=1.0cm,
  path edge/.style={cB,line width=1.3pt},
  forced/.style={black!42,line width=1.4pt},
  meshpt/.style={circle,fill=black,inner sep=1.3pt},
  downarrow/.style={->,black!55,line width=.9pt}
]

% Before: mesh point
\begin{scope}[shift={(0,3.05)}]
  \node[meshpt] at (0,.5) {};
  \draw[path edge] (-.55,.5)--(0,.5)--(0,1.15);
\end{scope}

% Expansion arrow
\draw[downarrow] (0,2.40)--(0,1.95);

% After: 2x2 block
\begin{scope}[shift={(-.5,-.1)}]
  \draw[base,dashed,rounded corners=2pt]
    (-.3,-.3) rectangle (1.3,1.3);
  \draw[base,black!12] (0,0) rectangle (1,1);
  \draw[path edge] (-.66,0)--(0,0)--(0,1.66);
  \draw[forced] (1,0)--(1,1);
  \node[wv] at (0,1) {};
  \node[bv] at (1,1) {};
  \node[bv] at (0,0) {};
  \node[wv] at (1,0) {};
\end{scope}

\node[font=\small] at (0,-.75) {Path block};

\end{tikzpicture}
\end{minipage}\hfill
%
% --------------------------------------------------
% Unused mesh point
% --------------------------------------------------
\begin{minipage}{.31\linewidth}
\centering
\begin{tikzpicture}[x=1.05cm,y=1.0cm,
  forced/.style={black!42,line width=1.4pt},
  meshpt/.style={circle,fill=black,inner sep=1.3pt},
  downarrow/.style={->,black!55,line width=.9pt}
]

% Before: unused mesh point
\begin{scope}[shift={(0,3.05)}]
  \node[meshpt] at (0,.5) {};
\end{scope}

% Expansion arrow
\draw[downarrow] (0,2.40)--(0,1.95);

% After: 2x2 block
\begin{scope}[shift={(-.5,-.1)}]
  \draw[base,dashed,rounded corners=2pt]
    (-.3,-.3) rectangle (1.3,1.3);
  \draw[base,black!12] (0,0) rectangle (1,1);
  \draw[forced] (0,0)--(1,0) (0,1)--(1,1);
  \node[wv] at (0,1) {};
  \node[bv] at (1,1) {};
  \node[bv] at (0,0) {};
  \node[wv] at (1,0) {};
\end{scope}

\node[font=\small] at (0,-.75) {Unused};

\end{tikzpicture}
\end{minipage}

\caption{The three kinds of mesh points and their
expansion into $2\times2$ blocks. At a vertex,
the edge in $M$ takes one step inside the block
before leaving, while the other edges leave directly.
At a point along a path, the path may turn, as shown,
or go straight through. The gray edges match the
remaining vertices.}
\label{fig:grid-blocks}
\end{figure}
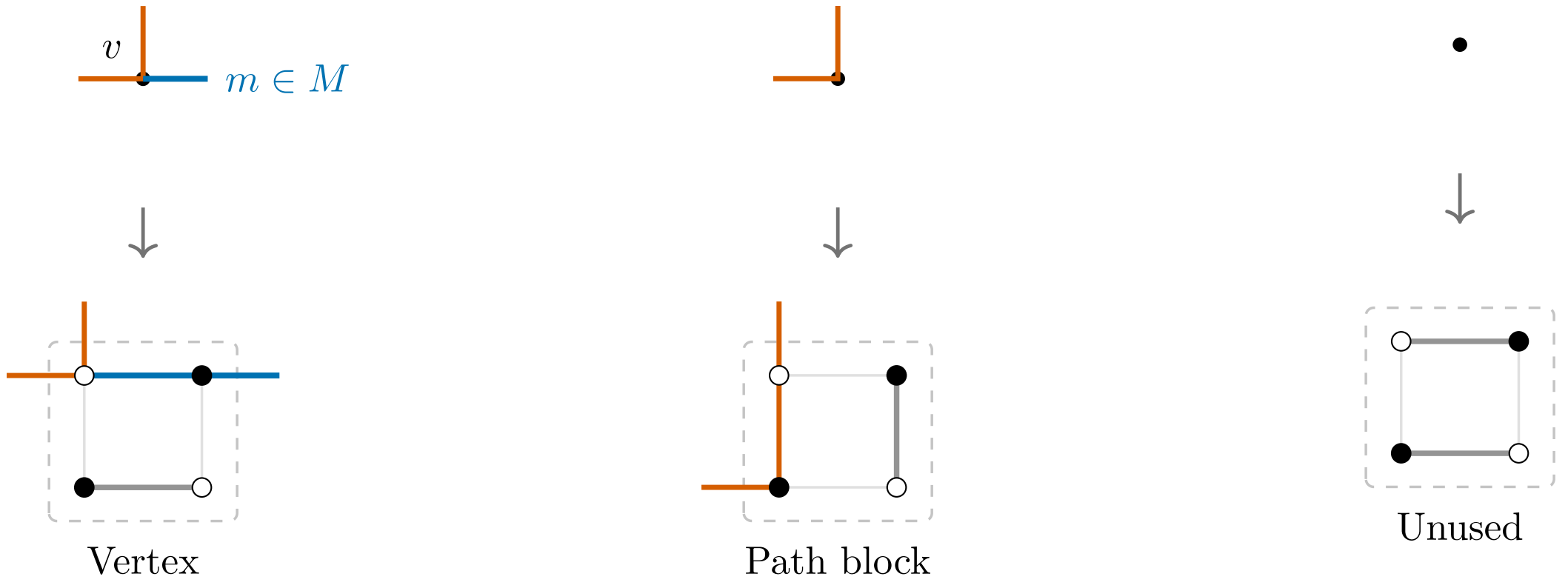

\begin{figure}[!htb]
\centering
\begin{tikzpicture}[
 x=.92cm,y=.92cm,
 me/.style={cA,line width=1.4pt},
 ue/.style={cC,line width=1.4pt},
 ve/.style={cB,line width=1.4pt},
 forced/.style={black!42,line width=1.4pt},
 seam/.style={black!48,densely dashed},
 stagearrow/.style={-{Stealth[length=2.3mm]},black!50,line width=.85pt},
 meshpoint/.style={circle,fill=black!55,inner sep=1.2pt}
]

% 1. Original cylindrical graph (and the two phases of one edge).
\begin{scope}[shift={(.25,11.15)}]
  \node[font=\small] at (2.3,2.07) {The three edges};
  \draw[seam] (0,-.05)--(0,1.77) (4.6,-.05)--(4.6,1.77);
  \draw[me] (1.28,.87)--(3.32,.87)
    node[midway,above,font=\small] {$y_m$};
  \draw[ve] (3.32,.87) to[out=22,in=180] (4.6,1.51);
  \draw[ve] (0,1.51) to[out=0,in=159] (1.28,.87);
  \draw[ue] (3.32,.87) to[out=-22,in=180] (4.6,.24);
  \draw[ue] (0,.24) to[out=0,in=201] (1.28,.87);
  \node[font=\small,above,text=cB] at (.36,1.46) {$v$};
  \node[font=\small,below,text=cC] at (.36,.27) {$u$};
  \node[bv] at (3.32,.87) {};
  \node[wv] at (1.28,.87) {};
\end{scope}
\begin{scope}[shift={(5.75,11.15)}]
  \node[font=\small] at (2.1,2.07) {The two phases of $m$};
  \foreach \y in {.52,1.42}{
    \draw[base] (0,\y)--(3.35,\y);
  }
  \draw[me] (0,1.42)--(.67,1.42)
    (1.34,1.42)--(2.01,1.42)
      node[midway,above=1pt,font=\small] {$y_m$}
    (2.68,1.42)--(3.35,1.42);
  \draw[me] (.67,.52)--(1.34,.52) (2.01,.52)--(2.68,.52);
  \foreach \y in {.52,1.42}{
    \foreach \i in {0,2,4}{
      \node[wv,inner sep=1.45pt] at ({.67*\i},\y) {};
    }
    \foreach \i in {1,3,5}{
      \node[bv,inner sep=1.45pt] at ({.67*\i},\y) {};
    }
  }
  \node[anchor=west,font=\small] at (3.56,1.42) {use $m$};
  \node[anchor=west,font=\small] at (3.56,.52) {omit $m$};
\end{scope}

% Show the passage to mesh routing.
\draw[stagearrow] (4.6,10.91)--(4.6,10.08);

% 2. Smaller intermediate mesh drawing, before expansion to blocks.
% Each dot is a mesh point; only the white and black points represent G.
\begin{scope}[shift={(1.80,6.35)},scale=.82]
  \node[font=\small] at (3.5,3.32) {The three edges routed on the mesh};
  \draw[seam] (-.5,-.22)--(-.5,3.04) (7.5,-.22)--(7.5,3.04);
  % All eight mesh points, including the unused one.
  \foreach \x in {.5,2.5,4.5,6.5}{
    \foreach \y in {.5,2.5}{
      \node[meshpoint] at (\x,\y) {};
    }
  }
  % The three routes; draw edges after the dots and redraw the true vertices.
  \draw[me] (.5,.5)--(4.5,.5);
  \draw[ue] (4.5,.5)--(7.5,.5) (-.5,.5)--(.5,.5);
  \draw[ve] (4.5,.5)--(4.5,2.5)--(7.5,2.5)
    (-.5,2.5)--(.5,2.5)--(.5,.5);
  \node[wv] at (.5,.5) {};
  \node[bv] at (4.5,.5) {};
\end{scope}

% Show the expansion of the mesh points to 2-by-2 square blocks.
\draw[stagearrow] (4.6,5.76)--(4.6,5.05);

% 3. The completed grid. Horizontal and vertical units are equal:
% every 2-by-2 block is drawn as a SQUARE.
\begin{scope}[shift={(1.15,.45)}]
  \node[font=\small] at (3.5,4.18) {After replacing mesh points by $2\times2$ blocks};
  \draw[seam] (-.5,-.5)--(-.5,3.5) (7.5,-.5)--(7.5,3.5);
  % Delineate all eight square blocks lightly.
  \foreach \x in {0,2,4,6}{
    \foreach \y in {0,2}{
      \draw[black!19,dashed,rounded corners=1pt]
        ({\x-.23},{\y-.23}) rectangle ({\x+1.23},{\y+1.23});
    }
  }
  % Paths representing m, u and v.
  \draw[me] (0,1)--(5,1);
  \node[above=1pt,font=\small,text=cA] at (2.5,1) {$y_m$};
  \draw[ue] (5,1)--(7.5,1) (-.5,1)--(0,1);
  \node[font=\small,anchor=east,text=cC] at (-.57,1) {$u$};
  \draw[ve] (5,1)--(5,3)--(7.5,3)
    (-.5,3)--(0,3)--(0,1);
  \node[font=\small,left=3pt,text=cB] at (0,1.5) {$v$};
  % Forced internal edges.
  \draw[forced] (0,0)--(1,0) (2,0)--(3,0)
    (4,0)--(5,0) (6,0)--(7,0)
    (1,2)--(1,3) (2,2)--(3,2)
    (2,3)--(3,3) (4,2)--(4,3) (6,2)--(7,2);
  \foreach \x in {0,2,4,6}{
    \foreach \y in {0,2}{
      \node[bv,inner sep=1.45pt] at (\x,\y) {};
      \node[wv,inner sep=1.45pt] at ({\x+1},\y) {};
      \node[wv,inner sep=1.45pt] at (\x,{\y+1}) {};
      \node[bv,inner sep=1.45pt] at ({\x+1},{\y+1}) {};
    }
  }
\end{scope}
\end{tikzpicture}
\caption{The three edges of the cylindrical graph, their routes on a
rectangular mesh, and the resulting grid completion on
$C_8\times P_4$. The dashed sides are identified. Each mesh
point becomes a square $2\times2$ block. The gray edges are
forced, and $y_m$ is placed on the third edge of the blue path.
The top right shows the two ways to match the path for $m$.}
\label{fig:grid-junction}
\end{figure}
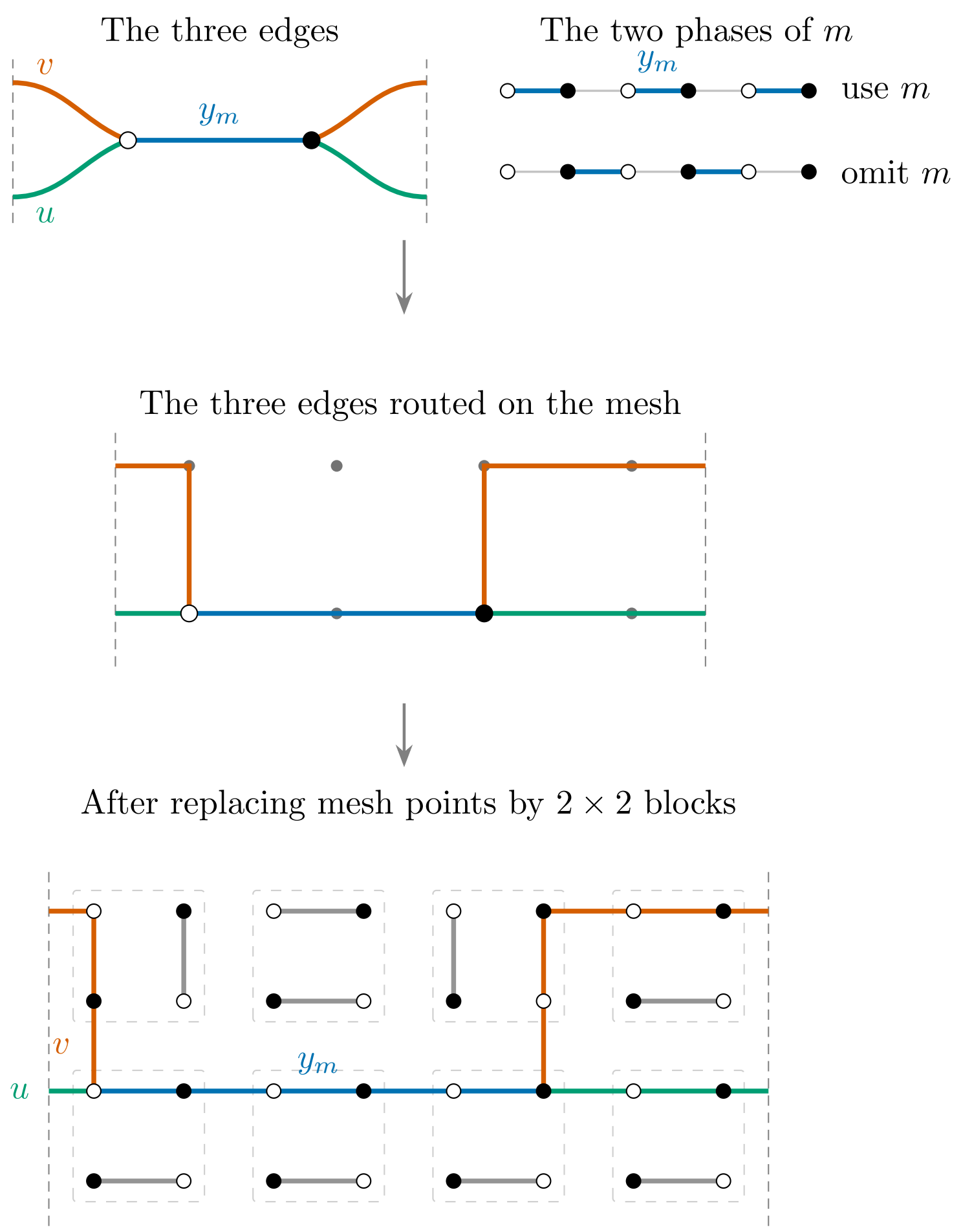

\begin{proof}

Replacing an edge by a path of odd length does not
change the perfect matchings. Indeed, there are exactly
two ways to match the internal vertices: use the odd
edges of the path, covering both endpoints, or use
the even edges, leaving the endpoints free
(Figure~\ref{fig:grid-junction}). Put the original
edge weight on one of the odd edges and give all
other edges weight $1$. We will embed such a subdivision of $G$ into
the grid, then match the remaining vertices
using edges that are forced in every surviving
perfect matching.

First, we make all vertex degrees at most three.
Choose a minimum-height matching $M$ of $G$.
At a vertex $v$ of degree $d\ge4$, replace $v$
by a short path $v'-x-v''$, where $v'$ and $v''$
have the same color as $v$. Attach two consecutive
edges of $v$ to $v'$ and the remaining $d-2$ to
$v''$. Every perfect matching uses exactly one
of the original incident edges, while $x$ is
matched to the other end of the new path.
Thus perfect matchings are preserved. Give the
two new edges weight $1$, and repeat until all
degrees are at most three, carrying $M$ along.

We first draw $G$ on a sufficiently fine periodic
rectangular mesh. At each vertex $v$, choose one
of the two corners of its color in the future
$2\times2$ block, as in Figure~\ref{fig:grid-blocks}.
The edges outside $M$ leave directly through
the two sides meeting at this corner. The edge
of $M$ first takes one step to an adjacent corner,
then leaves through a side not meeting $v$.

We can arrange these exits to respect the cyclic
order of the edges around $v$: the two edges
outside $M$ may be assigned to their sides in
either order. Unused exits are left empty.

By refining the mesh, we can now draw each edge
as a horizontal-vertical path connecting its
chosen exits, with no intersections. This works
across the periodic seam as well, and preserves
the winding numbers of cycles.

Now replace every mesh point by a $2\times2$ block,
as in Figure~\ref{fig:grid-blocks}. At an original vertex $v$, place $v$ at the corner
chosen above, and add the internal step for the
edge of $M$. Match the two
remaining corners by a gray edge of weight $1$.

We still need to connect these vertex blocks along
the paths drawn on the mesh. Between two neighboring
blocks there are two possible crossing edges, one
leaving a black vertex and one leaving a white
vertex. Follow an original edge $e$ from its black
endpoint to its white endpoint. At each crossing,
choose the edge leaving black if $e\notin M$, and
the one leaving white if $e\in M$. In every
intermediate block, the path then enters and
leaves at corners of opposite colors, which
are necessarily adjacent. Join these corners
by an internal edge, and match the other two
corners by a gray edge. At unused mesh points,
put in two gray edges, as in the figure.

Why does this work? The crossings fit at both ends: an edge of $M$ takes an
internal step at each endpoint, and its crossings leave
white, so they meet the corners reached by those steps;
an edge outside $M$ takes no internal step, and its
crossings leave black, so they meet the vertices themselves. Suppose the route of $e$
passes through $k$ intermediate blocks.
The resulting path has $k+1$ crossing edges
and $k$ internal edges, together with two
additional steps at its endpoints if $e\in M$.
Its length is therefore $2k+1$ or $2k+3$,
in either case odd.

Put the original weight of $e$ on one odd-numbered
edge of its path, and give all other edges in
the construction weight $1$. Set every remaining
grid edge weight to $0$. The gray edges are now
forced. Along each path, we can either match its
odd-numbered edges, corresponding to using $e$,
or its even-numbered edges, corresponding to
omitting $e$. Thus the surviving perfect matchings
are in weight-preserving bijection with those
of $G$.

It remains to check the heights. Since the
construction only subdivides edges and splits
vertices locally, corresponding alternating
cycles have the same winding numbers. Thus
relative heights are preserved.

To find the absolute shift, consider the lift
of $M$. It uses only edges inside the blocks:
on paths belonging to $M$ the odd-numbered
edges are internal, and on all other paths
the even-numbered edges are internal.
In particular, it does not cross a seam
chosen between two block columns. On the other
hand, the minimum matching of the full grid,
in which every black vertex is matched westward,
crosses this seam $q$ times westward. Hence
the lift of $M$ has height $q$, and every
surviving matching has its height increased
by $q$, as claimed.
\end{proof}

Cheng, Peng, and Srinivasan gave a height-preserving
reduction to a subgraph of a cylindrical grid
\cite[Lemma~3.1]{CPS}. For a planar version, see
Datta, Kulkarni, and Roy \cite[Lemma~2.9]{DKR}.
Here we complete the subgraph to a full grid by placing
forced edges on the unused vertices. This gives the exact
height shift $q$ when heights are measured from the
minimum matching of the full grid.

\Needspace{10\baselineskip}

\begin{example}[Three edges in eight blocks]\label{ex:grid-junction}
Consider two vertices joined by three edges of weights
$y_m,u,v$. Take $M=\{m\}$ as the height-$0$ matching.
The other two edges each wind once relative to $m$.

Figure~\ref{fig:grid-junction} shows the construction.
The three edges become paths of lengths $5,3,7$.
To obtain a matching, we use the odd-numbered edges
of one path and the even-numbered edges of the other two.
The lift of $M$ stays inside the blocks, while the
matchings using $u$ or $v$ cross the seam once
eastward. Since there are two block rows, the
specialized grid polynomial is
\[
z^2\bigl(y_m+(u+v)z\bigr).
\]
\end{example}

\subsection{Transferring positivity to dimers}

We can now transfer positivity from directed networks
to arbitrary cylindrical dimers. 
For $r\ge1$, expand
\begin{equation}\label{eq:superposition}
\prod_{i=1}^r A_G(x_i)
=\sum_{\mm}\yy^{\mm}F_{\mm}(\xx),
\qquad
\yy^{\mm}=\prod_e y_e^{m_e}.
\end{equation}
Thus $F_{\mm}$ counts ordered $r$-tuples of
matchings in which each edge $e$ is used exactly
$m_e$ times. The exponent of $x_i$ records the
height of the $i$th matching.

The polynomial $F_{\mm}$ is symmetric, since
we can permute the matchings. It is also
homogeneous: the sum of their heights depends
only on the superposition, not on how its edges
are divided among the matchings.

\begin{proposition}[Positivity of superpositions]
\label{prop:schur}
Every superposition polynomial $F_{\mm}$
is Schur positive. Moreover, if
$\mathcal S(\mu)$ holds for all partitions
with $\mu_1\le r$, then every $F_{\mm}$
is elementary positive.
\end{proposition}

\begin{proof}
We first work with a full cylindrical grid $K$.
By Proposition~\ref{prop:unique-minimum},
\[
A_K(z)=a_0\det(\mathrm{Id}+zR),
\]
where $a_0$ is the weight of its minimum
matching and $R$ is a disk path matrix,
after a Laurent monomial change of weights.

The two forms of the Cauchy identity give
\[
\begin{aligned}
\prod_{i=1}^r A_K(x_i)
 &=a_0^r\sum_{\ell(\lambda)\le r}
   s_\lambda(\xx)s_{\lambda'}(R)\\
 &=a_0^r\sum_{\mu_1\le r}
   e_\mu(\xx)m_\mu(R).
\end{aligned}
\]
By Corollary~\ref{cor:pathschur}, the Schur
coefficients are nonnegative Laurent
polynomials in the grid edge weights.
The same holds for the elementary coefficients
under the stated assumption on $\mathcal S(\mu)$.

These coefficients are actually polynomials in
the grid edge weights, not just Laurent polynomials.
Indeed, the product on the left is a polynomial
in those weights, and expanding it in the Schur
or elementary basis cannot introduce denominators.
Thus we can set some edge weights to zero, as in
Lemma~\ref{lem:grid-completion}, without losing
positivity.

The product becomes
\[
e_r(\xx)^q\prod_{i=1}^r A_G(x_i).
\]
In the elementary basis, multiplication by
$e_r^q$ simply increases the exponent of $e_r$.
In the Schur basis it adds a rectangle, since
$e_r^qs_\lambda=s_{\lambda+(q^r)}$.
Removing this common factor preserves positivity
in both bases.

Thus $\prod_i A_G(x_i)$ has the required
coefficientwise positivity. Taking the
coefficient of $\yy^{\mm}$ proves both assertions.
\end{proof}

\begin{proof}[Proof of Theorem~\ref{thm:dimer}]
By Proposition~\ref{prop:schur}, every
$F_{\mm}$ is Schur positive. Thus, for every $r$,
all Schur coefficients of
\[
\prod_{i=1}^r A_G(x_i)
=\sum_{\mm}\yy^{\mm}F_{\mm}(\xx)
\]
are polynomials in the edge weights with
nonnegative coefficients. By
Lemma~\ref{lem:jt}, the same holds for every
Toeplitz minor of $(a_k)$.
\end{proof}

The Schur expansion expresses Toeplitz-minor positivity
separately for each superposition. We can therefore fix
the multiset of edges and study its decompositions into
perfect matchings.

\begin{figure}[htbp]
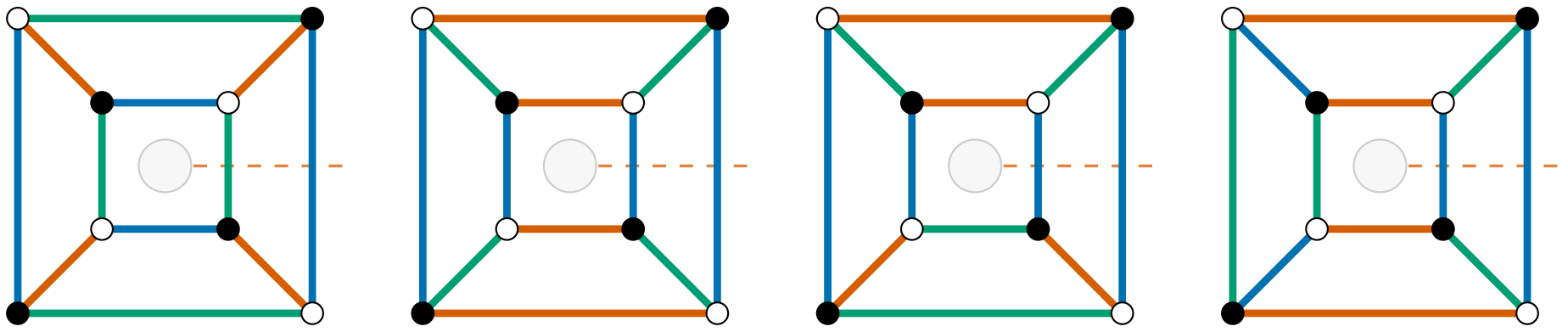

\centering
\ColoringFigure
\caption{The four decompositions of the graph
into three matchings, up to permuting colors.
Each represents six ordered colorings.
The first contributes $m_{21}$, and the
other three together contribute $18m_{111}$.}
\label{fig:colorings}
\end{figure}

\begin{example}
\label{ex:cubic}
Return to the running example, and take $r=3$, with every
edge used exactly once. A decomposition into three
matchings is then a coloring of the edges with
three colors, each color forming a perfect matching.

Up to permuting the colors, there are four such
colorings, shown in Figure~\ref{fig:colorings}.
One has matching heights $0,1,2$, while the other
three have heights $1,1,1$. Each gives six ordered
colorings, so
\begin{equation}\label{eq:cubeSchur}
F_{\mathbf1}
=m_{21}+18m_{111}
=s_{21}+16s_{111}.
\end{equation}
Here we used $s_{21}=m_{21}+2m_{111}$.

Let us see the coefficient $16$ directly.
By \eqref{eq:alternant}, the coefficient of
$\prod_e y_e$ in the cubic Toeplitz minor
\[
\det\begin{pmatrix}
a_1&a_2&0\\
a_0&a_1&a_2\\
0&a_0&a_1
\end{pmatrix}
=a_1^3-2a_0a_1a_2
\]
is $18-2=16$, exactly the coefficient of
$s_{111}$ above. Similarly, the coefficient
of $s_{21}$ is $1$, corresponding to the
minor $a_0a_1a_2$. The upper-right zero is $a_3=0$, and the lower-left zero
is $a_{-1}=0$.
\end{example}

\subsection{From matching tuples to edge-colorings}

There is a direct connection between the superposition
polynomials $F_{\mm}$ and the edge-coloring polynomials
$Z_H$ studied earlier.

\begin{proposition}\label{prop:superposition-coloring}
Fix a nonzero $F_{\mm}$, and replace each edge $e$
of $G$ by $m_e$ distinguished parallel copies.
This gives an $r$-regular bipartite multigraph $H$.
Let $b_{\mm}$ be the minimum height, measured in $G$,
of a perfect matching of $H$. Then
\begin{equation}\label{eq:factorial}
F_{\mm}(\xx)
=\frac{e_r(\xx)^{b_{\mm}}}{\mm!}Z_H(\xx),
\qquad
\mm!=\prod_e m_e!.
\end{equation}
\end{proposition}

\begin{proof}
A proper $r$-edge-coloring of $H$ is a decomposition
into $r$ perfect matchings, one of each color.
Forgetting the distinctions between parallel edges
gives an ordered tuple counted by $F_{\mm}$.

Conversely, each such tuple gives exactly $\mm!$
colorings of $H$: for every edge $e$, we can assign
its $m_e$ occurrences to the parallel copies in
$m_e!$ ways.

The only difference is the normalization of heights.
The $i$th matching has height $b_{\mm}$ more in $G$ than
in $H$, giving a factor of $x_i^{b_{\mm}}$. Multiplying
over all $r$ matchings gives the factor $e_r^{b_{\mm}}$
in \eqref{eq:factorial}.
\end{proof}

In particular, Theorem~\ref{thm:fourcolumns-intro}
implies that $F_{\mm}$ is elementary positive
for $r\le4$. For our ladder example, this gives
\[
F_{\mathbf1}=e_2e_1+15e_3,
\]
strengthening its Schur expansion in
\eqref{eq:cubeSchur}.

\begin{figure}[htbp]
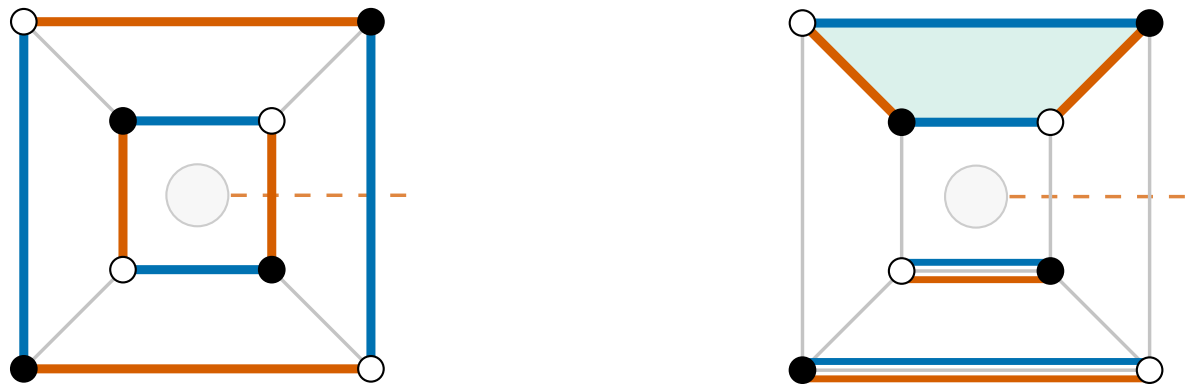

\centering
\TwoSuperpositionsFigure
\caption{Two superpositions. On the left,
two essential alternating cycles give $e_1^2$. On the
right, a contractible alternating cycle gives two
decompositions, both with heights $(1,1)$, giving $2e_2$.}
\label{fig:two-superpositions}
\end{figure}

For two colors we can say more. Remove the doubled
edges from a superposition, and let $c$ and $\ell$
be the numbers of contractible and essential
alternating cycles, respectively. Swapping colors
on a contractible cycle preserves both heights,
giving a factor $2$. On an essential cycle it
transfers one unit of height between the colors,
giving a factor $e_1$. Thus
\[
F_{\mm}(x_1,x_2)=2^c e_2^\beta e_1^\ell,
\qquad \beta\ge0.
\]
This is the inherited-height version of
\eqref{eq:two-color}. Figure~\ref{fig:two-superpositions}
illustrates both kinds of cycles.

\subsection{Elementary positivity and Stembridge's conjecture}
\label{sec:dictionary}

\begin{proof}[Proof of Theorem~\ref{thm:dictionary-intro}]
We first compare superpositions and edge-colorings.
By Proposition~\ref{prop:superposition-coloring},
\[
F_{\mm}=\frac{e_r^{b_{\mm}}}{\mm!}Z_H.
\]
Thus elementary positivity of all $Z_H$ implies
that of all $F_{\mm}$. The coefficients are
nonnegative rational numbers, and in fact integers,
since $F_{\mm}$ has integer coefficients and the
change to the elementary basis is integral.
Conversely, take $G=H$ and $m_e=1$ for every edge.
Then $F_{\mathbf1}=Z_H$. Hence these two
positivity statements are equivalent.

Now suppose $\mathcal S(\mu)$ holds for every
$\mu_1\le r$. The Cauchy expansion and grid
completion, as in the proof of
Proposition~\ref{prop:schur}, show that every
$F_{\mm}$ is elementary positive.

For the converse, suppose every $Z_H$ is
elementary positive. By the gadget lemma,
every interval closure is then elementary
positive in $r$ variables. Intervals of
length greater than $r$ give zero.
Equation~\eqref{eq:sign-twist} now gives
$\mathcal S(\mu)$ for every $\mu_1\le r$.
\end{proof}

Thus elementary positivity of dimer superpositions
is equivalent to Stembridge's conjecture. We can
state it entirely in terms of matchings.

\begin{conjecture}
\label{conj:epositive}
Every dimer superposition polynomial
$F_{\mm}(x_1,\ldots,x_r)$ belongs to
$\NN[e_1,\ldots,e_r]$.
\end{conjecture}

By Theorem~\ref{thm:fourcolumns-intro},
the conjecture holds for $r\le4$.

\subsection{From colorings back to closures}

The gadget lemma turns products of interval
antisymmetrizers into annular edge-coloring problems.
We can also go in the opposite direction.

\begin{proposition}\label{prop:grid-closure}
For every $r$-regular annular bipartite multigraph
$H$, there is a product $Q$ of interval
antisymmetrizers such that
\[
\Phi(Q)(x_1,\ldots,x_r)
=c\,e_r^t Z_H(x_1,\ldots,x_r),
\qquad c>0,\quad t\ge0.
\]
\end{proposition}

\begin{figure}[htbp]
\centering
\begin{tikzpicture}[
  x=.76cm,y=.85cm,
  keep/.style={cA,line width=1.5pt},
  link/.style={cB,arr},
  active/.style={circle,inner sep=1.5pt,fill=black}
]

% Before contraction
\draw[base] (0,0) grid (5,1);

\draw[keep]
  (0,1)--(1,1)
  (2,1)--(3,1)
  (4,1)--(5,1)
  (1,0)--(2,0)
  (3,0)--(4,0);

\foreach \x in {0,2}
  \draw[link] (\x,0)--({\x+1},0);

\foreach \x in {1,3}
  \draw[link] (\x,1)--({\x+1},1);

\foreach \x in {1,3}
  \draw[link] (\x,1)--(\x,0);

\foreach \x in {2,4}
  \draw[link] (\x,0)--(\x,1);

\node[bv] at (0,0) {};
\node[wv] at (1,0) {};
\node[bv] at (2,0) {};
\node[wv] at (3,0) {};
\node[bv] at (4,0) {};
\node[wv] at (5,0) {};

\node[wv] at (0,1) {};
\node[bv] at (1,1) {};
\node[wv] at (2,1) {};
\node[bv] at (3,1) {};
\node[wv] at (4,1) {};
\node[bv] at (5,1) {};

\node[font=\small] at (2.5,-.55)
  {All-west matching};

% Arrow between pictures
\node at (6,.5) {$\longrightarrow$};

% After contraction
\begin{scope}[shift={(7,0)}]

  \foreach \x in {.5,2.5}
    \draw[link] (\x,1)--({\x+2},1);

  \draw[link] (1.5,0)--(3.5,0);

  \foreach \x in {.5,2.5}
    \draw[link] (\x,1)--({\x+1},0);

  \foreach \x in {1.5,3.5}
    \draw[link] (\x,0)--({\x+1},1);

  \foreach \x in {.5,2.5,4.5}
    \node[active] at (\x,1) {};

  \foreach \x in {1.5,3.5}
    \node[active] at (\x,0) {};

  \node[font=\small] at (2.5,-.55)
    {After contraction};

\end{scope}
\end{tikzpicture}

\caption{Contracting the all-west matching
of the full grid. The remaining edges become
directed arrows advancing one or two columns,
so all arrows go forward around the cylinder.}
\label{fig:grid-west}
\end{figure}
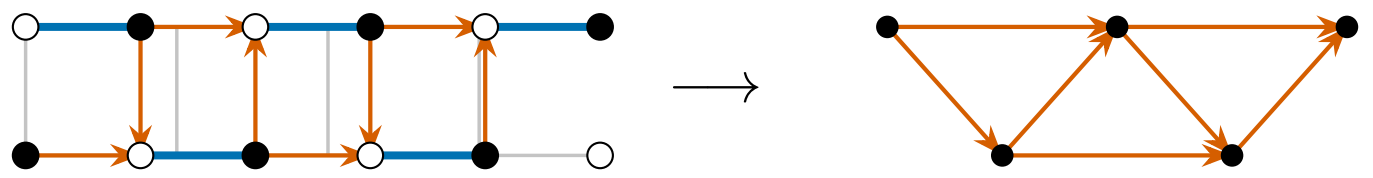

\begin{proof}

Apply grid completion to $H$, but now consider
all $r$ color classes at once. Choose a proper
$r$-edge-coloring of $H$, lift its matchings to
the grid, and record how many lifted matchings
use each grid edge.

Each edge of $H$ starts with multiplicity $1$. At a split
$v'-x-v''$, a color uses $v'-x$ exactly when its matching
uses an edge assigned to $v''$, and conversely for $x-v''$.
Thus each split determines the new multiplicities. A path
replacing an edge of multiplicity $a$ then has alternating
multiplicities $a$ and $r-a$, and each forced edge has
multiplicity $r$.

Thus all colorings of $H$ give decompositions
of the same grid superposition. Conversely,
the matching bijection of
Lemma~\ref{lem:grid-completion} shows that every
decomposition of this superposition comes from
a coloring of $H$.

We now turn this superposition into a product
of interval antisymmetrizers. Contract the
minimum matching of the full grid, as in
Figure~\ref{fig:grid-west}. The remaining edges
form a planar directed network: horizontal
edges advance two columns, and vertical edges
advance one. In particular, every arrow goes
forward around the cylinder. Replace each
directed edge by as many parallel strands
as its multiplicity in the superposition.

Imagine moving around the cylinder, passing
the vertices one at a time. Suppose $k$ strands
enter a vertex. They carry $k$ distinct colors,
and exactly the same colors must leave, possibly
in a different order. The other $r-k$ colors
use the contracted matching edge.

Perturb the vertices in each column slightly so that
the sweep encounters them one at a time. Just before the
sweep reaches a vertex, the incoming strands are consecutive
along the sweep line, since they converge to that vertex
and every other strand stays away from it. At each
vertex, use an interval antisymmetrizer on
the $k$ corresponding positions. Let $Q$ be
the product of these antisymmetrizers, in the
order the vertices are encountered.

By Lemma~\ref{lem:positive-closure}, each allowed
reordering contributes $1/k!$. Conversely,
coloring the strands determines a coloring of
the grid superposition, and hence of $H$.
The only extra choices are the orders of colors
on parallel strands. Their number depends on
the fixed edge multiplicities, not on the
coloring. Thus the two generating functions
agree up to a positive constant $c$.

Finally, the closure records the number of
strands crossing the cut. This is the height
of the lifted matching in the full grid.
By Lemma~\ref{lem:grid-completion}, it is
$q$ more than the original height on $H$.
Hence $t=q$, completing the proof.
\end{proof}

There is also a connection to annular webs.
Drawing each event as a merge followed by a
split gives a positive annular web, whose
coloring function is the specialization at
quantum parameter $1$ considered by
Oh--Yoon~\cite[Section~3]{OY}.
The results of Gorsky--Wedrich
\cite[Corollary~1.4]{GW} provide another
route to Schur positivity. The conjectured
elementary positivity of positive annular
webs~\cite[Conjecture~5.4]{QR},
as discussed in~\cite[Conjecture~1.1]{OY},
would imply the dimer form of Stembridge's
conjecture in every rank.

\FloatBarrier

\end{document}